\documentclass[review]{elsarticle}

\usepackage{lineno}
\modulolinenumbers[1]

\usepackage{amsmath}
\usepackage{dsfont}
\usepackage{algpseudocode}
\usepackage{amssymb}
\usepackage{graphicx,color}
\usepackage{extarrows}
\usepackage{mathrsfs}
\usepackage[ruled,linesnumbered]{algorithm2e}
\usepackage{booktabs}
\usepackage{multirow}
\usepackage{threeparttable}
\usepackage[justification=raggedright]{caption}
\usepackage{threeparttable}
\usepackage[misc]{ifsym}
\usepackage{array}
\usepackage{makecell}
\usepackage{float}
\usepackage{bm}
\usepackage{subcaption}

\makeatletter
\@addtoreset{equation}{section}
\makeatother

\graphicspath{ {./images/} }

\newtheorem{theorem}{Theorem}[section]
\newtheorem{remark}{Remark}
\newtheorem{proof}{Proof}
\newtheorem{lemma}{Lemma}[section]
\newtheorem{definition}{Definition}[section]

\usepackage[utf8]{inputenc}

\usepackage{xcolor,colortbl}
\usepackage{wrapfig}

\definecolor{Gray}{gray}{0.85}
\definecolor{LightCyan}{rgb}{0.88,1,1}
\definecolor{pink}{rgb}{0.6,0.1,0.2}
\newcolumntype{a}{>{\columncolor{Gray}}c}
\newcolumntype{b}{>{\columncolor{white}}c}
\modulolinenumbers[5]
\usepackage[titletoc,page]{appendix}
\usepackage{hyperref}

\begin{document}

\begin{frontmatter}

\title{{ A Priori} Error Estimation for Physics-Informed Neural Networks Solving Allen--Cahn and Cahn--Hilliard Equations}

\author[a,b]{$\text{Guangtao Zhang}$}
\author[a,c]{$\text{Jiani Lin}$}
\author[d]{$\text{Qijia Zhai}$}
\author[a,e]{$\text{Huiyu Yang}$}
\author[a,e]{$\text{Xujun Chen}$}
\author[b]{$\text{Ieng Tak Leong}$\corref{cor1}}
\author[f]{$\text{Fang Zhu}$\corref{cor1}}
\cortext[cor1]{Corresponding authors (emails: $<$itleong@um.edu.mo$>$ (Ieng Tak Leong); $<$zhufang147@gmail.com$>$ (Fang Zhu)). }
\address[a]{SandGold AI Research, Guangzhou, Guangdong, China}
\address[b]{Department of Mathematics, Faculty of Science and  Technology, University of Macau, Macau, China}
\address[c]{School of Mathematical Sciences, Xiamen University, Xiamen, China}
\address[d]{School of Mathematics, Sichuan University, Chengdu, China}
\address[e]{College of Mathematics and Informatics, South China Agricultural University, Guangdong, China}
\address[f]{School of Computer Science and Engineering, Faculty of Innovation Engineering, Macau University of Science and Technology, Macao Special Administrative Region of China}

\begin{abstract}
{Physics-Informed Neural Networks (PINNs) encounter accuracy limitations when solving the Allen--Cahn (AC) and Cahn--Hilliard (CH) partial differential equations (PDEs). To overcome this, we employ a novel loss function, Residuals-weighted Region Activation Evaluation (Residuals-RAE), featuring a { pre-training weight update scheme}. { Unlike conventional self-adaptive PINNs where weights evolve simultaneously with network parameters, Residuals-RAE-PINNs computes weights from current residuals before each training step and holds them constant during gradient updates. We establish weight convergence under standard neural network optimization assumptions, which justifies analyzing the converged network with constant weights.} Based on this theoretical framework, we derive the error estimation for PINNs with Residuals-RAE when solving AC and CH equations. {  The analysis is aligned with Monte-Carlo (Latin hypercube) sampling for the discretization of integrals, consistent with the numerical experiments.} Numerical experiments on one- and two-dimensional AC and CH systems confirm our theoretical results. Additionally, our analysis reveals that feedforward neural networks with two hidden layers and the tanh activation function bound the approximation errors of the solution, its temporal derivative, and the nonlinear term, constrained by the training loss and the number of collocation points.}
\end{abstract}
\vspace{1cm}
\end{frontmatter}

\section{Introduction}
\label{section::introduction}


In recent decades, the deep learning community has made significant progress in various fields of science and engineering, such as computer vision\cite{du2017fused, voulodimos2018deep, fraternali2023black},  material science\cite{choudhary2022recent, wiercioch2023dnn, song2023high}, and finance\cite{heaton2017deep, naveed2023artificial, sahu2023overview}. These achievements highlight the remarkable capabilities of deep neural networks.  One major advantage is their ability to learn complex features and patterns from complex, high-dimensional data.  Deep neural networks (DNN) possess automated representation learning capabilities that obviate tedious modeling and solving engineering. Since DNNs are versatile function approximators, it is natural to leverage them as approximation spaces for solving partial differential equations (PDEs) solutions. In recent years, there has been a significant increase in literature on using deep learning for numerically approximating PDE solutions. Notable examples include the application of deep neural networks to approximate high-dimensional semi-linear parabolic PDEs \cite{han2017deep}, linear elliptic PDEs \cite{schwab2019deep, kutyniok2022theoretical}, and nonlinear hyperbolic PDEs \cite{lye2020deep, lye2021iterative}.

{In recent years, specialized deep neural network architectures have gained traction within the scientific community by enabling physics-informed and reproducible machine learning.  A significant framework is Physics-Informed Neural Network (PINN) \cite{raissi2019physics}, which incorporate physical laws and domain expertise within the loss function.} By incorporating such prioris and constraints, PINN is capable of accurately representing complex dynamics even with limited data availability. PINN has gained substantial popularity in physics-centric fields like fluid dynamics, quantum mechanics, and combustion modeling.  These networks effectively combine data-driven learning with first-principles-based computation \cite{cai2021physics}. Noteworthy applications that showcase the potential of this approach include high-fidelity modeling of turbines, prediction of extreme climate events, and data-efficient control for soft robots. The inherent flexibility of the PINN formulation continues to foster innovative research at the intersection of physical sciences and deep learning \cite{chen2020physics, jagtap2020conservative, jagtap2020extended, yang2021b}.




In the domain of error estimation in Physics-Informed Neural Networks (PINN), the prevailing method relies on statistical principles, specifically the concept of expectation. These techniques enable us to decompose the error in the training process of PINNs into three constituent parts: approximation error, generalization error, and optimization error \cite{Pinkus1999Approximation, Kutyniok2022Mathematics, mishra2022estimates, mishra2023estimates, de2022error, qian2023physics}.

The approximation error quantifies the discrepancy between the approximate and exact solutions, and it reflects the network's capacity to accurately approximate the solution. Since the approximation error is directly related to the network architecture, a significant amount of research in PINN has been dedicated to understanding and characterizing this error. Initial analytical guarantees were provided by Shin et al. \cite{Shin2020convergence}, demonstrating that two-layer neural networks with tanh activations can universally approximate solutions in designated Sobolev spaces for a model problem. Comparable estimates have also been derived for advection equations, elucidating the approximation accuracy \cite{shin2020error}. More recently, stability estimates for the governing partial differential equations have been employed to formulate bounds on the total error, specifically the generalization error \cite{mishra2022estimates, mishra2023estimates}. Additionally, {  Monte-Carlo (LHS) sampling has been leveraged to estimate the generalization error based on the training error, with the analysis aligned to stochastic sampling and its $O(N^{-1/2})$ convergence} \cite{mishra2023estimates}.

The generalization error quantifies the disparity between the predictions made by a model on sampled data and the true solutions. In the context of PINN (Physics-Informed Neural Networks), the loss function represents an empirical risk, while the actual error is assessed using the { $L^2$} norm of the difference between the network's approximation and the analytical solution. The discretization inherent in PINN's physical loss, which takes integral form, contributes to the generalization error. Ongoing research on bounding the generalization error suggests that it is influenced by the number of collocation points, training loss, and the dimensionality of PDEs. Shin et al. \cite{Shin2020convergence} proposed a framework for estimating the error in linear PDEs, deriving bounds on the generalization error based on the training loss and number of training points. Ryck et al. \cite{Ryck2022Estimates} demonstrated that for problems governed by Navier-Stokes equations, the generalization error can be controlled by the training loss and the curse of the spatial domain's dimensionality. Specifically, reducing errors requires an exponential increase in collocation points and network parameters. Similar results were shown by Mishra et al. \cite{Mishra2021radiative} for radiative transfer equations, indicating that the generalization error primarily depends on the training error and the number of collocation points, with less impact from the dimensionality. 


Furthermore, the errors in PINN are influenced by various factors, including optimization algorithms, hyperparameter tuning, and other related aspects \cite{yuan2022pinn, pantidis2023error, qian2023physics, wang2024pinn}. For instance, when using PINN to solve the Allen--Cahn (AC) and Cahn--Hilliard (CH) equations, there exist few papers addressing these challenges, such as \cite{katbar2025data1} and \cite{katbar2025data}; however, these works do not provide error estimation for PINN in solving AC and CH equations. The optimization error arises from the algorithms utilized to update the network weights. Initially, PINN predominantly employs the Adam optimizer \cite{kingma2014adam} combined with the limited-memory Broyden–Fletcher–Goldfarb–Shanno approach \cite{zhu1997algorithm, berahas2016multi} to optimize the loss function and its gradients. However, it is important to note that the PINN loss landscape is generally non-convex. Therefore, even if the function space contains the analytical solution, optimization constraints only yield a local minimum rather than the global minimum. The gap between the theoretical and achieved optima represents the optimization error, which cannot be completely eliminated due to inherent limitations in the algorithm. While techniques like multi-task learning, uncertainty-weighted losses, and gradient surgery have been proposed in the context of PDE solutions \cite{thanasutives2021adversarial}, more refined methods are needed to balance the competing terms in the PINN loss and their gradients \cite{bischof2021multi}. In summary, the errors in PINN arise from multiple sources. The network architecture introduces approximation errors, while optimization and parameter choices lead to additional errors. Given that PINN applications often involve complex phenomena that demand high accuracy, it is crucial to conduct a thorough analysis of these errors. Understanding the interplay between approximation, generalization, and optimization errors is essential for characterizing PINN convergence and ensuring reliability, thus representing a significant challenge that still needs to be addressed.

{In this work, we explore the use of Physics-Informed Neural Networks (PINNs) to address the Allen--Cahn (AC) and Cahn--Hilliard (CH) equations, which are pivotal in material science for phase field modeling, particularly in understanding diffusion separation and multi-phase flows. These equations present significant challenges due to their high nonlinearity, with the CH equation further complicated by fourth-order derivatives. To improve the accuracy and robustness of PINNs in solving these partial differential equations (PDEs), we introduce a novel loss function termed Residuals-weighted Region Activation Evaluation (Residuals-RAE) \cite{zhangresiduals}. { Unlike conventional self-adaptive PINNs (e.g., SA-PINNs \cite{mcclenny2023self}) where weights are trainable parameters updated simultaneously with network parameters, Residuals-RAE-PINNs employs a pre-training weight update mechanism: at each training iteration, the weights are first computed from current residuals and then held constant during the gradient descent step. This decoupled approach} averages residuals over neighboring points to create a weight distribution that aligns closely with the residual distribution, thus ensuring balanced training and capturing the solution's behavior effectively while avoiding excessive weight fluctuations. { We provide a rigorous theoretical framework establishing that the weights converge to a stable limit as training progresses, under standard assumptions on neural network optimization. This convergence result, together with the pre-training weight update mechanism, justifies analyzing the converged network using the framework of weighted PINNs with constant weights.} A key finding of our study is that the error estimation for PINNs employing Residuals-RAE remains comparable to that of standard PINNs under specific conditions. For the AC equation, we establish bounds on the total error and approximation error between the PINN solution and the exact solution, as detailed in Theorems \ref{theorem::ac4} and \ref{theorem::ac5}, assuming Lipschitz continuity and smoothness. Similarly, for the CH equation, we derive analogous bounds in Theorems \ref{theorem::ch7} and \ref{theorem::CH2_TER_Numerical} under sufficient smoothness assumptions. Our accuracy analysis indicates that for the Allen--Cahn (AC) equation, the method's computational accuracy reaches $\ln^2 (N) N^{-k+2}$ as $N$ increases. For the Cahn--Hilliard (CH) equation, the accuracy is $\ln^2 (N) N^{-k+4}$. In these formulas, $k$ represents the order of continuity of the exact solution $u$. $N$ is a theoretical parameter from function approximation theory that controls the network's architectural complexity; a larger $N$ implies the theoretical existence of a more accurate and complex network. Numerical experiments further corroborate these theoretical insights, highlighting the efficiency and precision of the Residuals-RAE approach.}

This paper is organized as follows. In Section \ref{section::MP}, we provide a concise overview of neural networks used to approximate the solutions of PDEs, as well as { the Residuals-RAE weighting scheme and its theoretical convergence analysis}. Following that, we present an error analysis for both the AC equation and CH equation. Additionally, we carry out a series of numerical experiments on these two PDEs to complement and reinforce our theoretical analyses in Sections \ref{section::3} and \ref{section::4}. In Section \ref{section::5}, we offer some concluding remarks. Finally, we review some auxiliary results for our analysis and provide the proofs of the main Theorems in Sections \ref{section::3} and \ref{section::4} in the appendix.

\section{Preliminary Study on the Application of Deep Learning for Solving Partial Differential Equations (PDEs)} 
\label{section::MP}
\subsection{\textbf{Basic setup and notation of Generic PDEs}}
Let $D\subset \mathbb{R}^d$ be an open  bounded domain in $\mathbb{R}^d$ with a sufficiently smooth boundary $\partial D,$ and $\mathcal{T}=[0,T]$ be a time interval with  ${T}<+\infty.$ 
Let $\Omega =D\times \mathcal{T},$ consider the time-dependent partial differential equation with certain boundary conditions:
\begin{equation}
\left\{\begin{aligned}
   u_t =\mathcal{S}(u) \quad & \text{in} \quad  D \times [0,T],  \\
    u = u_0 \quad &  \text{on} \quad D \times \{0\},\\
   \mathcal{B}(u) = 0 \quad & \text{on} \quad \partial D \times [0,T]
   \end{aligned}\right.
   \label{eq::PDE}
\end{equation}
where $\mathcal{S}(u)=\mathcal{L} u+\mathcal{N}(u)$ is an operator composed of constant-coefficient linear differential operator $\mathcal{L}$, and constant-coefficient nonlinear differential (or non-differential) operator $\mathcal{N}$. Here, $\mathcal{B}(u)$  in \eqref{eq::PDE} denotes the boundary operator on $u$ in  either one of the following 3 types:  Dirichlet, Neumann, or periodic boundary conditions which will be discussed in detail in the following. Our goal is to construct a numerical solution of  PDE problem \eqref{eq::PDE} with the help of neural networks, that is, to find $u(\boldsymbol{x},t)$ in some proper function space $F$ of real-valued functions defined on $D\times \mathcal{T}$  which  satisfy the PDEs above in certain sense.  The Sobolev spaces $H^m(\Omega)$ or $H^m(D)$  will be  the natural choices of our PDE problems, and one can refer  the details in the Appendices \ref{Notation} and \ref{Auxiliary results}.
\subsection{\textbf{Deep neural networks architecture}}
Deep Neural Networks (DNNs) are a specific type of artificial neural network that consists of multiple hidden layers. They are  extensions of the perception-based neural networks. Unlike simpler neural networks, DNNs use the function composition to construct the function approximation $u_{dnn}(\boldsymbol{x}, t; \theta)$ instead of direct methods.

\begin{definition}
\label{def::net}
Let $R \in(0, \infty]$, $L, W \in \mathbb{N}$, $l_0, \cdots, l_L \in \mathbb{N}$. Given an input $z = (\boldsymbol{x}, t)$ and a twice differentiable activation function $\sigma: \mathbb{R} \rightarrow \mathbb{R}
$ with non-linearity {(specifically, we take $\sigma = \tanh$ for concreteness).} We call the affine mapping $u_{dnn}$ the network 
\begin{align}
u_{dnn}(z ; \theta)=f_L^{\theta_L} \circ f_{L-1}^{\theta_{L-1}} \circ \cdots \circ f_1^{\theta_1}(z)=f_L\left(f_{L-1}\left(\cdots f_1\left(z ; \theta_1\right)  ; \theta_{L-1}\right) ; \theta_L\right)
\end{align}
where
\begin{align}
f_k^\theta(z)= \begin{cases}\mathcal{A}_L^\theta(z) & k=L \\ \left(\sigma \circ \mathcal{A}_k^\theta\right)(z) & 1 \leq k<L\end{cases}
\end{align}
which consists of multiple layers of interconnected nodes, including an input layer $f_1$, an output layer $f_L$, and several hidden layers $\{f_2,..., f_{L-1}\}$ in the middle. The matrix multiplication $\mathcal{A}_k^\theta(z): \mathbb{R}^{l_{k-1}} \rightarrow \mathbb{R}^{l_k}$ in the $k^{th}$ layer parameterized by trainable parameters $\theta_{k} := \{W_k, b_k\}$, includes the weights matrix $W_k:= \{w_{i,j}\}_{l_{k} \times l_{k+1}}$, and biases vector $b_k:= \{b_{i,j}\}_{l_{k+1} \times 1}$ with real  entries in $[-R,R].$ Here, $l_0 = d+1$ for approximating the PDE problem (\ref{eq::PDE}). 
\end{definition}

\begin{figure}
\centering
\includegraphics[width=1.0\textwidth]{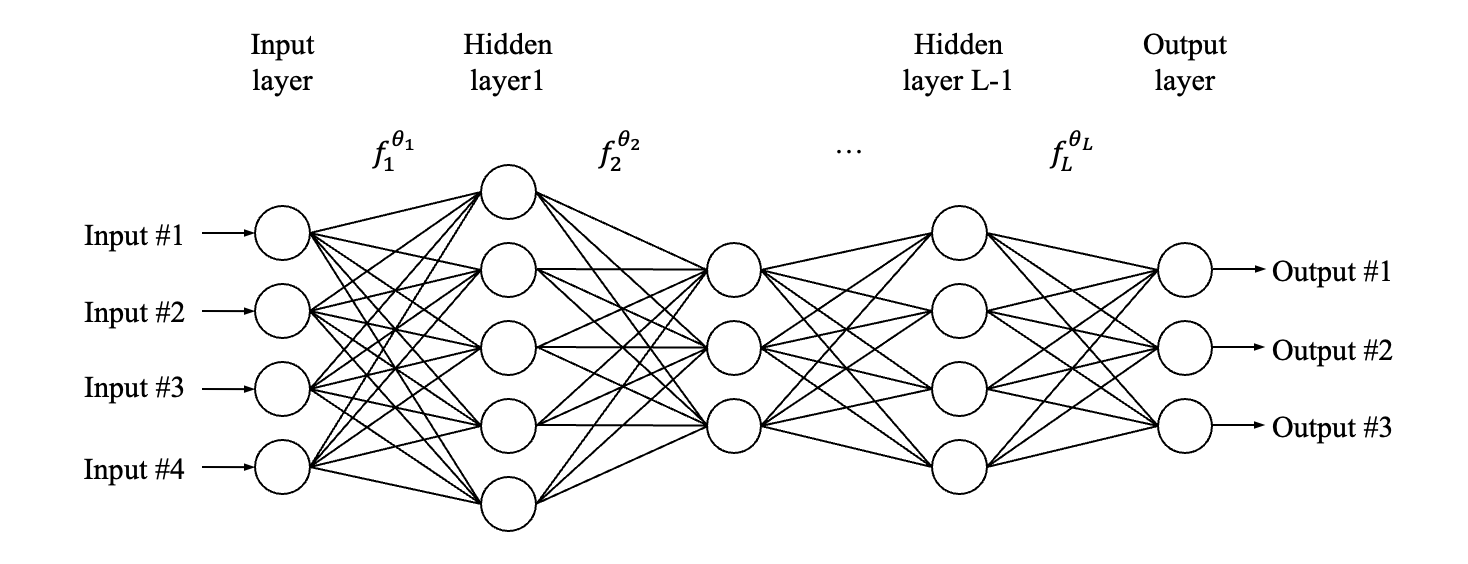}
\caption{\textbf{Structure of fully connected deep neural network.}}
\label{fig::ch2_DNN}
\end{figure}

These layers are introduced into the model using fully connected linear relations and each neuron in a layer is connected to every neuron in the preceding and following layers. This architecture allows the network to capture complex patterns in the data through the combination of these linear relationships and non-linear activation. Over the years, There have been many influential studies on DNN over the years, demonstrating that it can handle the multi-type of data and offer automation and high performance. Here, we can easily find that we could compute the derivatives of output $u$ at the input, such as $\frac{\partial u}{\partial t}$, $\frac{\partial u}{\partial x}$, $\ldots$, by back-propagation method (see Fig. \ref{fig::ch2_DNN}).


Having introduced the neural network architecture, we now establish its approximation properties in the following lemma.
\begin{lemma}
\label{lemma::Error_N}
(\cite{Ryck2022Estimates}) {Let $d \geq 2, n \geq 2, m \geq 3, \delta>0, a_i, b_i \in \mathbb{Z}$ }with $a_i<b_i$ for $1 \leq i \leq d, D=\prod_{i=1}^d\left[a_i, b_i\right]$ and $f \in H^m(D)$. Then for every $N \in \mathbb{N}$ with $N>5$ there exists a tanh neural network $u_{dnn}$ with two hidden layers, one of width at most $3\left\lceil\frac{m+n-2}{2}\right\rceil\left|P_{m-1, d+1}\right|+\sum_{i=1}^d\left(b_i-a_i\right)(N-1)$ and another of width at most $3\left\lceil\frac{d+n}{2}\right\rceil\left|P_{d+1, d+1}\right| N^d \prod_{i=1}^d\left(b_i-a_i\right)$, such that for $k=0,1,2$  the following inequalities hold
$$
\left\|f-u_{dnn}\right\|_{H^k(D)} \leq 2^k 3^d C_{k, m, d, f}(1+\delta) \ln ^k\left(\beta_{k, \delta, d, f} N^{d+m+2}\right) N^{-m+k},
$$
and where
$$
\begin{aligned}
C_{k, m, d, f} & =\max _{0 \leq l \leq k}\left(\begin{array}{c}
d+l-1 \\
l
\end{array}\right)^{1 / 2} \frac{((m-l) !)^{1 / 2}}{\left(\left\lceil\frac{m-l}{d}\right\rceil !\right)^{d / 2}}\left(\frac{3 \sqrt{d}}{\pi}\right)^{m-l}|f|_{H^m(\Omega)}, \\
\beta_{k, \delta, d, f} & =\frac{5 \cdot 2^{k d} \max \left\{\prod_{i=1}^d\left(b_i-a_i\right), d\right\} \max \left\{\|f\|_{W^{k, \infty}(\Omega)}, 1\right\}}{3^d \delta \min \left\{1, C_{k, m, d, f}\right\}} .
\end{aligned}
$$

Moreover, the weights of $u_{dnn}$ scale as $O\left(N^\gamma+N \ln N\right)$ with $\gamma=\max \left\{m^2 / n, d(2+m+d) / n\right\}$.

\end{lemma}


Lemma \ref{lemma::Error_N}, which has been established in \cite{Ryck2022Estimates}, shows that a two-hidden-layer tanh neural network is capable of achieving arbitrary smallness of $\|f-u_{dnn}\|_{H^2(D \times[0,T])}$ (i.e. $\left\|f-u_{dnn}\right\|_{H^m(D \times[0,T])}<\epsilon$). Furthermore, it provides explicit bounds on the required network width. This illustrates the network's strong ability to represent complex functions. In the subsequent section, we can utilize this approximation lemma to prove the upper bound on the residuals of the PDEs derived from the Physics-Informed Neural Network (PINN) (for more details, {refer to Theorem \ref{theorem::ac3} and Theorem \ref{theorem::ch6}).}


\subsection{\textbf{Physics-informed neural networks for approximating the PDEs solution}}



Let $u_{dnn}(\boldsymbol{x},t; \theta)$ represent the approximated solution of the partial differential equation (PDE) (\ref{eq::PDE}) obtained from a neural network, where the network's parameters $\theta=\{W^{L}, b^{L}\}$ correspond to the neurons. These parameters can be updated using various optimization algorithms. In the field of physics-informed neural networks (PINNs) \cite{raissi2019physics}, these networks are used to solve PDEs by incorporating the governing equations as soft constraints in the optimization objective. This objective, denoted as $\mathcal{J}(u(\cdot ; \theta))$, is defined as follows:
\begin{align}
\label{eq::base_loss}
\mathcal{J}(u(\cdot ; \theta)) = \gamma_{int} { \|R_{int}(\boldsymbol{x}_{int},t_{int} ; \theta)\|_{L^2(D \times [0,T])}^2} + \gamma_{tb} { \|R_{tb}(\boldsymbol{x}_{tb} ; \theta)\|_{L^2(D)}^2}+\gamma_{sb} { \|R_{sb}(\boldsymbol{x}_{sb},t_{sb} ; \theta)\|_{L^2(\partial D \times [0,T])}^2},
\end{align}
where $R_{int}(\boldsymbol{x}_{int},t_{int} ; \theta) = u_t(\boldsymbol{x}_{int},t_{int} ; \theta) - \mathcal{S}(u(\boldsymbol{x}_{int},t_{int} ; \theta))$, $R_{sb}(\boldsymbol{x}_{sb},t_{sb} ; \theta) = \mathcal{B}^k(u(\boldsymbol{x}_{sb},t_{sb} ; \theta))-0$, and $R_{tb}(\boldsymbol{x}_{tb} ; \theta) = u(\boldsymbol{x}_{tb} ; \theta)-u_0(\boldsymbol{x}_{tb})$ are the residual terms in the physics-informed neural network (PINN) framework. The hyperparameters, represented by   $\gamma = \{\gamma_{int}, \gamma_{sb}, \gamma_{tb}\}$, serve as penalty coefficients that help balance the learning rate of each individual loss term. {  In this context, we use the standard Sobolev (here $L^2$) norm notation consistent with the theoretical error analysis: $\|\cdot\|_{L^2(D\times [0,T])}^2 := \int_{0}^{T}\int_D (\cdot)^2 d\boldsymbol{x}dt$ and $\|\cdot\|_{L^2(D)}^2 := \int_D (\cdot)^2 d\boldsymbol{x}$, and $\|\cdot\|_{L^2(\partial D \times [0,T])}^2 := \int_0^T\int_{\partial D} (\cdot)^2 ds(\boldsymbol{x})dt$ denote the squared $L^2$ norms over the spatiotemporal interior, the spatial domain $D$, and the spatial boundary in time, respectively.} In construct a reasonable numerical solution of  \eqref{eq::PDE} represented by the neural network, we consider the least squares problems on the loss function $\mathcal{J}(u(\cdot ; \theta))$ with network parameter $\theta,$ this introduces a minimized optimization goal defined as
\begin{align}
\theta^*=\arg \min _{\theta} \mathcal{J}(u(\cdot ; \theta)),
\end{align}
where $u\left(\cdot ; \theta^*\right)$ can be understood as the solution to the PDEs with optimized parameters. In practice, it is necessary to discretize the integral norms in the objective function (\ref{eq::base_loss}) using {  Monte-Carlo (LHS) sampling} (or other appropriate sampling methods), leading to an empirical approximation of the training loss as follows
\begin{subequations}
\begin{align}
{ \|R_{int}(\boldsymbol{x}_{int},t_{int} ; \theta)\|_{L^2(D \times [0,T])}^2} & \approx \frac{1}{M_{int}} \sum_{i=1}^{M_{int}} R_{int}^2(\boldsymbol{x}_{int}^{(i)}, t_{int}^{(i)}; \theta), \\
{ \|R_{sb}(\boldsymbol{x}_{sb},t_{sb} ; \theta)\|_{L^2(\partial D \times [0,T])}^2} & \approx \frac{1}{M_b} \sum_{i=1}^{M_{sb}} R_{sb}^2(\boldsymbol{x}_{sb}^{(i)}, t_{sb}^{(i)}; \theta), \\
{ \|R_{tb}(\boldsymbol{x}_{tb} ; \theta)\|_{L^2(D)}^2} & \approx \frac{1}{M_{tb}} \sum_{i=1}^{M_{tb}} R_{tb}^2(\boldsymbol{x}_{tb}^{(i)}; \theta).
\end{align}
\end{subequations}

\begin{figure}
\centering
\includegraphics[width=1.0\linewidth, height=0.40\linewidth]{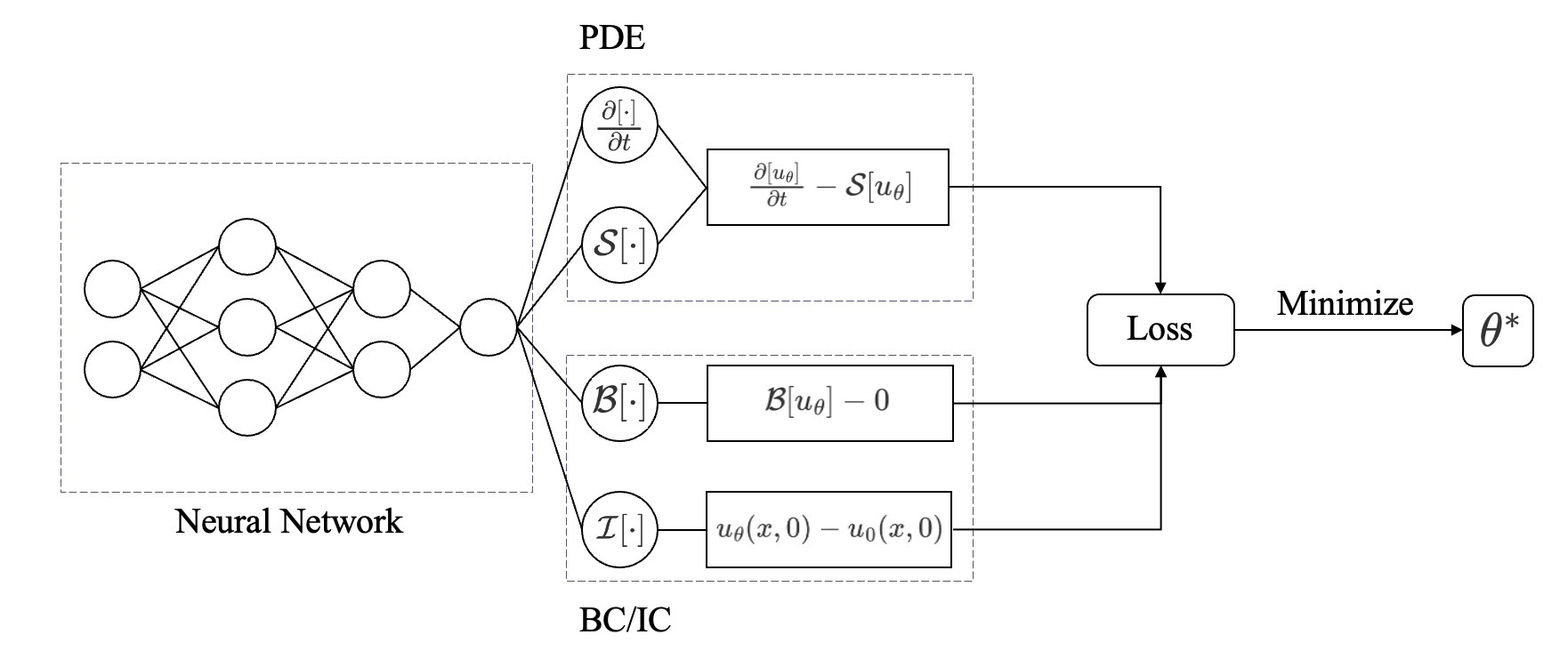}
\caption{\textbf{Schematic diagram of Physics-informed neural networks (PINNs).}}
\label{fig::ch2_uniformsample}
\end{figure}

Here, the training points $\{(\boldsymbol{x}_{int}^{(i)}, t_{int}^{(i)})\}_{i=1}^{M_{int}}$, $\{(\boldsymbol{x}_{sb}^{(i)}, t_{sb}^{(i)})\}_{i=1}^{M_{sb}}$ and $\{(\boldsymbol{x}_{tb}^{(i)})\}_{i=1}^{M_{tb}}$ are located in the interior, boundary and initial conditions of the computational domain, respectively. These points could be sampled using various methods such as uniform or Latin hypercube sampling (LHS). Afterward, the accuracy of the neural network's predictions can be evaluated by comparing them to the exact solution. Let $u_{exact}$, $u_{dnn}$ represent the exact and predicted solutions from PINNs, respectively. In this work, we measure the accuracy using the relative $L^2$ error and the relative $L^{\infty}$ error. {  We define these errors in the continuous (integral) form, consistent with the $L^2$ norm notation used in the theoretical error analysis:}
\begin{equation}
\begin{split}
{ \text{Relative } L^2 \text{ error}} &:= \frac{\|u_{dnn}-u_{exact}\|_{L^2(D\times[0,T])}}{\|u_{exact}\|_{L^2(D\times[0,T])}} = \frac{\left(\int_0^T\int_D \left|u_{dnn}(\boldsymbol{x},t)-u_{exact}(\boldsymbol{x},t)\right|^2 \mathrm{d}\boldsymbol{x}\,\mathrm{d}t\right)^{1/2}}{\left(\int_0^T\int_D \left|u_{exact}(\boldsymbol{x},t)\right|^2 \mathrm{d}\boldsymbol{x}\,\mathrm{d}t\right)^{1/2}},\\
{ \text{Relative } L^{\infty} \text{ error}} &:= \frac{\|u_{dnn}-u_{exact}\|_{L^\infty(D\times[0,T])}}{\|u_{exact}\|_{L^\infty(D\times[0,T])}} = \frac{\operatorname*{ess\,sup}_{(\boldsymbol{x},t)\in D\times[0,T]}\left|u_{dnn}(\boldsymbol{x},t)-u_{exact}(\boldsymbol{x},t)\right|}{\operatorname*{ess\,sup}_{(\boldsymbol{x},t)\in D\times[0,T]}\left|u_{exact}(\boldsymbol{x},t)\right|}.
\end{split}
\end{equation}
In practice, these integrals and essential suprema are approximated using a set of testing points $\left\{(\boldsymbol{x}^{(i)}, t^{(i)})\right\}_{i=1}^{M_{total}} \subset D \times [0, T]$ (e.g., by Monte-Carlo or uniform sampling). In order to compare our numerical solutions with the  another one from the other methods,  our reference solutions for the Allen--Cahn and Cahn--Hilliard equations are the ones obtained by  a Chebyshev polynomial-based numerical algorithm \cite{driscoll2014chebfun}, which ensures highly accurate solutions.

\subsection{\textbf{Self-adaptive learning using the Residuals-RAE weighting scheme}}
\begin{figure}
\centering
\subcaptionbox{Results for 1D Allen--Cahn equation (with initial condition $\# 1$)  using vanilla PINNs.\label{fig::1dAC_vanilla}}{\includegraphics[width=1.0\linewidth]{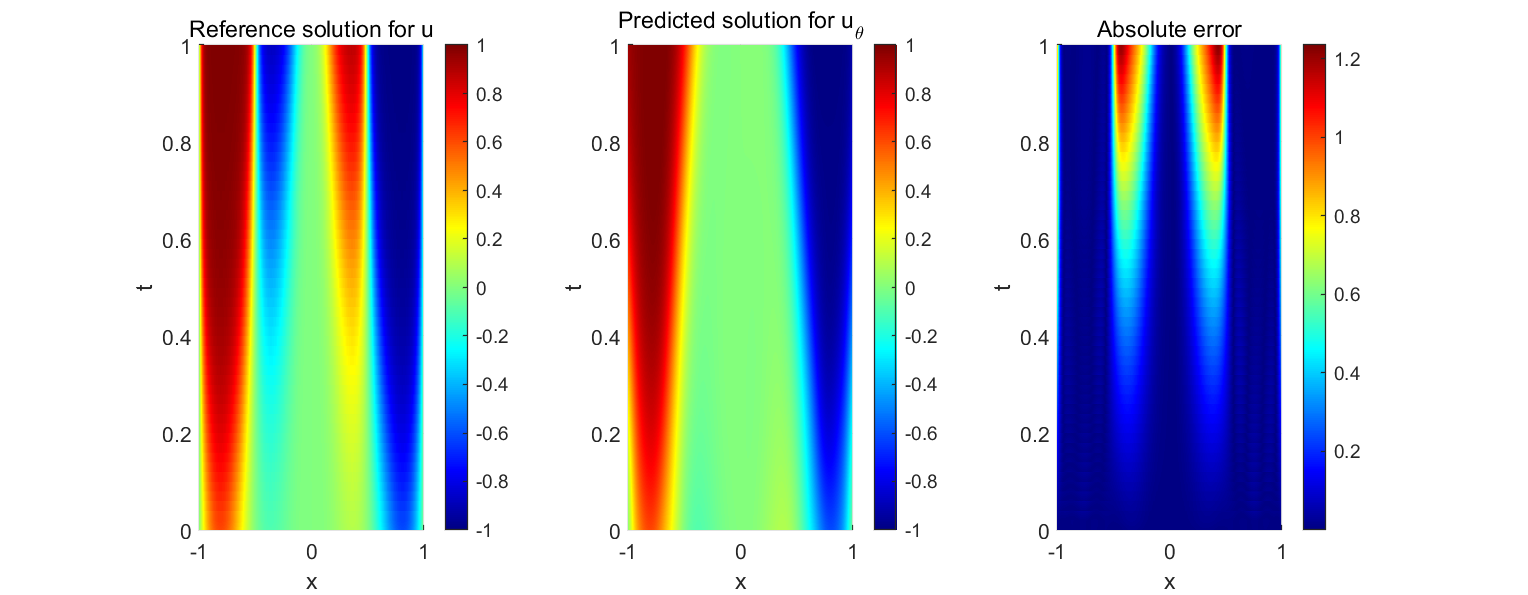}}
\centering
\subcaptionbox{Results for 1D Cahn--Hilliard equation using vanilla PINNs.\label{fig::1dCH_vanilla} \label{fig::1dCH_rae_i2}}{\includegraphics[width=1.0\linewidth]{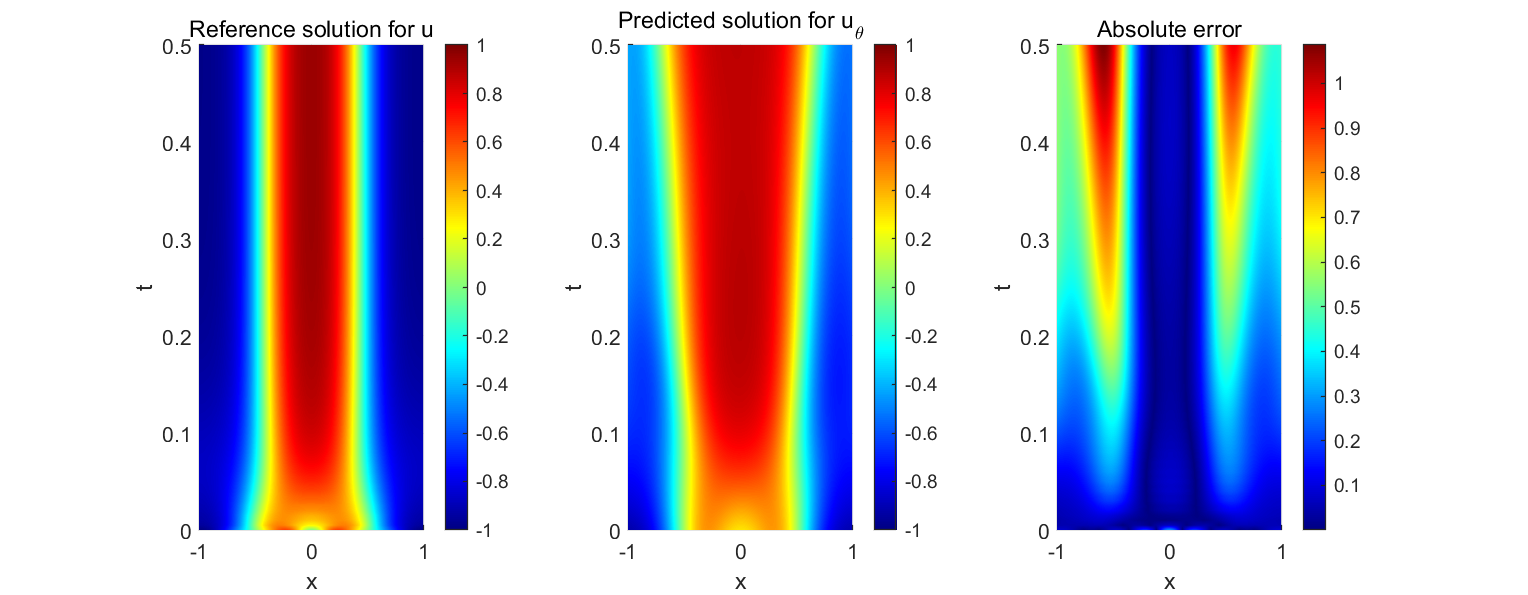}}
\caption{\textbf{Results for solving 1D Allen--Cahn and 1D Cahn--Hilliard equations using vanilla PINN.}\label{fig::failure}}
\end{figure}

According to the study conducted by Mattey et al. (2022) \cite{mattey2022novel}, the std-PINN encounters challenges when dealing with highly nonlinear PDEs, such as the Allen--Cahn equation and Cahn--Hilliard equation. Fig. \ref{fig::failure} illustrates the predicted solution for these two equations using the vanilla PINNs with the specific forms of equations (\ref{eq::ac1dnum}) and (\ref{eq::ch1dnum}). It is evident that the vanilla PINNs fail to provide accurate predictions. To address this issue, recent research has explored potential solutions through the implementation of self-adaptive sample points or self-adaptive weights. The objective is to achieve fair training across the entire computational domain. However, it is important to note that previous experiments primarily utilized deep or deeper neural networks with multiple layers ($L \geq 3$). In contrast, our work needs to adopt a shallow neural network ($u_{\theta}$ with $L=2$, i.e., a single hidden layer) combined with a self-adaptive weighting scheme.

{ We proposed a new formulation of the weighting scheme known as Residuals-RAE-PINNs \cite{zhangresiduals}, which incorporates the K-nearest algorithm. It is crucial to emphasize that Residuals-RAE-PINNs differs fundamentally from conventional self-adaptive PINNs (such as SA-PINNs \cite{mcclenny2023self}). In SA-PINNs, the weights are treated as trainable parameters that evolve simultaneously with the network parameters $\theta$ during gradient descent, creating a coupled optimization problem. In contrast, Residuals-RAE-PINNs employs a \textbf{pre-training weight update} mechanism: at the beginning of each training iteration, the weights are first computed based on the current residuals and then held fixed throughout that iteration's parameter update. This decoupled procedure ensures that within each gradient descent step, the network training operates with constant weights, thereby making the theoretical analysis tractable. This approach considers the PDE residuals from neighboring points when calculating the weights. The research showcases that this algorithm provides exceptional stability, resulting in enhanced accuracy for difficult learning problems. Moreover, during the process of error analysis, a crucial issue arises when attempting to approximate the general loss using discrete samples from the domain. Therefore, it is essential to ensure that the weights are bounded and stable. This is a vital consideration in our decision to choose this framework.}
The primary difference in the modified training loss for Residuals-RAE-PINNs lies in its utilization of a weighted approach to evaluate PDE residuals, as opposed to the conventional MSE loss. More specifically, the modified training loss for interior points can be defined as follows:
\begin{equation}
    \label{eq::modifiedloss}
\mathcal{J}(u(\cdot ; \theta)) = \gamma_{int} \|\hat{R}_{int}(\boldsymbol{x}_{int},t_{int} ; \theta)\|_{L^2(D \times [0,T])}^2 + \gamma_{tb} \|R_{tb}(\boldsymbol{x}_{tb}; \theta)\|_{L^2(D)}^2
+ \gamma_{sb} \|R_{sb}(\boldsymbol{x}_{sb},t_{sb} ; \theta)\|_{L^2(\partial D \times [0,T])}^2.
\end{equation}

In Eq. \eqref{eq::modifiedloss}, we introduce the squared $L^2$ norm of the modified internal residual { $\|\hat{R}_{int}(\boldsymbol{x}_{int},t_{int} ; \theta)\|_{L^2(D \times [0,T])}^2$}. The value of $\hat{R}_{int}^{(i)}$ can be estimated by multiplying the self-adaptive pointwise weights, $\lambda_{int}^{(i)}$, and the original internal residual, ${R}_{int}^{(i)}$. Here, $(\bullet)^{(i)}$ represents the index of the residual point sampled from the interior domain. The Residuals-RAE weighting scheme is designed to ensure stability in training and encourage the network to assign higher weights to sub-domains that are more difficult to learn. Initially, the Residuals-RAE-PINNs considers a simple weighting scheme:
\begin{align}
& w_{int}(\boldsymbol{x}_{int}^{(i)},t_{int}^{(i)}; \theta) = \frac{|{R}_{int}(\boldsymbol{x}_{int}^{(i)},t_{int}^{(i)}; \theta)|}{ {\textstyle \sum_{i=1}^{M_{int}}} |{R}_{int}(\boldsymbol{x}_{int}^{(i)},t_{int}^{(i)}; \theta)|}\cdot M_{int} , \ i=1,2,\dots,M_{int} ,
\label{eq::simpleRAE}
\end{align}
where the normalization is applied here to obtain a normalized weights. However, the weights obtained from a simple weighting scheme may not be stable because the PDE residuals can change rapidly, potentially leading to failed training \cite{han2022residual}. Inspired by the XPINNs \cite{jagtap2020extended}, {the Residuals-RAE PINNs \cite{zhangresiduals} }construct a set of near points for each interior points to average the weights from its neighboring points:
\begin{align}
& \lambda_{Knear,r}^{(i)} = \frac{1}{k_{int}}\sum_{(\boldsymbol{x}_{int}^{(j)},t_{int}^{(j)}) \in \boldsymbol{S}_{k_{int}}^{(i)}} w_{int}(\boldsymbol{x}_{int}^{(j)},t_{int}^{(j)}; \theta), \ i=1,2,\dots,M_{int}, \\ \label{eq::lknear}
& \boldsymbol{\lambda}_{int}^k = \beta  \boldsymbol{\lambda}_{Knear,r}^k + (1-\beta) \boldsymbol{\lambda}^{k-1}_{int}.
\end{align}

{  In this paper, \textbf{bold} symbols denote \textbf{vectors}; non-bold symbols denote scalars. Thus $\boldsymbol{\lambda}_{int}^k$, $\boldsymbol{\lambda}_{Knear,r}^k \in \mathbb{R}^{M_{int}}$ are the vectors of interior and K-near-normalized weights at iteration $k$, with $i$-th components $\lambda_{int}^{(i)}$ and $\lambda_{Knear,r}^{(i)}$ respectively.}

Here, $\boldsymbol{S}_{k_{int}}^{(i)}$ represents the set of near points for the $i$-th point $(\boldsymbol{x}_{int}^{(i)}, t_{int}^{(i)})$, where  $k_{int}$ denotes the number of elements in the set. The scalar $\lambda_{Knear,r}^{(i)}$ is a middle weight (after normalization from neighboring points) at the $i$-th collocation point. The superscript $(\bullet)^k$ stands for the training iteration index. {  \textbf{From pointwise to global (Eq.\,\eqref{eq::lknear}).} The \textbf{pointwise} quantity $\lambda_{Knear,r}^{(i)}$ in the first line is defined for each $i=1,\ldots,M_{int}$; it depends on $\theta$ through $w_{int}(\cdot;\theta)$. At iteration $k$, we evaluate this formula for \textbf{all} $i$ using the current parameters $\theta^k$ (in practice, $\theta^{k-1}$ when weights are updated before the $k$-th gradient step). The \textbf{global} vector $\boldsymbol{\lambda}_{Knear,r}^k$ in the second line is defined as the collection of these pointwise values at iteration $k$: $\boldsymbol{\lambda}_{Knear,r}^k = \bigl( \lambda_{Knear,r}^{(1)}(\theta^k), \lambda_{Knear,r}^{(2)}(\theta^k), \ldots, \lambda_{Knear,r}^{(M_{int})}(\theta^k) \bigr)^\top$. The second line of Eq.\,\eqref{eq::lknear} is then the vector update: the new weight vector $\boldsymbol{\lambda}_{int}^k$ is a convex combination of the previous vector $\boldsymbol{\lambda}_{int}^{k-1}$ and the current K-near vector $\boldsymbol{\lambda}_{Knear,r}^k$.}

{ 
To explicitly illustrate the implementation details of the mathematical formulation described above, the complete procedure is outlined in \textbf{Algorithm \ref{alg:residuals_rae}} (provided in \textbf{Appendix \ref{sec:appendix_algorithm}}). This algorithmic flowchart visually reinforces the decoupled nature of our training strategy, demonstrating step-by-step how the weights are computed and frozen prior to the network's parameter update in each iteration.
}

Further information on the design of a weighting scheme can be found in our upcoming publication titled "Residuals-RAE-PINNs." {  A rigorous proof that the weights $\boldsymbol{\lambda}^k_{int}$ converge to a limit $\boldsymbol{\lambda}^*$ as training progresses is given in Theorem~\ref{thm:weight_convergence} (Section~\ref{sec:weight_convergence}, Weight Convergence Analysis), under standard assumptions on the optimizer and the residual continuity. For the reader's convenience, we summarize the key intermediate step here. Iterating the update rule in Eq.~\eqref{eq::lknear}, $\boldsymbol{\lambda}^{k}_{int} = \beta \boldsymbol{\lambda}_{Knear,r}^k + (1-\beta) \boldsymbol{\lambda}^{k-1}_{int}$, yields
\begin{equation}
\label{eq:weight_expansion}
\boldsymbol{\lambda}^{k}_{int} = \beta \sum_{i=1}^{k} (1-\beta)^{k-i} \boldsymbol{\lambda}_{Knear,r}^{i} + (1-\beta)^k \boldsymbol{\lambda}^{0}_{int}.
\end{equation}
Since $\beta \in (0,1)$, the factor $(1-\beta)^k \to 0$ as $k \to \infty$, so the initial term vanishes. Under Assumptions A1--A3, the sequence $\boldsymbol{\lambda}_{Knear,r}^{i}$ converges as $\theta^i \to \theta^*$ (Assumption A2 and continuity of residuals in $\theta$), so the weighted sum on the right-hand side converges to a limit $\boldsymbol{\lambda}^*$. Thus the weights achieve stability when the training process approaches convergence; stability can also be promoted in practice by using a small or decaying learning rate so that $\theta^k$ evolves smoothly toward $\theta^*$. This ensures that the pointwise weights in the interior domain remain stable as the number of training epochs increases. The full statement and proof are in Theorem~\ref{thm:weight_convergence}.}

Furthermore, the adjusted generalization error of a Residuals-RAE-PINN can be determined by
\begin{equation}
\label{eq::Rae_SA}
   \begin{aligned}
\mathcal{E}_G(\theta)^2=  \int_{D \times [0,T]}\left|R_{i n t }(\boldsymbol{x}_{i n t },t_{i n t })\right|^2 \mathrm{~d} \boldsymbol{x} \mathrm{~d} t  +\int_{ D }\left|R_{tb}(\boldsymbol{x}_{tb})\right|^2\mathrm{~d} \boldsymbol{x} +\int_{\partial D \times [0,T]}\left|R_{sb}(\boldsymbol{x}_{sb}, t_{sb})\right|^2  \mathrm{~d} s(\boldsymbol{x})\mathrm{~d} t.
\end{aligned} 
\end{equation}

To simplify the calculation of the generalization error, we can approximate the integrals in Definition (\ref{eq::Rae_SA}) above using {  Monte-Carlo (LHS) sampling}, as explained in Lemma \ref{lemma::quadrules}. As a result, we define the training loss as follows:
\begin{align} 
\mathcal{E}_T(\theta,S)^2=\mathcal{E}_T^{int}(\theta,S_{int})^2+\mathcal{E}_T^{tb}(\theta,S_{tb})^2+\mathcal{E}_T^{sb}(\theta,S_{sb})^2,
\end{align}
where 
\begin{equation*}
    \begin{aligned}
    \mathcal{E}_T^{int}(\theta,S_{int})^2=&~\sum_{i=1}^{M_{int}}\lambda_{int}^{(i)}|R_{int}(\boldsymbol{x}_{int}^{(i)},t_{int}^{(i)})|^2,\\
    \mathcal{E}_T^{tb}(\theta,S_{tb})^2=&~\sum_{i=1}^{M_{tb}}\lambda_{tb}^{(i)}|R_{tb}(\boldsymbol{x}_{tb}^{(i)})|^2,\\
    \mathcal{E}_T^{sb}(\theta,S_{sb})^2=&~\sum_{i=1}^{M_{sb}}\lambda_{sb}^{(i)}|R_{sb}(\boldsymbol{x}_{sb}^{(i)},t_{sb}^{(i)})|^2.
\end{aligned}
\end{equation*}

Here, the weights $\lambda_{int}^{(i)}$ in \eqref{eq::lknear} are obtained from the Residuals-RAE-PINNs, while {  the other weights $\lambda_{sb}^{(i)}$ and $\lambda_{tb}^{(i)}$ are related to the Monte-Carlo (LHS) sampling}.  The definition of the residuals remains the same as in Definition (\ref{eq::base_loss}).

{ 
\subsection{Weight Convergence Analysis of Residuals-RAE-PINNs}
\label{sec:weight_convergence}

In this subsection, we establish the theoretical foundation for the convergence of the weights in Residuals-RAE-PINNs. We first state the necessary assumptions regarding neural network training, and then prove that the weights converge to a stable limit under these assumptions.

The following assumptions are standard in the analysis of neural network optimization and are supported by extensive theoretical and empirical studies.

\begin{description}
\item[\textbf{Assumption A1} (Smoothness of the Loss Landscape).] The loss function $\mathcal{L}(\theta)$ is continuously differentiable with respect to the network parameters $\theta$, and its gradient $\nabla_\theta \mathcal{L}(\theta)$ is Lipschitz continuous with constant $L_g > 0$:
$$
\|\nabla_\theta \mathcal{L}(\theta_1) - \nabla_\theta \mathcal{L}(\theta_2)\| \leq L_g \|\theta_1 - \theta_2\|, \quad \forall \theta_1, \theta_2.
$$
This assumption is commonly adopted in the convergence analysis of gradient-based optimization methods \cite{nesterov2018lectures, bottou2018optimization}.

\item[\textbf{Assumption A2} (Optimizer Convergence).] The optimization algorithm (e.g., Adam with L-BFGS refinement) produces a sequence of parameters $\{\theta^k\}_{k=1}^{\infty}$ that converges to a stationary point $\theta^*$:
$$
\lim_{k \to \infty} \theta^k = \theta^*, \quad \text{with} \quad \nabla_\theta \mathcal{L}(\theta^*) = 0.
$$
The convergence of Adam and L-BFGS for smooth, non-convex objectives has been established under mild conditions \cite{reddi2019convergence, liu1989limited}.

\item[\textbf{Assumption A3} (Continuity of Residual Functions).] The residual function $R_{int}(\boldsymbol{x}, t; \theta)$ is continuous with respect to $\theta$ for all $(\boldsymbol{x}, t) \in D \times [0, T]$. This property holds naturally for neural networks with smooth activation functions such as $\tanh$.
\end{description}

Under the above assumptions, we establish the convergence of the weights in Residuals-RAE-PINNs. Since Residuals-RAE-PINNs employs a pre-training weight update mechanism (weights are computed before each training step and held constant during the gradient update), the following theorem demonstrates that the weight sequence converges as the network training converges.

\begin{theorem}
\label{thm:weight_convergence}
(Weight Convergence in Residuals-RAE-PINNs) Let Assumptions A1--A3 hold. Consider the weight update rule defined in Eq.~\eqref{eq::lknear}:
$$
\boldsymbol{\lambda}^{k}_{int} = \beta  \boldsymbol{\lambda}_{Knear,r}^k + (1-\beta) \boldsymbol{\lambda}^{k-1}_{int}, \quad \beta \in (0,1).
$$
If the network parameters $\theta^k$ converge to $\theta^*$ as $k \to \infty$ (Assumption A2), then the weights $\boldsymbol{\lambda}^k_{int}$ converge to a limit $\boldsymbol{\lambda}^*$:
$$
\lim_{k \to \infty} \boldsymbol{\lambda}^k_{int} = \boldsymbol{\lambda}^*.
$$
\end{theorem}

\begin{proof}
We expand the recursive relation iteratively. Applying the update rule $k$ times yields:
$$
\boldsymbol{\lambda}^{k}_{int} = \beta \sum_{i=1}^{k} (1-\beta)^{k-i} \boldsymbol{\lambda}_{Knear,r}^{i} + (1-\beta)^k \boldsymbol{\lambda}^{0}_{int}.
$$
Since $\beta \in (0,1)$, we have $(1-\beta)^k \to 0$ as $k \to \infty$, so the initial condition term vanishes asymptotically.

By Assumption A3, the residual function $R_{int}(\cdot; \theta)$ is continuous in $\theta$. Combined with Assumption A2 ($\theta^k \to \theta^*$), {  the vector $\boldsymbol{\lambda}_{Knear,r}^{k}$ (whose $i$-th component is the pointwise $\lambda_{Knear,r}^{(i)}$ in Eq.~\eqref{eq::lknear} evaluated at $\theta = \theta^k$) satisfies, for each $i$,
$$
\bigl(\boldsymbol{\lambda}_{Knear,r}^{k}\bigr)_i = \lambda_{Knear,r}^{(i)}(\theta^k) = \frac{1}{k_{int}}\sum_{(\boldsymbol{x}_{int}^{(j)},t_{int}^{(j)}) \in \boldsymbol{S}_{k_{int}}^{(i)}} w_{int}(\boldsymbol{x}_{int}^{(j)},t_{int}^{(j)}; \theta^k),
$$
and therefore} the intermediate weight vectors $\boldsymbol{\lambda}_{Knear,r}^{k}$ converge to a stable limit $\boldsymbol{\lambda}_{Knear,r}^*$ as $k \to \infty$.

For sufficiently large $k$, we have $\|\boldsymbol{\lambda}_{Knear,r}^{k} - \boldsymbol{\lambda}_{Knear,r}^*\| < \epsilon$ for any $\epsilon > 0$. The series $\sum_{i=1}^{\infty} (1-\beta)^{i-1}$ converges to $1/\beta$ due to the geometric decay. Therefore:
$$
\boldsymbol{\lambda}^* = \beta \sum_{i=1}^{\infty} (1-\beta)^{i-1} \boldsymbol{\lambda}_{Knear,r}^* = \boldsymbol{\lambda}_{Knear,r}^*.
$$
We conclude that $\boldsymbol{\lambda}^k_{int} \to \boldsymbol{\lambda}^*$ as $k \to \infty$.
\end{proof}

\begin{remark}
\label{rem:fixed_weight_validity}
Theorem~\ref{thm:weight_convergence} provides the theoretical justification for analyzing Residuals-RAE-PINNs using the framework of weighted PINNs with constant weights. Since the weights converge to $\boldsymbol{\lambda}^*$ as training progresses, the error analysis conducted with the converged weights accurately characterizes the approximation quality of the final trained network. The pre-training weight update mechanism of Residuals-RAE-PINNs ensures that within each training iteration, the weights remain constant, which is consistent with our theoretical framework.
\end{remark}

\begin{remark}
The convergence rate of $\boldsymbol{\lambda}^k_{int}$ is governed by both the optimizer convergence rate and the discount factor $\beta$. Smaller values of $\beta$ accelerate weight convergence but may reduce responsiveness to recent residual changes. In practice, $\beta \in [0.9, 0.999]$ provides a suitable balance.
\end{remark}

\begin{remark}
The spatial smoothing via the K-nearest algorithm and the temporal smoothing via the exponential moving average in Eq.~\eqref{eq::lknear} together prevent the rapid weight oscillations observed in simpler weighting schemes (Eq.~\eqref{eq::simpleRAE}). This stability is essential for the convergence guarantee established above.
\end{remark}
}

The following sections establish a priori error bounds for the Residuals-RAE-PINN solutions of the Allen--Cahn and Cahn--Hilliard equations. The proofs rely on the weight convergence result (Theorem~\ref{thm:weight_convergence}) and the Monte-Carlo (LHS) sampling lemma (Lemma~\ref{lemma::quadrules}).

\section{Physics Informed Neural Networks for Approximating the Allen--Cahn Equation}\label{section::3}

In this section we present the error analysis for the Allen--Cahn equation: we state the main error bounds (Theorems~\ref{theorem::ac4} and \ref{theorem::ac5}) and their proofs, and then report numerical examples.

\subsection{\textbf{Allen--Cahn Equation}}
Let us consider the general form of the Allen--Cahn equation defined in the domain $\Omega := D\times[0,T]$, where $D\subset\mathbb{R}^d$, $d\in\mathbb{N}$ and $\partial D$ denotes the boundary of the domain $D$. 
The Allen--Cahn equation is given by:
\begin{subequations}\label{eq::ac}
\begin{align}
    & u_t-\epsilon^2 \nabla^2 u+f(u)=0, \quad t \in[0, T], \quad \boldsymbol{x} \in D  \\ 
    & u(\boldsymbol{x}, 0)=\psi(\boldsymbol{x}), \quad \boldsymbol{x} \in D  \\
    & u(\boldsymbol{x}, t)=u(-\boldsymbol{x}, t), \quad t \in[0, T], \quad \boldsymbol{x} \in \partial D \\ 
    & \nabla u(\boldsymbol{x}, t)=\nabla u(-\boldsymbol{x}, t) , \quad t \in[0, T], \quad \boldsymbol{x} \in \partial D
\end{align}
\end{subequations}

Here, the diffusion coefficients $\epsilon$ is a constant,  the nonlinear term $f(u)$ is a polynomial in $u,$ and   $u(\boldsymbol{x}, t)$ represents the unknown field solution,  and $\psi$ denotes the initial distributions for $u$ at $t=0$. The boundary conditions have been assumed to be periodic.

Furthermore, we are able to derive the following regularity results.
{
\begin{lemma}[\cite{fukao2004some}, \cite{bartels2015numerical}]\label{lemma::ac_exist}
Let $r \geq 1$, $\psi \in H^r(D)$, and assume $f \in C^1(\mathbb{R})$ with $|f'(s)| \leq C(1 + |s|^{p-1})$ for some $p \geq 1$ and $C > 0$. Then there exists a unique solution $u$ to the Allen-Cahn equation such that $u \in C([0,T]; H^r(D)) \cap C^1([0,T]; H^{r-2}(D))$ for some $T > 0$.
\end{lemma}

\begin{lemma}\label{acle2}
Let $k \in \mathbb{N}$, $\psi \in H^{r+2k}(D)$ with $r > \frac{d}{2} + 2k$, and assume $f \in C^k(\mathbb{R})$ with bounded derivatives up to order $k$. If $f^{(j)}$ has polynomial growth of degree at most $p-j$ for $j = 0, 1, \ldots, k$ and some $p \geq 1$, then there exists $T > 0$ and a solution $u$ to the Allen-Cahn equation such that $u(\cdot, 0) = \psi$, $u \in H^{r+2k}(D \times [0,T])$, and $u_t \in H^{r+2k-2}(D \times [0,T])$.
\end{lemma}}


The proof for Lemma \ref{acle2} is provided in the Appendix \ref{Proofs}.

\subsection{\textbf{Physics Informed Neural Networks}}
In this section, we present a physics-informed neural network that can be used to estimate the solution of partial differential equations (PDEs). The network, denoted as $u_{\theta}: \Omega\rightarrow \mathbb{R}$, is parameterized by $\theta$ and operates on the domain $\Omega:= D\times[0, T]$. Subsequently, we can derive the following: 
\begin{subequations}
\label{eq::ac_{int}esidual1}
\begin{align}
& R_{\text {int }}^{AC}\left[u_\theta\right](\boldsymbol{x}, t)=\partial_tu_\theta - \epsilon^2\nabla^2u_\theta+f(u_\theta), \\
& R_{t b}^{AC}\left[u_\theta\right](\boldsymbol{x})=u_\theta(\boldsymbol{x}, 0)-\psi(\boldsymbol{x}), \\
& R_{s b 1}^{AC}\left[u_\theta\right](\boldsymbol{x}, t)=u_\theta(\boldsymbol{x}, t)-u_\theta(-\boldsymbol{x}, t), \\
& R_{s b 2}^{AC}\left[u_\theta\right](\boldsymbol{x}, t)=\nabla u_\theta(\boldsymbol{x}, t)-\nabla u_\theta(-\boldsymbol{x}, t).
\end{align}
\end{subequations}

The equations above can be derived promptly. PINNs convert the task of solving equations into an approximation problem, resulting in the generalization error as below
\begin{equation}\label{eq::ac_GenError}
\begin{aligned}
\mathcal{E}^{AC}_G(\theta)^2=
&\int_D|R_{tb}^{AC}|^2d\boldsymbol{x}
+\int_0^T\int_D|R_{int}^{AC}|^2d\boldsymbol{x}dt\\
&+\epsilon^2\tilde{C}T^{1/2}\left(\int_0^T\int_{\partial D}|R_{sb1}^{AC}|^2 ds(\boldsymbol{x})dt \right)^{1/2}
+\epsilon^2\tilde{C}T^{1/2}\left(\int_0^T\int_{\partial D}|R_{sb2}^{AC}|^2 ds(\boldsymbol{x})dt\right)^{1/2},
\end{aligned}
\end{equation}
where $\tilde{C}=|\partial D|^{\frac{1}{2}}(||{u}||_{C^1(\partial D\times[0,T])}+||{u}_\theta||_{C^1(\partial D\times[0,T])})$.

In practice, the estimation of the generalization error $\mathcal{E}_G^{AC}(\theta)$ can be approximated using {  Monte-Carlo (LHS) sampling} and the training points from the {  Monte-Carlo (LHS) sampling} method. The training set, denoted as $\mathcal{S}^{AC}$, consists of {  sample points from Monte-Carlo (LHS) sampling} for the inner domain, boundary, and initial conditions. This set can be further divided into three subsets: $\mathcal{S}^{AC}_{int}=\{(\boldsymbol{x}^{(i)}_{int},t^{(i)}_{int})\}_{i=1}^{M_{int}}$, $\mathcal{S}^{AC}_{tb}=\{(\boldsymbol{x}^{(i)}_{tb})\}_{i=1}^{M_{tb}}$, and $\mathcal{S}^{AC}_{sb}=\{(\boldsymbol{x}^{(i)}_{sb},t^{(i)}_{sb})\}_{i=1}^{M_{sb}}$. Here, $M_{int}$, $M_{sb}$, and $M_{tb}$ represent the number of training points for the governing equation, boundary conditions, and initial conditions, respectively. The training error can then be expressed as follows:
\begin{equation}
\label{eq::AC_trainerror}
\begin{aligned}
\mathcal{E}^{AC}_T(\theta,\mathcal{S}^{AC})^2:=&\sum_{i=1}^{M_{tb}}\lambda_{tb}^{(i)}|R_{tb}^{AC}(\boldsymbol{x}_{tb}^{(i)})|^2+\sum_{i=1}^{M_{int}}\lambda_{int}^{(i)}|R_{int}^{AC}(\boldsymbol{x}_{int}^{(i)},t_{int}^{(i)})|^2\\
&+\epsilon^2
\tilde{C}T^{1/2}\left(\sum_{i=1}^{M_{sb}}\lambda_{sb}^{(i)}|R_{sb1}^{AC}(\boldsymbol{x}_{sb}^{(i)},t_{sb}^{(i)})|^2\right)^{1/2}
+\epsilon^2
\tilde{C}T^{1/2}\left(\sum_{i=1}^{M_{sb}}\lambda_{sb}^{(i)}|R_{sb2}^{AC}(\boldsymbol{x}_{sb}^{(i)},t_{sb}^{(i)})|^2\right)^{1/2}.
\end{aligned}
\end{equation}

{ Here, the constant $\epsilon$ is the physical diffusion coefficient from the governing Allen-Cahn equation (Eq. \ref{eq::ac_{int}esidual1}). The temporal factor $T^{1/2}$ and the $1/2$ power on the boundary integrals---representing their $L^2$ norms---both result from applying the Cauchy-Schwarz inequality to bound the boundary terms during the derivation, as detailed in the proof of Theorem~\ref{theorem::ac4}. }


Next, we use $\hat{u} = u_{\theta} - u$ to represent the numerical error between the predicted solution from PINNs and the exact solution. 
Then we define the total error of PINNs with the following form with $\epsilon$ being the physical diffusion coefficient of the Allen-Cahn equation:
\begin{definition}[Total Error for Allen--Cahn System]
    The total error for the Allen–Cahn (AC) equation, denoted by $\mathcal{E}^{AC}(\theta)^2$, is a comprehensive metric that quantifies both the deviation in solution values and the discrepancy in spatial gradients between the predicted and exact solutions. It is defined as follows:
    \begin{equation}
    \label{eq::ac_total}
    \mathcal{E}^{AC}(\theta)^2=\|\hat{u}(\boldsymbol{x},t)\|^2_{L^2{(D\times [0,T])}}+{ 2\epsilon^2\|\nabla\hat{u}\|^2_{L^2(D \times [0,T])}}.
\end{equation}
Here: 
\begin{itemize}
    \item $\hat{u} = u_{\theta} - u$ denotes the residual error between PINNs' prediction $u_\theta$ and the exact solution $u$. 
    \item $\epsilon$ is the diffusion coefficient in the Allen-Cahn equation, controlling the interface width.
    \item The first term, $\|\hat{u}(\boldsymbol{x},t)\|^2_{L^2(D\times [0,T])}$, represents the standard $L^2$ error over the full spatiotemporal domain, quantifying the average discrepancy between the predicted solution and the exact reference.
    \item The second term, { $2\epsilon^2 \|\nabla \hat{u}\|^2_{L^2(D\times[0,T])}$}, measures the $L^2$ error of the spatial gradients over the spatiotemporal domain, emphasizing the fidelity of the predicted interface structure and smoothness.
\end{itemize}
\end{definition}


{


\textbf{Remark 3.1.} The structure of the total error, $\mathcal{E}^{AC}(\theta)^2$, is a direct consequence of the analytical technique for parabolic PDEs, as detailed in the proof of Theorem \ref{theorem::ac4}. The second term, which captures the gradient error, is not arbitrary; its form arises from the proof methodology. The derivation establishes a differential inequality for the error's norm, which is then resolved using Grönwall's inequality. This integration over time naturally yields the final structure of the gradient term, making this error metric fundamental to the subsequent stability analysis.

In the next section, we will undertake a detailed analysis of  $\mathcal{E}^{AC}(\theta)^2$.}

\subsection{\textbf{Error Analysis}}
{ The subsequent error analysis is conducted for the converged network, where the weights have stabilized to their limiting values as established in Theorem~\ref{thm:weight_convergence}. Since Residuals-RAE-PINNs employs a pre-training weight update mechanism, the weights remain constant within each training iteration, which is consistent with the theoretical framework presented below.}

The Allen--Cahn equations and the definitions for different residuals can be re-expressed as follows:
\begin{subequations}
\label{eq::ac_residual2}
\begin{align}
& R_{\text {int}}^{AC}(\boldsymbol{x}, t)=\hat{u}_t-\epsilon^2\nabla^2\hat{u}-f(u)+f(u_\theta),\\
& R_{t b}^{AC}(\boldsymbol{x})=\hat{u}(\boldsymbol{x}, 0), \\
& R_{s b 1}^{AC}(\boldsymbol{x}, t)=\hat{u}(\boldsymbol{x}, t)-\hat{u}(-\boldsymbol{x}, t), \\
& R_{s b 2}^{AC}(\boldsymbol{x}, t)=\nabla\hat{u}(\boldsymbol{x}, t)-\nabla\hat{u}(-\boldsymbol{x}, t), 
\end{align}
\end{subequations}
where $\hat{u} = u_\theta - u$ is substituted into Eq. (\ref{eq::ac_{int}esidual1}) to obtain the above equations.

\subsubsection{Bound on the Residuals}
{\begin{theorem}\label{theorem::ac3} 
Let $n\geq 2$, $d, r, k \in \mathbb{N}$ with $k \geq 3$. Let $\psi \in H^{r+2k}(D)$ with $r > \frac{d}{2} + 2k$, and assume $f \in C^k(\mathbb{R})$ with bounded derivatives up to order $k$. If $f^{(j)}$ has polynomial growth of degree at most $p-j$ for $j = 0, 1, \ldots, k$ and some $p \geq 1$. For every integer $N>5$, there exists a tanh neural network $u_\theta$ with two hidden layers, of width at most $3\left\lceil\frac{k+n-2}{2}\right\rceil\left|P_{k-1, d+2}\right|+\lceil N T\rceil+d(N-1)$, and   $3\left\lceil\frac{d+n+1}{2}\right\rceil\left|P_{d+2, d+2}\right|\lceil N T\rceil N^d$ such that      
\begin{subequations}     
\label{eq::ac_{int}esidualBound}     
\begin{align}      
& \left\|R_{\text {int }}^{AC}\right\|_{L^2(D)} ,\left\|R_{s b 2}^{AC}\right\|_{L^2(\partial D \times [0, T])}\lesssim \ln^2 \left(N\right) N^{-k+2}, \\      
& \left\|R_{s b 1}^{AC}\right\|_{L^2(\partial D \times[0, T])},\left\|R_{\text {tb }}^{AC}\right\|_{L^2(D)}  \lesssim \ln \left(N\right) N^{-k+1} .     
\end{align}     
\end{subequations} 
\end{theorem}
}

\textbf{Proof of Theorem \ref{theorem::ac3}.}
{
Given that $f \in C^k(\mathbb{R})$ with bounded derivatives up to order $k$, we first seek to bound the nonlinear term $f(u_\theta) - f(u)$. By the Mean Value Theorem, there exists a value $\xi$ between $u$ and $u_\theta$ such that $|f(u_\theta) - f(u)| = |f'(\xi)| |u_\theta - u|$. Since $f \in C^k(\mathbb{R})$ with bounded derivatives and the solution $u$ is bounded in $H^{r+2k}(D \times [0,T])$ (according to Lemma \ref{acle2}), the derivative $f'$ is bounded on the relevant domain. We can therefore define a positive constant $M = \sup |f'(s)|$, which yields the key inequality: 
\begin{equation}     
|f(u_\theta) - f(u)| \le M|\hat{u}|, \quad \text{where } \hat{u} = u_\theta - u. 
\end{equation}}
Hence, the interior residual term can be bounded as follows:
\begin{equation}
 R_{\text {int}}^{AC}(\boldsymbol{x}, t)=\hat{u}_t-\epsilon^2\nabla^2\hat{u} + \left(f(u_\theta)-f(u)\right)\leq\hat{u}_t+M|\hat{u}|-\epsilon^2\nabla^2\hat{u}.
\end{equation}

Based on Lemma \ref{acle2}, and considering Lemma \ref{lemma::Error_N}, we can conclude that there is a neural network denoted as $u_\theta$, which has two hidden layers and widths of $3\left\lceil\frac{k+n-2}{2}\right\rceil\left|P_{k-1, d+2}\right|+\lceil N T\rceil+d(N-1)$, such that for every $0 \leq l \leq 2$, we have the inequality
\begin{align}\label{eq::uErrorAboutN}
    \left\|\hat{u}\right\|_{H^l(\Omega)}=\left\|u_\theta-u\right\|_{H^l(\Omega)} \leq C_{l, k, d+1, u} \lambda_{l, u}\left(N\right) N^{-k+l}.
\end{align}

Here, $\lambda_{l, u}=2^l 3^{d+1}(1+\sigma) \ln ^l\left(\beta_{l, \sigma, d+1, u} N^{d+k+3}\right)$, $\sigma=\frac{1}{100}$, , and the definition for the other constants can be found in Lemma \ref{lemma::Error_N}.

In light of Lemma \ref{lemma::trace},we can bound the terms of the PINN residual as follows:
\begin{subequations}
\begin{align}
& \left\|\hat{u}_t\right\|_{L^2(\Omega)} \leq\|\hat{u}\|_{H^1(\Omega)},   \|\nabla^2 \hat{u}\|_{L^2(\Omega)} \leq\|\hat{u}\|_{H^2(\Omega)}, \\
& \|\hat{u}\|_{L^2(D)} \leq\|\hat{u}\|_{L^2(\partial \Omega)} \leq C_{h_{\Omega}, d+1, \rho_{\Omega}}\|\hat{u}\|_{H^1(\Omega)}, \\
& \|\nabla \hat{u}\|_{L^2(\partial D \times[0, T])} \leq\|\nabla \hat{u}\|_{L^2(\partial \Omega)} \leq C_{h_{\Omega}, d+1, \rho_{\Omega}}\|\hat{u}\|_{H^2(\Omega)}.
\end{align}
\end{subequations}

By combining these relationships with inequality (\ref{eq::uErrorAboutN}), we can obtain
$$
\begin{aligned}
 \left\|R_{\text {int } }^{AC}\right\|_{L^2(\Omega)}
 &\leq \left\|\hat{u}_t+M|\hat{u}|-\epsilon^2\nabla^2\hat{u}\right\|_{L^2(\Omega)} \\
&\leq\|\hat{u}_t\|_{L^2(\Omega)}+M\|\hat{u}\|_{H^0(\Omega)}+\epsilon^2\|\nabla^2\hat{u}\|_{L^2(\Omega)} \\
&\leq\|\hat{u}\|_{H^1(\Omega)}+M\|\hat{u}\|_{H^0(\Omega)}+\epsilon^2\|\hat{u}\|_{H^2(\Omega)} \\
&\leq C_{1, k, d+1, u} \lambda_{1, u}(N) N^{-k+1}+ MC_{0, k, d+1, u} \lambda_{0, u}(N) N^{-k}\\
& \quad\quad\quad\quad\quad\quad\quad\quad\quad\quad\quad\quad\quad\quad\quad+\epsilon^2C_{2, k-1, d+1, u} \lambda_{2, u}(N) N^{-k+2} \\
&\lesssim \ln^2 \left(N\right) N^{-k+2}, \\
\end{aligned}
$$
$$
\begin{aligned}
&\left\|R_{t b }^{AC}\right\|_{L^2(D)} 
\leq C_{h_{\Omega}, d+1,
\rho_{\Omega}}\|\hat{u}\|_{H^1(\Omega)} \lesssim \ln \left(N\right) N^{-k+1},\\ 
&\left\|R_{s b 1}^{AC}\right\|_{L^2(\partial D \times[0, T])} \leq 2\cdot C_{h_{\Omega}, d+1, \rho_{\Omega}}\|\hat{u}\|_{H^1(\Omega)} \lesssim \ln \left(N\right) N^{-k+1},\\ 
&\left\|R_{s b 2}^{AC}\right\|_{L^2(\partial D \times[0, T])} \leq 2\cdot C_{h_{\Omega}, d+1, \rho_{\Omega}}\|\hat{u}\|_{H^2(\Omega)} \lesssim \ln ^2 \left(N\right) N^{-k+2} .\\
\end{aligned}
$$

\subsubsection{Bounds on the Total Approximation Error}
In this section, we will demonstrate that the total error $\mathcal{E}^{AC}(\theta)^2$ can be controlled by the generalization error $\mathcal{E}^{AC}_G(\theta)^2$. Furthermore, we will establish that the total error $\mathcal{E}^{AC}(\theta)^2$ can be made arbitrarily small, provided that the training error $\mathcal{E}^{AC}_T(\theta,S^{AC})^2$ is kept sufficiently small and the sample set sufficiently large.
{\begin{theorem}
\label{theorem::ac4}
Let $d \in \mathbb{N}$ and $u \in C^1(\Omega)$ be the classical solution to the Allen--Cahn equations. Let $u_\theta$ denote the PINN approximation with parameter $\theta$. Then the following relation holds,
    \begin{subequations}
    $$
{  \mathcal{E}^{AC}(\theta)^2\leq C_G\,\frac{\exp((1+2M)T)-1}{1+2M},}
$$
where
$$
\begin{aligned}
{  C_G=\int_0^T\int_D|R_{int}^{AC}|^2d\boldsymbol{x}dt}
    &+\epsilon^2\tilde{C}T^{1/2}\left(\int_0^T\int_{\partial D}|R_{sb1}^{AC}|^2 ds\left(\boldsymbol{x}\right)dt\right)^{1/2}\\
&+\epsilon^2\tilde{C}T^{1/2}\left(\int_0^T\int_{\partial D}|R_{sb2}^{AC}|^2 ds\left(\boldsymbol{x}\right)dt\right)^{1/2},\\
\end{aligned}
$$
    \end{subequations}
    $M$ is constant and $\tilde{C}=|\partial D|^{\frac{1}{2}}(\|{u}\|_{C^1\left(\partial D\times[0,T]\right)}+\|{u}_\theta\|_{C^1\left(\partial D\times[0,T]\right)})$
.
\end{theorem}}

\textbf{Proof of Theorem \ref{theorem::ac4}.}
Considering
\begin{equation}\label{eq::1}
    \begin{aligned}
        \int_D R_{int}^{AC}\cdot\hat{u} d\boldsymbol{x}
        &=\int_D\hat{u_t}\hat{u}d\boldsymbol{x}-\epsilon^2\int_{D}\nabla^2\hat{u}\cdot\hat{u}d\boldsymbol{x}+\int_D\left(f(u_\theta)-f(u)\right)\hat{u}d\boldsymbol{x}\\
       & ={ \frac{1}{2}\frac{d}{dt}}\int_D|\hat{u}|^2d\boldsymbol{x}+\epsilon^2\int_D|\nabla\hat{u}|^2d\boldsymbol{x}-\epsilon^2\int_{\partial D}\hat{u}\nabla\hat{u}\cdot\bm{n}ds(\boldsymbol{x})+\int_D\left(f(u_\theta)-f(u)\right)\hat{u}d\boldsymbol{x}.
    \end{aligned}
\end{equation}

Due to the presence of the nonlinear term $f(u)$ in the Allen--Cahn equation  and Lemma \ref{acle2}, there exists a constant $M>0$ such that $|f(u_1)-f(u_2)|\leq M|u_1-u_2|$, $ \forall u_1,u_2\in C^2\left(\Omega\right)$. By combining this with Eq. (\ref{eq::1}), we have
\begin{equation}
\label{eq::dtint}
    \begin{aligned}
         &\quad{ \frac{1}{2}\frac{d}{dt}}\int_D
         |\hat{u}|^2d\boldsymbol{x}+
         \epsilon^2\int_D
         |\nabla\hat{u}
          |^2d\boldsymbol{x}  \\
          & =\epsilon^2\int_{\partial D}\hat{u}\nabla\hat{u}\cdot\bm{n}ds(\boldsymbol{x})-\int_D(f(u_\theta)-f(u))\hat{u}d\boldsymbol{x}+\int_D R_{int}^{AC}\cdot\hat{u} d\boldsymbol{x} \\
        & \leq \epsilon^2\int_{\partial D}\hat{u}\nabla\hat{u}\cdot\bm{n}ds(\boldsymbol{x})+M\int_D|\hat{u}|^2d\boldsymbol{x}+\int_D R_{int}^{AC}\cdot\hat{u} d\boldsymbol{x}\\
        & \leq \epsilon^2C_1\left(\int_{\partial D}|R_{sb1}^{AC}|^2ds(\boldsymbol{x})\right)^{1/2}+\epsilon^2C_2\left(\int_{\partial D}|R_{sb2}^{AC}|^2ds(\boldsymbol{x})\right)^{1/2}+M\int_D|\hat{u}|^2d\boldsymbol{x}+\int_DR_{int}^{AC}\cdot\hat{u} d\boldsymbol{x}\\
        & \leq \epsilon^2C_1\left(\int_{\partial D}|R_{sb1}^{AC}|^2ds(\boldsymbol{x})\right)^{1/2}+\epsilon^2C_2\left(\int_{\partial D}|R_{sb2}^{AC}|^2ds(\boldsymbol{x})\right)^{1/2}+\frac{1}{2}\int_D|R_{int}^{AC}|^2d\boldsymbol{x}+\frac{1+2M}{2}\int_D|\hat{u}|^2d\boldsymbol{x}\\
        & \leq \epsilon^2\tilde{C}\left(\int_{\partial D}|R_{sb1}^{AC}|^2ds(\boldsymbol{x})\right)^{1/2}+\epsilon^2\tilde{C}\left(\int_{\partial D}|R_{sb2}^{AC}|^2ds(\boldsymbol{x})\right)^{1/2}+\frac{1}{2}\int_D|R_{int}^{AC}|^2d\boldsymbol{x}+\frac{1+2M}{2}\int_D|\hat{u}|^2d\boldsymbol{x},
    \end{aligned}
\end{equation}
where $C_1=|\partial D|^{\frac{1}{2}}(||{u}||_{C^1(\partial D\times[0,T])}+||{u}_\theta||_{C^1(\partial D\times[0,T])}),C_2=|\partial D|^{\frac{1}{2}}(\|{u}\|_{C(\partial D\times[0,T])}+||{u}_\theta||_{C(\partial 
 D\times[0,T])})$,$\tilde{C}=C_1=|\partial D|^{\frac{1}{2}}(||{u}||_{C^1(\partial D\times[0,T])}+||{u}_\theta||_{C^1(\partial D\times[0,T])})\geq\max\{C_1,C_2\}$.

{The derivation from the third to the fourth line, which bounds the boundary integral terms, follows the approach presented in \cite{qian2023physics}.}

By integrating inequality (\ref{eq::dtint}) over $[0,T']$ for any $T'\leq T$ and applying Cauchy-Schwarz inequality, we can obtain 
{ 
$$
\begin{aligned}
    &\int_D|\hat{u}(\bm{x},T')|^2d\bm{x}+2\epsilon^2\int_0^{T'}\int_D|\nabla\hat{u}|^2d\bm{x}dt\\
    &\leq\int_0^T\int_D|R_{int}^{AC}|^2d\bm{x}dt+\epsilon^2\tilde{C}T^{1/2}\left(\int_0^T\int_{\partial D}|R_{sb1}^{AC}|^2 ds(\bm{x})dt\right)^{1/2}\\
    &\quad\quad\quad\quad+\epsilon^2\tilde{C}T^{1/2}\left(\int_0^T\int_{\partial D}|R_{sb2}^{AC}|^2 ds(\bm{x})dt\right)^{1/2}
    +(1+2M)\int_0^{T'}\int_D|\hat{u}|^2d\bm{x}dt.\\
\end{aligned}
$$
}
{ (The term $\int_D|R_{tb}^{AC}|^2d\bm{x}$ has been removed from the bound as it is superfluous.)}

{For convenience of notation, we define the quantity $C_G$ as follows:}
{ 
\begin{equation}
\begin{aligned} C_G=\int_0^T\int_D|R_{int}^{AC}|^2d\bm{x}dt
    &+\epsilon^2
    \tilde{C}T^{1/2}\left(\int_0^T\int_{\partial D}|R_{sb1}^{AC}|^2 ds(\bm{x})dt\right)^{1/2}\\
    &+\epsilon^2\tilde{C}T^{1/2}\left(\int_0^T\int_{\partial D}|R_{sb2}^{AC}|^2 ds(\bm{x})dt\right)^{1/2},\\
\end{aligned}
\end{equation}
}
(The term $\int_D|R_{tb}^{AC}|^2d\bm{x}$ is omitted from $C_G$ as it is superfluous for the bound.)
and setting
$$y(\xi)=\int_D|\hat{u}(\bm{x},\xi)|^2d\bm{x}+2\epsilon^2\int_0^{\xi}\int_D|\nabla\hat{u}|^2d\bm{x}d\tau,$$
we obtain the inequality: 
\begin{equation*}
    y(T')\leq C_G+(1+2M)\int_0^{T'}y(t)dt.
\end{equation*}

{  Applying the integral form of the Grönwall inequality (Lemma~\ref{lemma::Gronwall}) to $y(T')\leq C_G+(1+2M)\int_0^{T'}y(t)dt$, we obtain that the solution satisfies $y(t)\leq C_G\cdot \exp((1+2M)t)$ for each $t\in[0,T]$ (with $C_1=1+2M$, $C_2=C_G$ in the lemma; the standard form yields $y(t)\leq C_G\exp((1+2M)t)$). Integrating this bound over $t\in[0,T]$ gives
\begin{equation*}
    \int_0^T y(t)\,dt \leq C_G\int_0^T \exp((1+2M)t)\,dt = C_G\,\frac{\exp((1+2M)T)-1}{1+2M}.
\end{equation*}
The total error $\mathcal{E}^{AC}(\theta)^2$ is bounded above by $\int_0^T y(t)\,dt$ (since $y(t)$ dominates both the $L^2$ norm of $\hat{u}$ and the gradient term over $[0,t]$). Hence
\begin{equation*}
    \mathcal{E}^{AC}(\theta)^2 \leq C_G\,\frac{\exp((1+2M)T)-1}{1+2M}.
\end{equation*}
This yields the Theorem.}

{\begin{theorem}\label{theorem::ac5}
Let $d \in \mathbb{N}$ and $T>0$. Let $u \in C^4(\Omega)$  be the classical solution of the Allen--Cahn equations, and let $u_\theta$ denote the PINN approximation solution. The total error satisfies
\begin{equation}
\begin{aligned}
\mathcal{E}^{AC}(\theta)^2 
&{  \leq C_T\,\frac{\exp((1+2M)T)-1}{1+2M}}\\
& = \mathcal{O}\left(\mathcal{E}^{AC}_T\left(\theta,\mathcal{S}^{AC}\right)^2 + {  M_{t b}^{-\frac{1}{2}} + M_{i n t}^{-\frac{1}{2}} + M_{s b}^{-\frac{1}{2}}}\right).
\end{aligned}
\end{equation}
Here, the constant $C_T$ is defined as
\begin{equation}
\begin{aligned}
C_T= 
& {  C_{\left[\left(R_{t b }^{AC}\right)^2\right]} M_{t b}^{-\frac{1}{2}}}+\mathcal{Q}_{M_{t b}}^D\left[\left(R_{t b }^{AC}\right)^2\right]  +{  C_{\left[\left(R_{i n t }^{AC}\right)^2\right]} M_{i n t}^{-\frac{1}{2}}}+\mathcal{Q}_{M_{i n t}}^{\Omega}\left[\left(R_{i n t }^{AC}\right)^2\right] \\
& +\epsilon^2\tilde{C}T^{1/2}\left({  C_{\left[\left(R_{sb1}^{AC}\right)^2\right]} M_{s b}^{-\frac{1}{2}}}+\mathcal{Q}_{M_{s b}}^{\partial D\times \left[0,T\right]}\left[\left(R_{s b 1}^{AC}\right)^2\right]\right)^{1/2}\\
& +\epsilon^2
\tilde{C}T^{1/2}\left({  C_{\left[\left(R_{s b 2}^{AC}\right)^2\right]} M_{s b}^{-\frac{1}{2}}}+\mathcal{Q}_{M_{s b}}^{\partial D \times \left[0,T\right]}\left[\left(R_{s b 2}^{AC}\right)^2\right]\right)^{1/2} ,\\
\end{aligned}
\end{equation}
where 
$\tilde{C}=|\partial D|^{\frac{1}{2}}(||{u}||_{C^1(\partial D\times[0,T])}+||{u}_\theta||_{C^1(\partial D\times[0,T])})$ and $M_{tb}$, $M_{int}$ and $M_{sb}$ denote the number of {  Monte-Carlo (LHS) sample} points on the initial boundary, interior domain, and spatial boundary, respectively.
\end{theorem}}

\textbf{Proof of Theorem \ref{theorem::ac5}.} {  By combining Theorem \ref{theorem::ac4} with the Monte-Carlo (LHS) sampling error formula (\ref{equ::quadrules}), we can establish a bound for the generalization error constant $C_G$ in terms of the training error.} Recall from Theorem \ref{theorem::ac4} that $C_G$ is defined by the true integrals of the residuals. We now bound each of these integral terms as follows:
\begin{equation}\label{eq:quadrature_bounds_block}
\begin{gathered}
\begin{aligned}
\int_D\left|R_{t b }^{AC}\right|^2 \mathrm{d} \boldsymbol{x} & =\int_D\left|R_{t b }^{AC}\right|^2 \mathrm{d} \boldsymbol{x}-\mathcal{Q}_{M_{t b}}^D\left[\left(R_{t b }^{AC}\right)^2\right]+\mathcal{Q}_{M_{t b}}^D\left[\left(R_{t b }^{AC}\right)^2\right] \\
& \leq {  C_{\left[\left(R_{t b }^{AC}\right)^2\right]} M_{t b}^{-\frac{1}{2}}}+\mathcal{Q}_{M_{t b}}^D\left[\left(R_{t b }^{AC}\right)^2\right]
\end{aligned}
\\
\begin{aligned}
    \int_0^T\int_{D}\left|R_{i n t }^{AC}\right|^2 \mathrm{d} \boldsymbol{x} \mathrm{d} t & =\int_0^T\int_{D}\left|R_{i n t }^{AC}\right|^2 \mathrm{d} \boldsymbol{x} \mathrm{d} t-\mathcal{Q}_{M_{i n t}}^{\Omega}\left[\left(R_{i n t }^{AC}\right)^2\right]+\mathcal{Q}_{M_{i n t}}^{\Omega}\left[\left(R_{i n t }^{AC}\right)^2\right] \\
& \leq {  C_{\left[\left(R_{i n t }^{AC}\right)^2\right]} M_{i n t}^{-\frac{1}{2}}}+\mathcal{Q}_{M_{i n t}}^{\Omega}\left[\left(R_{i n t }^{AC}\right)^2\right]
\end{aligned}
\\
\begin{aligned}
    \int_0^T\int_{\partial D}\left|R_{sb1 }^{AC}\right|^2 ds(\boldsymbol{x}) \mathrm{d} t & =\int_0^T\int_{\partial D}\left|R_{sb1 }^{AC}\right|^2 ds(\boldsymbol{x}) \mathrm{d} t-\mathcal{Q}_{M_{sb}}^{\partial D\times[0,T]}\left[\left(R_{sb1 }^{AC}\right)^2\right]+\mathcal{Q}_{M_{sb}}^{\partial D\times[0,T]}\left[\left(R_{sb1 }^{AC}\right)^2\right] \\
& \leq {  C_{\left[\left(R_{sb1 }^{AC}\right)^2\right]} M_{sb}^{-\frac{1}{2}}}+\mathcal{Q}_{M_{sb}}^{\partial D\times[0,T]}\left[\left(R_{sb1}^{AC}\right)^2\right]
\end{aligned}
\\
\begin{aligned}
    \int_0^T\int_{\partial D}\left|R_{sb2 }^{AC}\right|^2 ds(\boldsymbol{x}) \mathrm{d} t & =\int_0^T\int_{\partial D}\left|R_{sb2 }^{AC}\right|^2 ds(\boldsymbol{x}) \mathrm{d} t-\mathcal{Q}_{M_{sb}}^{\partial D\times[0,T]}\left[\left(R_{sb2 }^{AC}\right)^2\right]+\mathcal{Q}_{M_{sb}}^{\partial D\times[0,T]}\left[\left(R_{sb2 }^{AC}\right)^2\right]\\
& \leq {  C_{\left[\left(R_{sb2 }^{AC}\right)^2\right]} M_{sb}^{-\frac{1}{2}}}+\mathcal{Q}_{M_{sb}}^{\partial D\times[0,T]}\left[\left(R_{sb2}^{AC}\right)^2\right]
\end{aligned}
\end{gathered}
\end{equation}
By substituting these individual bounds back into the definition of $C_G$ from Theorem \ref{theorem::ac4},

Applying this to each component of $C_G$, we find that $C_G$ is bounded by a new quantity, which we define as $C_T$:
\begin{equation}
    C_G \le C_T,
\end{equation}
where $C_T$ is given by
\begin{equation*}
    \begin{aligned}
C_T :=~
& \left( {  C_{\left[\left(R_{t b }^{AC}\right)^2\right]} M_{t b}^{-\frac{1}{2}}}+\mathcal{Q}_{M_{t b}}^D\left[\left(R_{t b }^{AC}\right)^2\right] \right) + \left( {  C_{\left[\left(R_{i n t }^{AC}\right)^2\right]} M_{i n t}^{-\frac{1}{2}}}+\mathcal{Q}_{M_{i n t}}^{\Omega}\left[\left(R_{i n t }^{AC}\right)^2\right] \right) \\
& +\epsilon^2\tilde{C}T^{1/2}\left({  C_{\left[\left(R_{s b 1}^{AC}\right)^2\right]} M_{s b}^{-\frac{1}{2}}}+\mathcal{Q}_{M_{s b}}^{\partial D\times[0,T]}\left[\left(R_{s b 1}^{AC}\right)^2\right]\right)^{1/2}\\
&+\epsilon^2\tilde{C}T^{1/2}\left({  C_{\left[\left(R_{s b 2}^{AC}\right)^2\right]} M_{s b}^{-\frac{1}{2}}}+\mathcal{Q}_{M_{s b}}^{\partial D\times[0,T]}\left[\left(R_{s b 2}^{AC}\right)^2\right]\right)^{1/2}.\\
\end{aligned}
\end{equation*}

Furthermore, {  the sampling constants} $C_{[\cdot]}$ can be bounded. For instance, for the term $C_{({R_{tb}^{AC}}^2)}$, we have:
\begin{equation*}
    C_{\left({R_{tb}^{AC}}^2\right)} \lesssim\|\hat{u}\|_{C^2}^2 \leq 2\left(\|u\|_{C^2}^2+\left\|u_\theta\right\|_{C^2}^2\right) \lesssim\|u\|_{C^2}^2+16^{2 L}(d+1)^8\left(e^2 2^4 W^3 R^2\|\sigma\|_{C^2}\right)^{4 L}.
\end{equation*}

Since similar bounds hold for the other {  sampling} constants, all terms in $C_T$ are bounded by quantities related to the training loss and network parameters. Substituting the resulting inequality $C_G \le C_T$ back into the main estimate from Theorem \ref{theorem::ac4} yields the desired result:
{ 
$$
\mathcal{E}^{AC}(\theta)^2=\|\hat{u}(\boldsymbol{x},t)\|^2_{L^2(D \times [0,T])}+2\epsilon^2\|\nabla\hat{u}\|^2_{L^2(D \times [0,T])} \leq C_T\,\frac{\exp((1+2M)T)-1}{1+2M}.
$$
}

\subsection{\textbf{Numerical Examples}}
\subsubsection{1D Allen--Cahn equation}
\label{sec:1D Allen--Cahn equation}
In this section, we will simulate the Allen--Cahn equation, which is a challenging partial differential equation (PDE) for classical Physics-Informed Neural Networks (PINNs) due to its stiffness and sharp transitions in the spatio-temporal domain.
Considering the space-time domain  
$D\times[0,T]=[-1,1]\times[0,1]$, the Allen--Cahn equation we will be discussing is for a specific one-dimensional time variation. It is accompanied by the following initial and periodic boundary conditions, as defined below:
\begin{subequations}
\begin{align}
\label{eq::ac1dnum}
& u_{t}-0.001u_{xx}+3u^3-3u=0,\quad (x,t)\in D\times [0,T]  \\
& u(x,0)=x^2 sin(2 \pi x),\quad (x,t)\in D\times\{0\} \label{eq::aci1} \\
& u(-x,t)=u(x,t),\quad  (x,t)\in \partial D\times[0,T]  \\
& \nabla u(-x,t)=\nabla u(x,t),\quad  (x,t)\in\partial D\times [0,T] .
\end{align}
\end{subequations}

The initial condition in Eq. (\ref{eq::aci1}) will be denoted as initial condition $\# 1$ in the following discussions. What's more, we also consider the initial condition $\# 2$ given by:
\begin{equation}
u(x,0)= cos\left(\pi x\right) - exp\left(-4\left(\pi x\right)^2\right).
\end{equation}

The Allen--Cahn equations will be denoted as AC-I1 and AC-I2, corresponding to initial conditions $\# 1$ and $\# 2$ respectively. To apply the Residuals-RAE method in solving the equation, we define the loss function as follows (using AC-I1 as an example):
\begin{align}
\nonumber
\mathcal{L}(\theta) &=  \gamma_{int} \cdot \mathcal{L}_{int}(\theta) +  \gamma_{tb} \cdot \mathcal{L}_{tb}(\theta)
+\gamma_{sb} \cdot \mathcal{L}_{sb}(\theta)\\
\nonumber
&= \gamma_{int} \cdot \frac{1}{M_{int}} \sum_{i=1}^{M_{int}} \hat{\lambda}_{int}^{(i)} | u_{t}(x_{int}^{(i)}, t_{int}^{(i)})+0.001 u_{xx}(x_{int}^{(i)}, t_{int}^{(i)}) + 3u^3(x_{int}^{(i)}, t_{int}^{(i)})-3u(x_{int}^{(i)}, t_{int}^{(i)})| ^2\\
\nonumber
&\quad+  \gamma_{tb} \cdot \frac{1}{M_{tb}} \sum_{i=1}^{M_{tb}} |u(x_{tb}^{(i)}, 0) - u_{0}(x_{tb}^{(i)})| ^2\\
\nonumber
&\quad+  \gamma_{sb} \cdot \left( \frac{1}{M_{sb}} \sum_{i=1}^{M_{sb}} (|u(-x_{sb}^{(i)}, t_{sb}^{(i)}) - u(x_{sb}^{(i)}, t_{sb}^{(i)})| ^2 \right)^{1/2}\\
\nonumber
&\quad + \gamma_{sb} \cdot \left(\frac{1}{M_{sb}} \sum_{i=1}^{M_{sb}} |\nabla u(-x_{sb}^{(i)}, t_{sb}^{(i)}) - \nabla u(x_{sb}^{(i)}, t_{sb}^{(i)})|^2\right)^{1/2},
\end{align}
where the weights for different terms are $\gamma_{int}=1$, $\gamma_{sb}=1$ and $\gamma_{tb}=100$. When the point-wise weights $\hat{\lambda}_{int}^{(i)}$ of interior points are set to 1, the loss function reduces to the classical MSE loss. During training, a fully-connected neural network architecture with a depth of 2 and a width of 128 is used. To optimize the modified loss function, a combination of $30,000$ iterations using the Adam algorithm and $1,000$ iterations using the L-BFGS algorithm are utilized. Alternatively, the training may terminate prematurely if the loss decreases by less than $10^{-7}$ between two consecutive epochs. We randomly sample {$M_{int} = 10,000$, $M_{sb} = 256$, $M_{tb} = 512$} collation points from the computational domain using the latin hypercube sampling (LHS) approach. 

\begin{figure}[H]
\centering
\subcaptionbox{Results for 1D Allen--Cahn equation (with initial condition $\# 1$ ($u(0,x) = x^2 sin(2 \pi x)$)  using Residuals-RAE PINNs.\label{fig::1dAC_rae_i1}}{\includegraphics[width=1.0\linewidth]{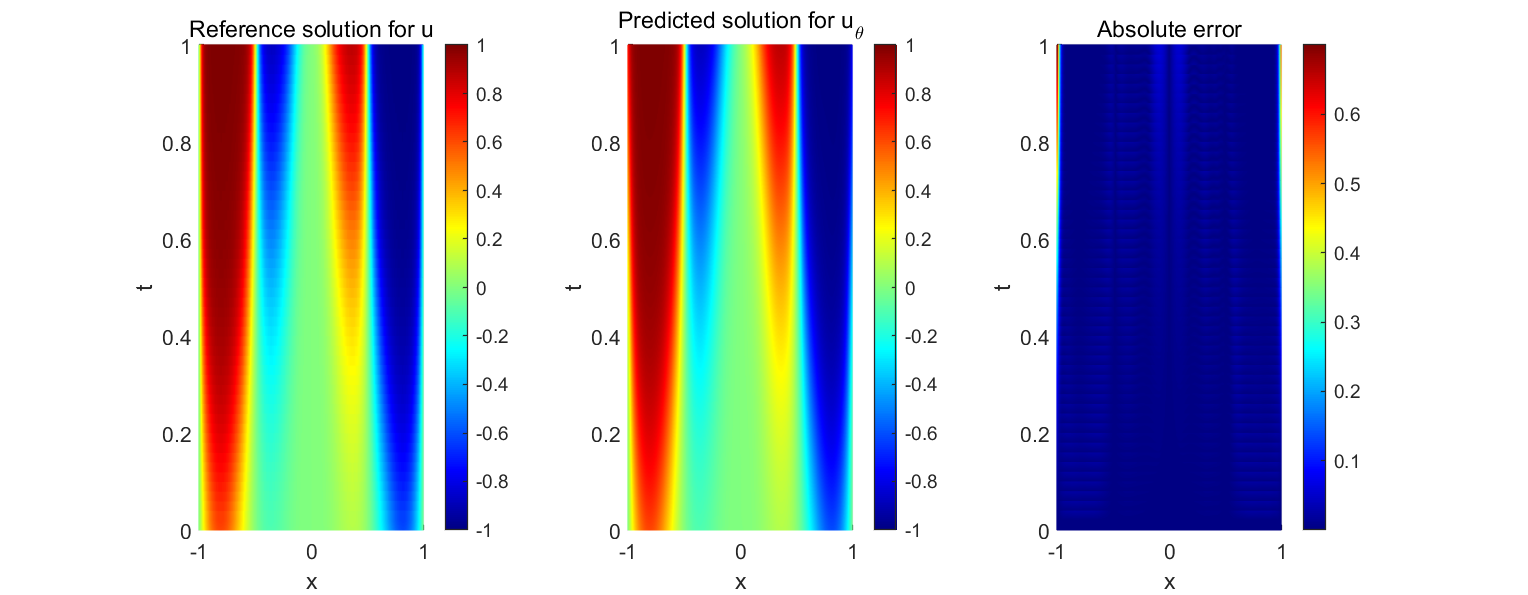}}
\centering
\subcaptionbox{Results for 1D Allen--Cahn equation (with initial condition $\# 2$ $u(0,x) = u(x,0)= cos(\pi x) - exp(-4(\pi x)^2$) using Residuals-RAE PINNs.\label{fig::1dAC_DASA} \label{fig::1dAC_rae_i2}}{\includegraphics[width=1.0\linewidth]{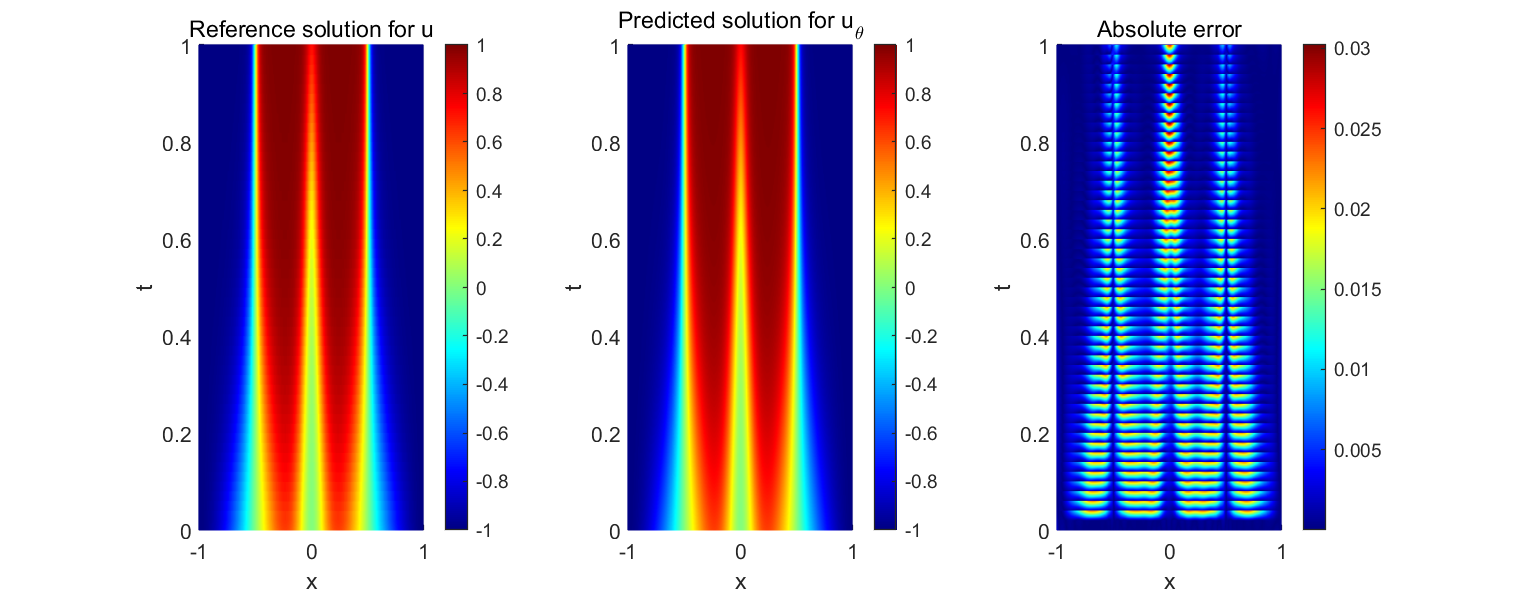}}
\caption{\textbf{Results for solving 1D Allen--Cahn equation using Residuals-RAE.} $1^{st}$ column: Reference $u$, $2^{nd}$ column: the predicted $u$ from Residuals-RAE, $3^{rd}$ column: absolute pointwise error (relative $L^2$ error: $8.092e-03$ (initial condition $\# 1$), and $6.467e-03$ (initial condition $\# 2$).\label{fig::ac1d_result}}
\end{figure}

To implement the Residuals-RAE weighting scheme for pointwise weights $\hat{\lambda}_{int}^{(i)}$ of residual points, the hyper-parameter $k_{int}$ in Eq. (\ref{eq::lknear}) is set to 50. This hyper-parameter determines the influence from the nearest points when measuring the pointwise weights. The relative $L^2$ errors of Residuals-RAE in the Allen--Cahn equation with initial conditions $\#1$ and $\#2$ are  $8.092e-03$, $7.751e-03$, respectively (see Fig. \ref{fig::ac1d_result}). The reference data are obtained using spectral methods with the Chebfun package. Moreover, the point-wise absolute error distributions in the domain tend to accumulate in sub-regions where the PDE solution exhibits sharp transitions. The phenomenon is often observed in physics-informed neural networks (PINNs).

\begin{figure}[H]
\centering
\subcaptionbox{Reference and Residuals-RAE-PINN predicted solutions (neurons per layer $\#1$ and  $\#2$) of the 1D Allen--Cahn equation (initial condition $\#$ 1) at different time snapshots (a) $t=0$, (b) $t=0.5$, (c) $t=1$.	\label{fig::1dACi1_difftime}}{\includegraphics[width=1.0\linewidth, height=0.55\linewidth]{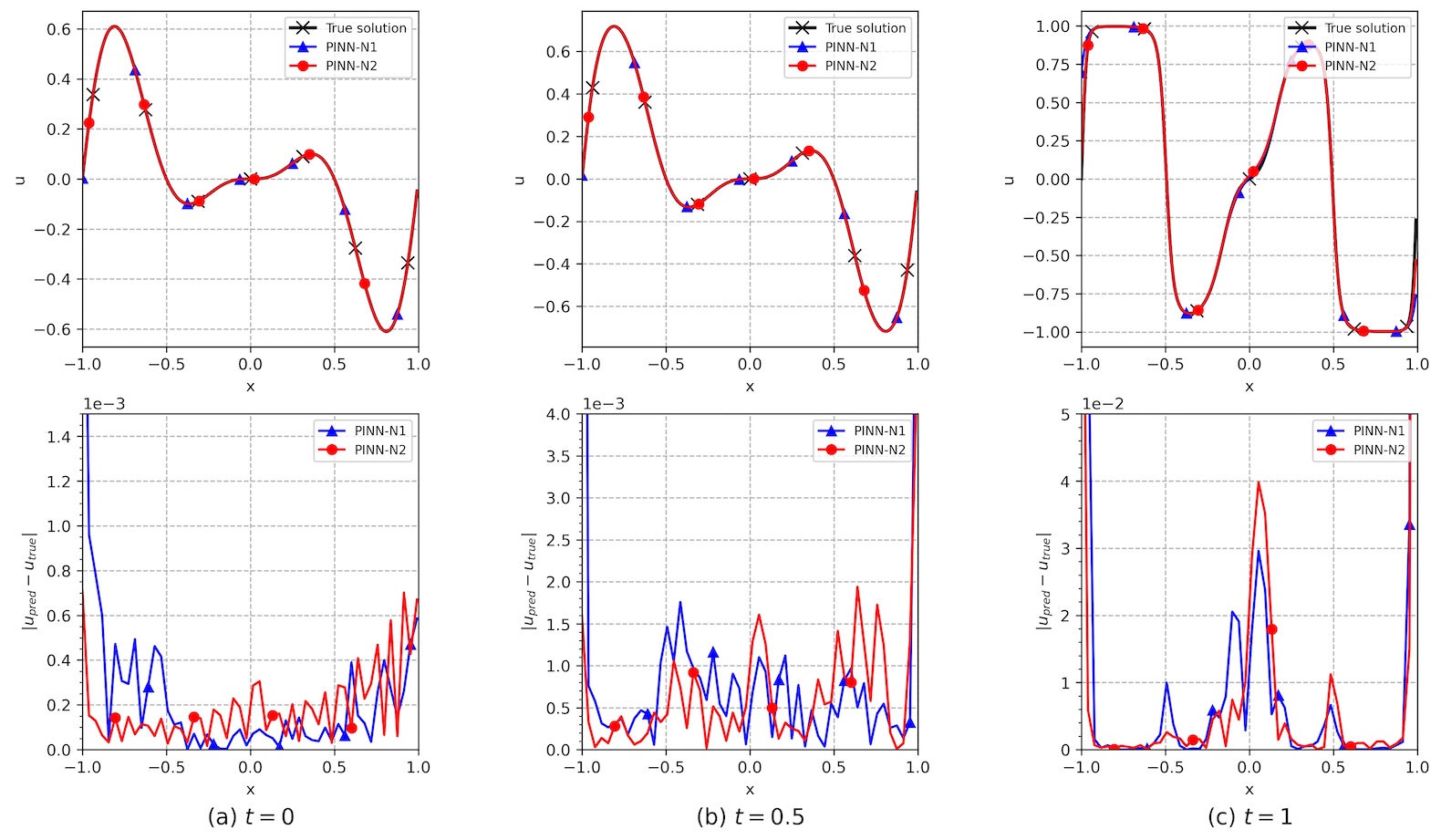}}
\end{figure}

\begin{figure}[H]
\centering
\subcaptionbox{Reference and Residuals-RAE-PINN predicted solutions (neurons per layer $\#1$ and  $\#2$) of the 1D Allen--Cahn equation (initial condition $\#$ 2) at different time snapshots (a) $t=0$, (b) $t=0.5$, (c) $t=1$. \label{fig::1dCH_{int}esample} }{\includegraphics[width=0.95\linewidth, height=0.55\linewidth]{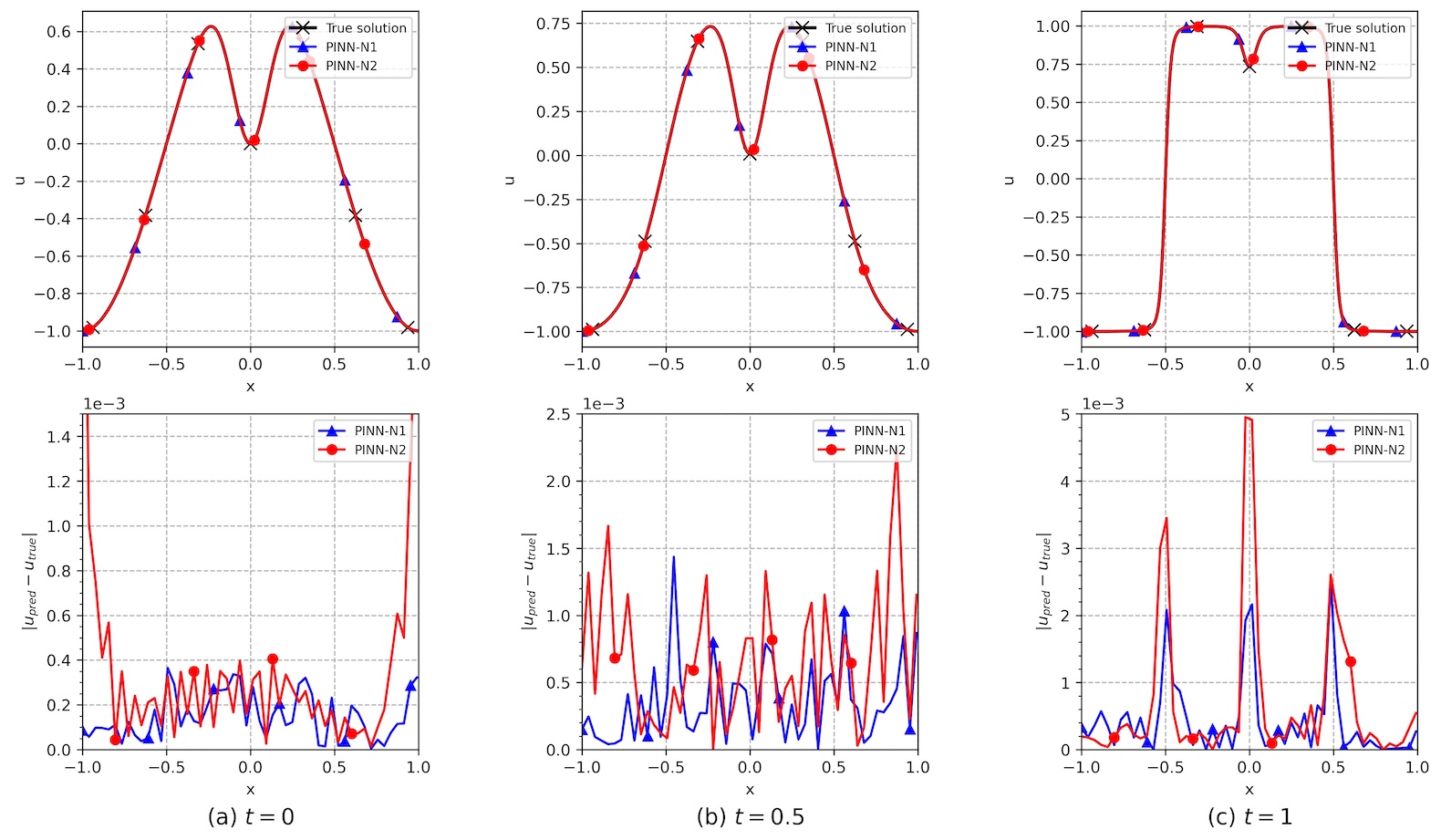}}
\caption{\textbf{Results for solving 1D Allen--Cahn equation using Residuals-RAE.}\label{fig::1dAC_diff}}
\end{figure}

\begin{figure}[H]
  \centering
  \includegraphics[scale=0.42]{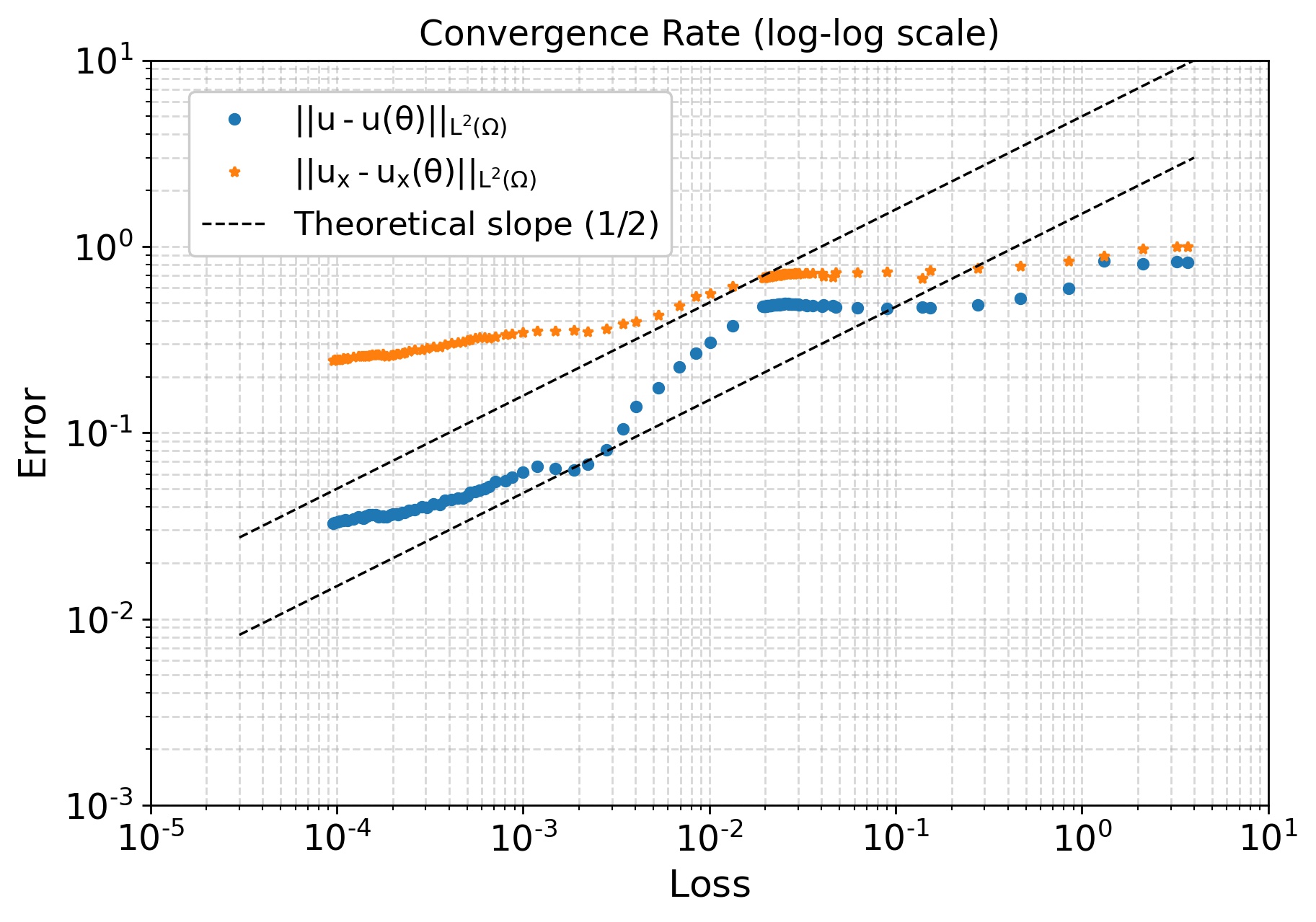}
  \hspace{0.02in}
  \includegraphics[scale=0.42]{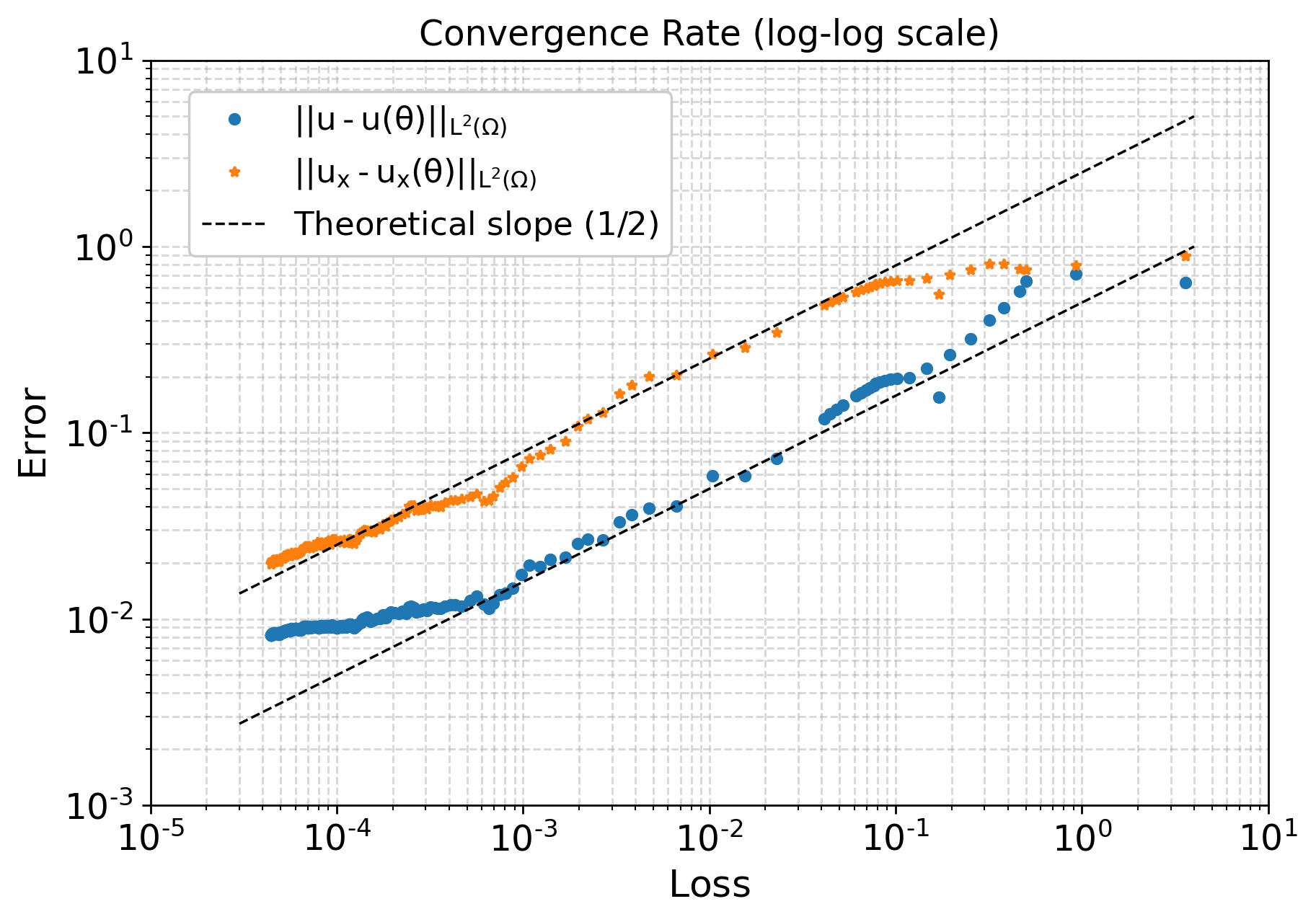}
  \caption{\textbf{1D Allen--Cahn equation: The evolution of $l^2$ solution errors for $u$ and $u_x$ as a function of the training loss value using Residuals-RAE.} Left: 1D Allen--Cahn equation with initial condition $\# 1$. Right: 1D Allen--Cahn equation with initial condition $\# 2$.}\label{fig::acerrordynamics}
\end{figure}


Fig. \ref{fig::1dAC_diff} presents the numerical outcomes achieved through the Residuals-RAE-PINNs method at different time snapshots ($t = 0, 0.5,$ and $1.0$). In order to examine the impact of the width of the neural network, we also provide the numerical results for this problem by varying the number of neurons in the 2-layer neural network. Specifically, PINN-NN1 employs a neural network with 256 neurons per layer, while PINN-NN2 adopts a neural network with 128 neurons per layer. The performance of these networks is evaluated against the reference solution, with the top row of graphs illustrating the solutions at the specified time intervals and the bottom row showcasing the corresponding errors. The numerical findings indicate that the wider neural network (PINN-NN1) performs better, displaying a closer fit to the true solution, and smaller errors compared to the narrower network (PINN-NN2). This is a common outcome in machine learning, where increasing the capacity of a model (such as increasing the number of neurons) can lead to better approximation of complex functions, up to a point of diminishing returns where overfitting can become an issue.

\begin{table}[h]\caption{{\bf {\color{black}{{Description of training data for 1D Allen--Cahn equation (AC-I1, AC-I2).}}} \label{tab::data_ac1d}}}
$$
\begin{array}{lll}
\hline M_{tb} & \text { Initial collocation points } & 512 \\
M_{sb} & \text { Boundary collocation points } & 256 \\
M_{int} & \text { Residual collocation points } & 10,000  \\
\hline
\end{array}
$$
\end{table}

In Fig. \ref{fig::acerrordynamics}, we analyze the behavior of errors for $u$,  $\frac{\partial u}{\partial x}$ in order to assess the convergence rate during training. The errors are plotted against the corresponding values of the loss function used in training, with a logarithmic scale applied. In the log-log plot, lines with slope $1/2$ are plotted as references to indicate the expected convergence rate. Based on the figure, the error dynamics of PINNs exhibit a convergence trend that is broadly consistent with the theoretical square-root scaling with respect to the training loss, although some deviations and fluctuations are observed, particularly for initial condition $\#1$. It is worth noting that in some regions—particularly in the upper-right corner of the first plot in Fig. \ref{fig::acerrordynamics}—the solution error is already relatively small despite a comparatively large training loss. This phenomenon arises because we apply Xavier initialization to the neural network parameters, which places the model in a favorable parameter space at the beginning of training. As training progresses, the loss gradually decreases and the error correspondingly reduces, resulting in an overall trend that increasingly aligns with the reference line. Our numerical results are consistent with our theoretical analysis. Theorem \ref{theorem::ac5} justifies that the total error is bounded by the training error, and the training dynamics in Fig. \ref{fig::acerrordynamics} empirically demonstrate that our proposed method effectively reduces this training error, leading to a corresponding and stable reduction in the true solution error. This validates the effectiveness of our strategy.

{
To validate the accuracy of the proposed Residuals-RAE-PINNs method, we first compared it with the reference solution obtained using the Chebfun method (Fig.  \ref{fig::ac1d_result}). Furthermore, we conducted a fair comparison with baseline PINNs \cite{raissi2019physics}, SA-PINNs \cite{mcclenny2023self}, and PINNs-WE \cite{liu2024discontinuity}. Both SAPINNs and PINNs-WE  aim to improve network training by adaptively adjusting weights, but they employ fundamentally different strategies, and both have demonstrated strong performance. SAPINNs assigns a trainable weight to each point and updates these weights via a gradient ascent approach, encouraging the network to focus on regions that are difficult to learn. This method has shown excellent results in various problems. On the other hand, PINNs-WE takes a complementary approach by reducing the weights in regions where the solution exhibits discontinuities or sharp transitions. This indirect strategy helps the network better capture such features. In the subsequent comparison involving 2D scalar equations, we adopt the magnitude of the gradient to quantify local variations, and compute the point-wise weights in the PINNs-WE method using the formula $\lambda_{1} = \frac{1}{k_{1} | \nabla u | + 1}$. { To ensure reproducibility and a strictly fair comparison, we established a unified set of hyperparameters for all baseline methods (SA-PINNs \cite{mcclenny2023self} and PINNs-WE \cite{liu2024discontinuity}). These settings were optimized based on the best practices reported in their respective literature and were consistently applied across all experiments in this work. The complete list of these unified parameters, including the specific mask learning rates and update frequencies, is detailed in Table \ref{tab:hyperparameters_detailed} in the Appendix \ref{sec:Detailed Experimental Settings}.}

As shown in Table \ref{tab::comp_ac1d}, the PINNs methods with adaptive weighting achieve better accuracy compared to the baseline PINNs, as measured by the relative $l^2$-error, especially under the equation with Initial Condition $\#1$. Both SA-PINNs and PINNs-WE exhibit improved capability in handling local sharp transitions. In addition, the results of Residuals-RAE-PINNs demonstrate enhanced robustness to different initial conditions, as reflected in its consistently low error across both cases. This indicates that the proposed method is less sensitive to variations in initial conditions and possesses stronger generalization ability.
}

\begin{table}[htbp]
\centering
\caption{{Comparative analysis of $l^2$- and $l^\infty$-errors for the 1D Allen-Cahn equation: Performance evaluation of PINN variants under initial conditions $\#1$ and $\#2$.} \label{tab::comp_ac1d}}
\begin{tabular}{ccccc}
\toprule
\textnormal{PINN method}
& \multicolumn{2}{@{\hskip 0pt}l}{\hspace{0.8em}$l_2$-error} 
& \multicolumn{2}{@{\hskip 0pt}l}{\hspace{0.8em}$l_\infty$-error} \\
\cmidrule(lr){2-3} \cmidrule(lr){4-5} 
& \multicolumn{1}{@{}l}{\hspace{0.8em}PINN-AC-I1 $u_\theta$} & \multicolumn{1}{@{}l}{\hspace{0.8em}PINN-AC-I2 $u_\theta$} 
& \multicolumn{1}{@{}l}{\hspace{0.8em}PINN-AC-I1 $u_\theta$} & \multicolumn{1}{@{}l}{\hspace{0.8em}PINN-AC-I2 $u_\theta$} \\
\midrule
baseline \cite{raissi2019physics}   & 1.255e-01 & 4.746e-02 & 4.355e-02 & 1.071e-02 \\
SA-PINNs \cite{mcclenny2023self}   & 6.535e-02 & 4.804e-02 & 3.623e-02 & 1.013e-02 \\
PINNs-WE \cite{liu2024discontinuity}  & 5.603e-02 & 4.582e-02 & 3.300e-02 & 1.091e-02 \\
{\bf Residuals-RAE-PINNs}  & 8.092e-03 & 7.751e-03 & 6.467e-03 & 1.002e-03 \\
\bottomrule
\end{tabular}
\end{table}

\subsubsection{2D Allen--Cahn equation}
Here, we also consider the 2D Allen--Cahn equation defined on the computational domain as $D \times [0,T] = [0, 1]^2 \times [0, 5]$, and the boundary conditions have been considered to be periodic. Then we can have the following formulation 
{\begin{subequations}
\begin{align}
&u_t(x,y,t)-\lambda \varepsilon^2 \Delta u(x,y,t)+\lambda u(x,y,t)^3-\lambda u(x,y,t)=0,\qquad (x,y,t)\in D\times [0,T] \\
&u(x,y,0)=u_0,\qquad (x,y,t)\in D\times \{0\} \\
&u^{(d-1)}(x,y,t)=u^{(d-1)}(-x,-y,t),\qquad (x,y,t)\in \partial D\times[0,T]\qquad(d=1,2)
\end{align}
\end{subequations}}
where the physical parameters $\lambda = 10, \varepsilon = 0.025$, and the initial condition for $u$ was given by $u_0=\tanh \left(\frac{0.35-\sqrt{(x-0.5)^2+(y-0.5)^2}}{2 \varepsilon}\right)$. In this problem, we consider the modified loss function with self-adaptive weights used here, which could be specified as below
\begin{align}
\nonumber
\mathcal{L}(\theta) &=  \gamma_{int} \cdot \mathcal{L}_{int}(\theta)  + \gamma_{tb} \cdot \mathcal{L}_{tb}(\theta)+ \gamma_{sb} \cdot \mathcal{L}_{sb}(\theta) \\
\nonumber
&= \gamma_{int} \cdot \frac{1}{M_{int}} \sum_{i=1}^{M_{int}} \hat{\lambda}_{int}^{(i)} \left| u_t(x_{int}^{(i)}, y_{int}^{(i)}, t_{int}^{(i)})-\lambda \varepsilon^2 \Delta u(x_{int}^{(i)}, y_{int}^{(i)}, t_{int}^{(i)})+\lambda u(x_{int}^{(i)}, y_{int}^{(i)}, t_{int}^{(i)})^3-\lambda u(x_{int}^{(i)}, y_{int}^{(i)}, t_{int}^{(i)})\right| ^2\\
\nonumber
&\quad+\gamma_{tb} \cdot \frac{1}{M_{tb}} \sum_{i=1}^{M_{tb}} \left|u(x_{tb}^{(i)}, y_{tb}^{(i)}, 0) - u_{0}(x_{tb}^{(i)}, y_{tb}^{(i)}, 0) \right| ^2\\
\nonumber
&\quad+\left( \gamma_{sb} \cdot \frac{1}{M_{sb}} \sum_{i=1}^{M_{sb}} (\left|u^{(1)} (x_{sb}^{(i)}, y_{sb}^{(i)}, t_{sb}^{(i)})-u^{(1)}(-x_{sb}^{(i)}, -y_{sb}^{(i)}, t_{sb}^{(i)})\right|^2 \right)^{1/2}\\
\nonumber
&\quad+\left( \gamma_{sb} \cdot \frac{1}{M_{sb}}\sum_{i=1}^{M_{sb}}\left|u^{(0)} (x_{sb}^{(i)}, y_{sb}^{(i)}, t_{sb}^{(i)})-u^{(0)}(-x_{sb}^{(i)}, -y_{sb}^{(i)}, t_{sb}^{(i)}) \right|^2\right)^{1/2}.
\end{align}

Here, the weights for the loss terms, namely $\gamma_{int}, \gamma_{sb}, \gamma_{tb}$, are set to $1, 1, 100$, respectively. To obtain point-wise weights, the residuals-RAE weighting scheme is employed. When considering the influence from nearby points during training, the number of k-nearest points is set to 50 for interior points. The input layer of the PINN structure consists of 3 neurons. For solving the PDEs, a neural network with a depth of 2 and a width of 128 is utilized, { mirroring the configuration used in the 1D experiments.} In addition, we sample a total of $25,600$ points from the interior domain, $512$ points from the boundary, and $1,024$ points from the initial conditions to train the network. The training process consists of $100,000$ ADAM iterations, with an initial learning rate of $0.001$, { incorporating an exponential decay rate of $0.9$ every $5,000$ iterations to ensure stability. Crucially, consistent with the 1D experiments, the baseline methods in this 2D case were strictly implemented using the same unified hyperparameter settings (e.g., specific mask learning rates for SA-PINNs) as detailed in Table \ref{tab:hyperparameters_detailed} (see Appendix \ref{sec:Detailed Experimental Settings}), ensuring a rigorous and fair comparison.}

\begin{table}[h]\caption{{\bf {\color{black}{{Description of training data for 2D Allen--Cahn equation.}}} \label{tab::data_ac2d}}}
$$
\begin{array}{lll}
\hline M_{tb} & \text { Initial collocation points } & 1024 \\
M_{sb} & \text { Boundary collocation points } & 512 \\
M_{int} & \text { Residual collocation points } & 25,600 \\
\hline
\end{array}
$$
\end{table}

\begin{figure}[H]
    \centering
    \subcaptionbox{Results for 2D Allen--Cahn equation ($t=0$) using Residuals-RAE PINNs.}{\includegraphics[width=1.0\linewidth]{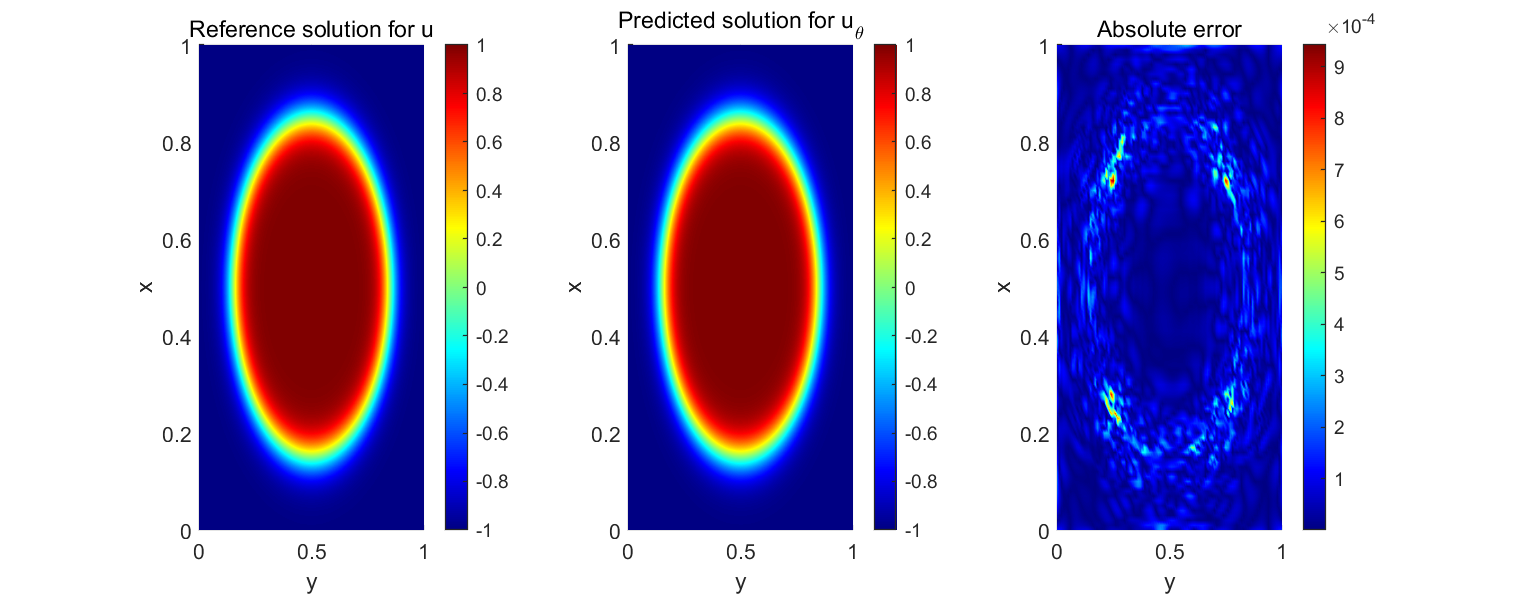}}
    
    \subcaptionbox{Results for 2D Allen--Cahn ($t=2.5$) equation using Residuals-RAE PINNs.}{\includegraphics[width=1.0\linewidth]{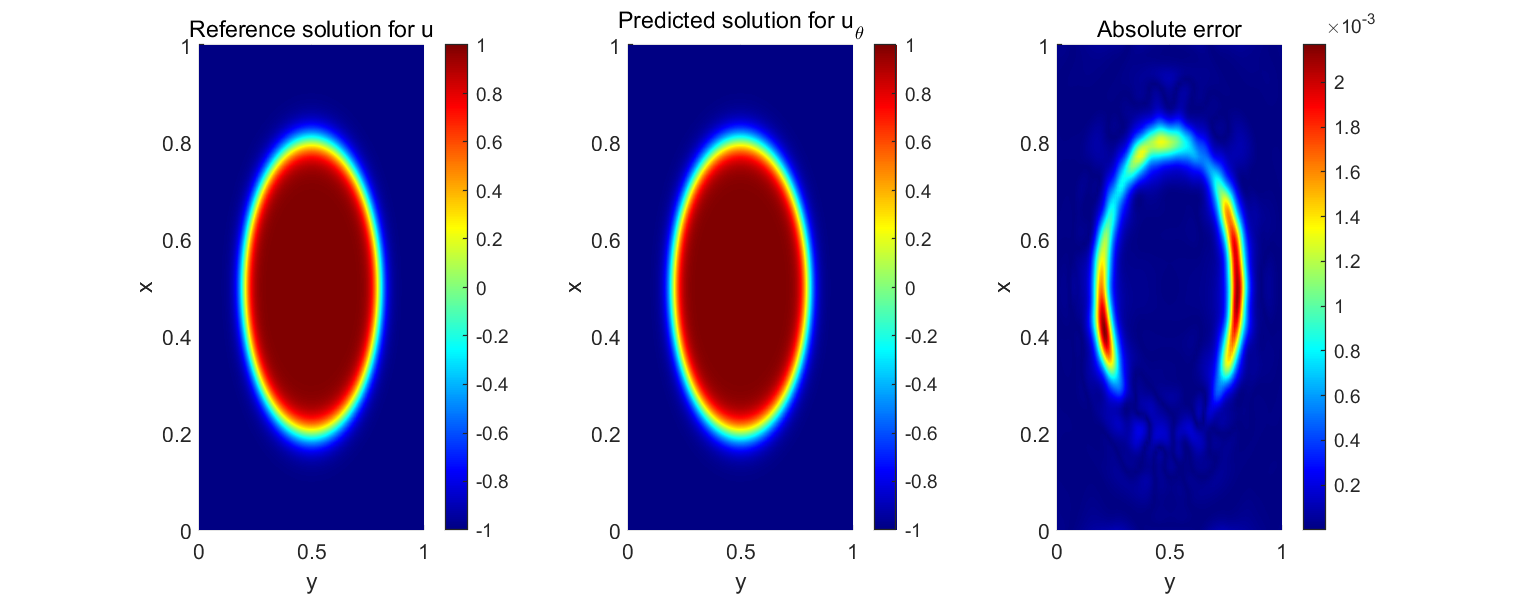}}
    \subcaptionbox{Results for 2D Allen--Cahn equation ($t=5.0$) using Residuals-RAE PINNs.}{\includegraphics[width=1.0\linewidth]{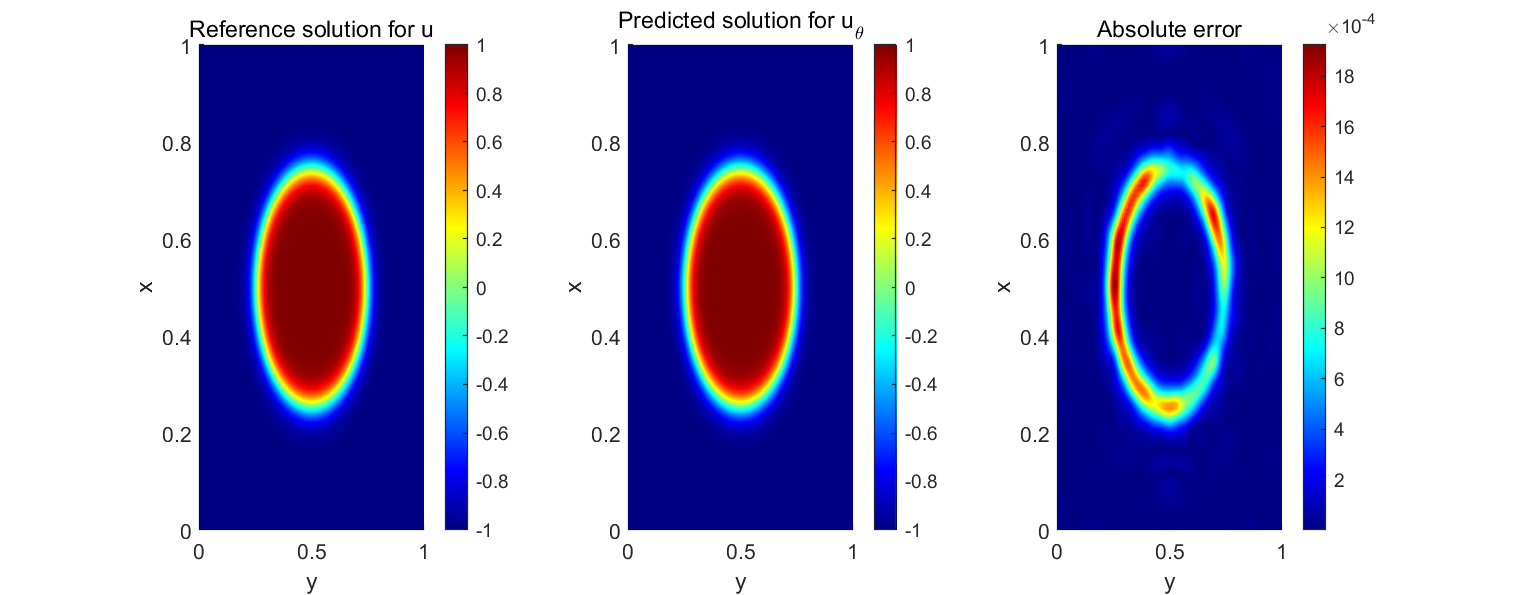}}
    
    \caption{\textbf{Results for solving 2D Allen--Cahn equation using Residuals-RAE.} $1^{st}$ column: Reference $u$, $2^{nd}$ column: the predicted $u$ from Residuals-RAE, $3^{rd}$ column: absolute pointwise error (relative $L^2$ error: $6.42e-04$)  at different time snapshots (a) t=0, (b) t=2.5, (c) t=5.}
    \label{fig::2dac_difft}
\end{figure}

To assess the accuracy of the PINN method, we also present the evolution of the function $u(x,y,t)$ over time by visualizing snapshots at t = 0, 2.5, and 5 for the 2D Allen--Cahn equations (see Fig. \ref{fig::2dac_difft}). The residuals-RAE PINN successfully captures the solution and shows a close match with the reference solution obtained from the spectral method in the Chebfun package. It is worth noting that the accurate characterization of the weights is crucial for the PINN algorithm to produce an accurate approximation.

\begin{figure}[H]
\centering
    \centering
    \includegraphics[width=0.6\linewidth, height = 0.40\linewidth]{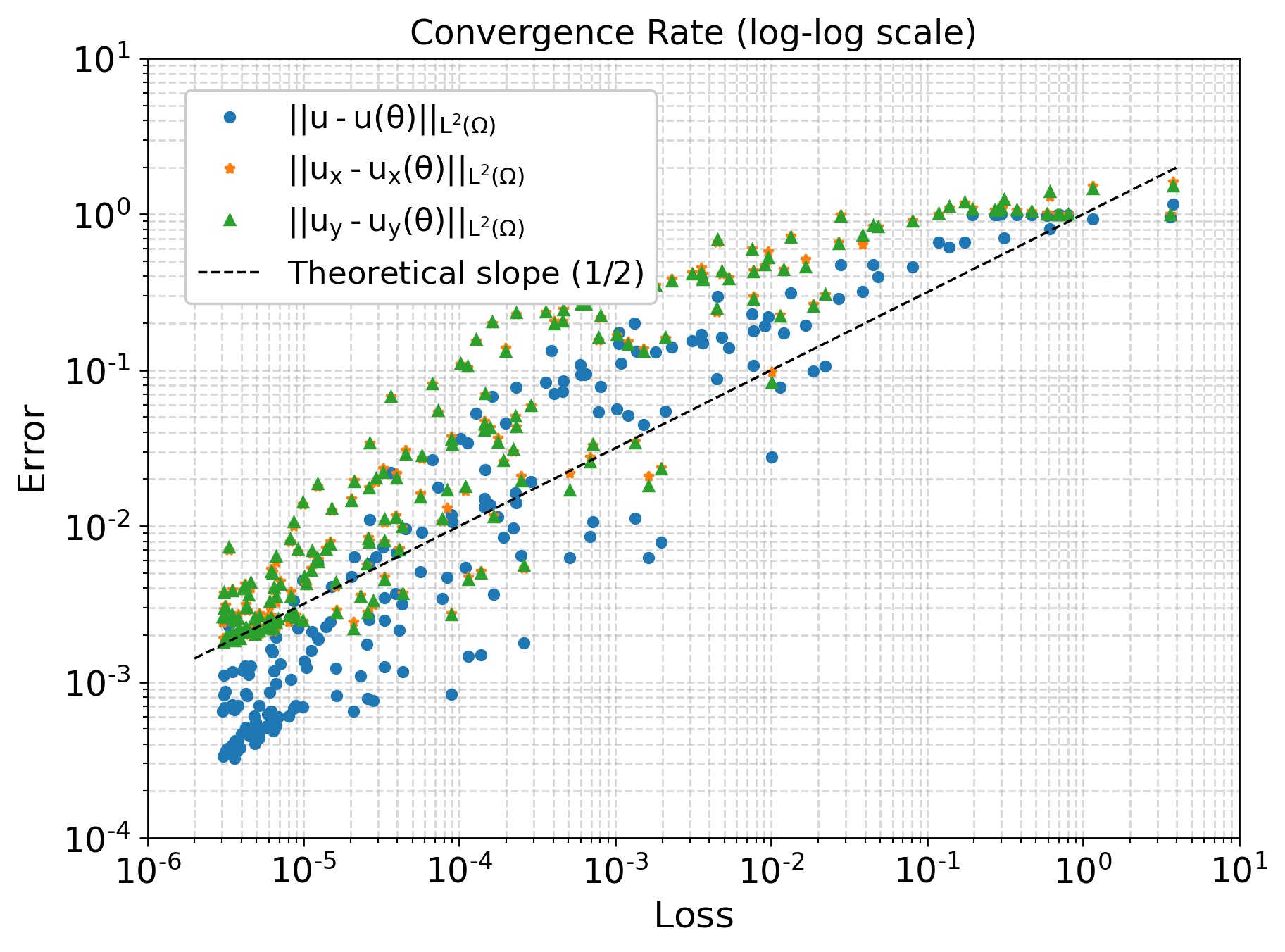}
    \caption{\textbf{2D Allen--Cahn equation: The evolution of $l^2$ solution errors for $u$, $u_x$ and $u_y$ as a function of the training loss value using Residuals-RAE.}\label{fig::2dAC_loss}}
\end{figure}

{Fig. \ref{fig::2dAC_loss} shows the evolution of errors in $u$, $u_x$, and $u_y$ with respect to the training loss during the optimization process. The plot includes a reference line with slope $1/2$, representing the theoretical square-root scaling between the training loss and the solution error. While the convergence curves do not perfectly follow the reference line, they exhibit an overall trend that is broadly consistent with this theoretical behavior. Such plots are commonly used to evaluate the convergence efficiency of neural networks in solving PDEs, particularly when the accuracy of both the solution and its derivatives is of interest. These observations provide numerical evidence supporting the theoretical estimate given in Theorem~\ref{theorem::ac5}.}

{We also present a comparison of relative errors in Table \ref{tab::comp_ac2d}. It can be observed that even shallow neural networks with only two layers are not lacking in expressivity. This highlights the importance of appropriate weight adjustment mechanisms, which help prevent the network from getting trapped in local minima when PDEs are imposed as soft constraints. From the numerical results, we find that both PINNs-WE and Residuals-RAE-PINNs are capable of effectively capturing sharp transitions in the solution, demonstrating their robustness in handling complex solution structures.}

\begin{table}[htbp]
\centering
\caption{{ Comparative analysis of $l^2$- and $l^\infty$-errors for the 2D Allen--Cahn equation: Performance evaluation of PINN variants.} \label{tab::comp_ac2d}}
\begin{tabular}{ccccc}
\toprule
\textnormal{PINN method}
& \multicolumn{1}{@{\hskip 0pt}l}{\hspace{0.8em}Relative $l_2$-error } 
& \multicolumn{2}{@{\hskip 0pt}l}{\hspace{0.8em}Relative $l_\infty$-error} \\
\midrule
baseline \cite{raissi2019physics} & 8.162e-01 & \hspace{1.2em} 7.908e-01 \\
SA-PINNs \cite{mcclenny2023self}  & 9.947e-01 & \hspace{1.2em}  1.100e+00 \\
PINNs-WE \cite{liu2024discontinuity}   & 2.585e-03 & \hspace{1.2em} 1.135e-03 \\
{\bf Residuals-RAE-PINNs}   & 6.423e-04 & \hspace{1.2em} 2.962e-04  \\
\bottomrule
\end{tabular}
\end{table}

\section{Physics Informed Neural Networks for Approximating the Cahn--Hilliard Equation}\label{section::4}

In this section we present the error analysis for the Cahn--Hilliard equation: we state the main error bounds (Theorems~\ref{theorem::ch7} and \ref{theorem::CH2_TER_Numerical}) and their proofs, and then report numerical examples.

\subsection{\textbf{Cahn--Hilliard Equation}}
We first review the general form of Cahn--Hilliard equation in the domain  $\Omega:= D\times[0,T]$, where $D\subset\mathbb{R}^d$, $d\in\mathbb{N}$ and $\partial D$ denotes the boundary of the domain $D$. The equation can be represented as follows:
\begin{subequations}
\label{eq::CH}
    \begin{align}
&u_t-\nabla^2\left(\kappa f(u)-(\alpha \kappa) \nabla^2 u(\boldsymbol{x}, t)\right)=0, \quad t \in[0, T], \quad \boldsymbol{x} \in D \\
& u(\boldsymbol{x}, 0)=\psi(\boldsymbol{x}), \quad \boldsymbol{x} \in D  \\
    & u(\boldsymbol{x}, t)=u(-\boldsymbol{x}, t), \quad t \in[0, T], \quad \boldsymbol{x} \in \partial D \\ 
    & \nabla u(\boldsymbol{x}, t)=\nabla u(-\boldsymbol{x}, t) , \quad t \in[0, T], \quad \boldsymbol{x} \in \partial D
\end{align}
\end{subequations}
where  the  coefficients $\kappa$ and $\alpha$ are non-negative constants. The nonlinear term $f(u)$ is a polynomial of degree 3  with  positive dominant coefficient. $\psi$ as the initial distribution for $u$ belongs to $L^q(D)\cap C(D)$ for $q\geq4$. The boundary conditions are assumed to be periodic.

{
\begin{lemma}[\cite{cardon2001cahn}]\label{chle21}
Let $r \geq 2$, $\psi \in H^r(D)$, and assume $f \in C^1(\mathbb{R})$ with $|f'(s)| \leq C(1 + |s|^{p-1})$ for some $p \geq 1$ and $C > 0$. Let $\kappa > 0$ and $\alpha > 0$. Then there exists $T > 0$ and a unique solution $u$ to the Cahn-Hilliard equation such that $u \in C([0,T]; H^r(D)) \cap C^1([0,T]; H^{r-4}(D))$.
\end{lemma}

\begin{lemma}\label{chle2}
Let $k \in \mathbb{N}$, $\psi \in H^{r+4k}(D)$ with $r > \frac{d}{2} + 4k$, and assume $f \in C^k(\mathbb{R})$ with bounded derivatives up to order $k$. If $f^{(j)}$ has polynomial growth of degree at most $p-j$ for $j = 0, 1, \ldots, k$ and some $p \geq 1$, then there exists $T > 0$ and a solution $u$ to the Cahn-Hilliard equation such that $u(\cdot, 0) = \psi$, $u \in H^{r+4k}(D \times [0,T])$, and $u_t \in H^{r+4k-4}(D \times [0,T])$.
\end{lemma}

The proof for Lemma \ref{chle2} is provided in the Appendix \ref{Proofs}.}

\subsection{\textbf{Physics Informed Neural Networks}}
In order to solve the equation mentioned above, we can utilize a neural network denoted as $u_\theta:\Omega\rightarrow\mathbb{R}$, where $\Omega := D\times[0,T]$ and is parameterized by $\theta$. Then, we can have the residuals: 
\begin{subequations}
\label{eq::CH_{int}esidual1}
\begin{align}
& R_{\text {int }}^{CH}\left[u_\theta\right](\boldsymbol{x}, t)=\partial_t u_\theta-\nabla^2(-(\alpha \kappa) \nabla^2u_\theta+\kappa f(u_\theta)), \\
& R_{t b}^{CH}\left[u_\theta\right](\boldsymbol{x})=u_\theta(\boldsymbol{x}, 0)-\psi(\boldsymbol{x}), \\
& R_{s b 1}^{CH}\left[u_\theta\right](\boldsymbol{x}, t)=u_\theta(\boldsymbol{x}, t)-u_\theta(-\boldsymbol{x}, t), \\
& R_{s b 2}^{CH}\left[u_\theta\right](\boldsymbol{x}, t)=\nabla u_\theta(\boldsymbol{x}, t)-\nabla u_\theta(-\boldsymbol{x}, t).
\end{align}
\end{subequations}
The generalization error, denoted as $\mathcal{E}_G^{CH}(\theta)^2$, is defined as follows:
\begin{equation}
\label{eq::CHGenError}
    \begin{aligned}
        \mathcal{E}_G^{CH}(\theta)^2
        &= \int_D|R_{tb}^{CH}|^2d\bm{x}
        +\int_0^{T}\int_D|R_{int}^{CH}|^2 d\bm{x} dt\\
        &+2(\alpha\kappa)\hat{M} \left(\left(\int_0^T\int_{\partial{D}}|{R_{sb1}^{CH}}|^2ds(\bm{x})dt\right)^{1/2}+\left(\int_0^T\int_{\partial{D}}|{\Delta R_{sb1}^{CH}}|^2ds(\bm{x})dt\right)^{1/2}\right)\\
        &+2(\alpha\kappa)\hat{M}\left(
        \left(\int_0^T\int_{\partial{D}}|{R_{sb2}^{CH}}|^2 ds(\bm{x})dt\right)^{1/2}+
        \left(\int_0^T\int_{\partial{D}}|{\Delta{R}_{sb2}^{CH}}|^2 ds(\bm{x})dt\right)^{1/2}\right),
    \end{aligned}
\end{equation}
where $\hat{M}=|\partial{D}|^{1/2}\cdot(\|{u}\|_{C^3(\partial{D}\times \left[0,T\right])}+\|{u_\theta}\|_{C^3(\partial{D}\times \left[0,T\right])})$ . 

Let $\mathcal{S}^{CH}=\mathcal{S}^{CH}_{int}\cup\mathcal{S}^{CH}_{tb}\cup\mathcal{S}^{CH}_{sb}$ be the set of training points based on {  Monte-Carlo (LHS) sampling} of definition (\ref{eq::CHGenError}), where $\mathcal{S}^{CH}_{int}=\{(\boldsymbol{x}_{int}^{(i)},t_{int}^{(i)})\}_{i=1}^{M_{int}}$, $\mathcal{S}^{CH}_{tb}=\{(\boldsymbol{x}_{tb}^{(i)})\}_{i=1}^{M_{tb}}$, $\mathcal{S}^{CH}_{sb}=\{(\boldsymbol{x}_{sb}^{(i)},t_{sb}^{(i)})\}_{i=1}^{M_{sb}}$ correspond to control equation, initial, and boundary conditions, respectively. Then we define the training error $\mathcal{E}_T^{CH}(\theta,\mathcal{S})^2$ to estiamte the generalization error as follow,
\begin{equation}\label{eq::ch_Gen_error}
    \begin{aligned}
        \mathcal{E}_T^{CH}(\theta,\mathcal{S})^2
        & = \sum_{i=1}^{M_{tb}}\lambda^{(i)}_{tb}|R_{tb}^{CH}(\boldsymbol{x}_{tb}^{(i)})|^2 + \sum_{i=1}^{M_{int}}\lambda^{(i)}_{int}|R_{int}^{CH}(\boldsymbol{x}_{int}^{(i)},t_{int}^{(i)})|^2\\
        & + 2(\alpha\kappa)\hat{M}\left(\left(\sum_{i=1}^{M_{sb}}\lambda^{(i)}_{sb}|R_{sb1}^{CH}(\boldsymbol{x}_{sb}^{(i)},t_{sb}^{(i)})|^2\right)^{1/2} + \left(\sum_{i=1}^{M_{sb}}\lambda^{(i)}_{sb}|\Delta R_{sb1}^{CH}(\boldsymbol{x}_{sb}^{(i)},t_{sb}^{(i)})|^2\right)^{1/2}\right)\\
        &+ 2(\alpha\kappa)\hat{M}\left(  \left(\sum_{i=1}^{M_{sb}}\lambda^{(i)}_{sb}|R_{sb2}^{CH}(\boldsymbol{x}_{sb}^{(i)},t_{sb}^{(i)})|^2\right)^{1/2} + \left(\sum_{i=1}^{M_{sb}}\lambda^{(i)}_{sb}|\Delta R_{sb2}^{CH}(\boldsymbol{x}_{sb}^{(i)},t_{sb}^{(i)})|^2\right)^{1/2}\right).
    \end{aligned}
\end{equation}

{Similarly, the training error for the CH equation is a weighted sum of residuals with tunable hyperparameters. These serve as practical substitutes for the more complex coefficients, such as $2\alpha \kappa\hat{M}$ , found in the theoretical generalization error in the equation \eqref{eq::CHGenError}. The theoretical coefficients themselves reveal the error's origins: $\alpha$ and $\kappa$ are physical parameters from the governing CH equation. As with the AC equation, the $1/2$ power on the boundary integral terms arises from the  proof derivation. 

Next, we use $\hat{u} = u_{\theta} - u$ to represent the numerical error between the predicted solution from PINNs and the exact solution. 
Then we define the total error of PINNs with the following form with $\alpha $ and $\kappa$ being the coefficient of the Cahn-Hilliard equation:}

\begin{definition}[Total Error for Cahn--Hilliard System]
The total error for the Cahn--Hilliard (CH) system, denoted by $\mathcal{E}^{CH}(\theta)^2$, quantifies both the solution discrepancy and the mismatch in second-order spatial derivatives between the predicted and exact solutions. It is defined as:
\begin{equation}
\label{eq::ch_total}
\mathcal{E}^{CH}(\theta)^2 = \|\hat{u}(\boldsymbol{x}, t)\|^2_{L^2(D \times [0, T])} 
+ { 2\alpha\kappa \|\Delta \hat{u}\|^2_{L^2(D\times[0,T])}},
\end{equation}
where:
\begin{itemize}
    \item $\hat{u} = u_{\theta} - u$ denotes the residual error between the predicted solution $u_\theta$ and the exact solution $u$.
    \item $\alpha > 0$ is a scaling constant, and $\kappa$ is the physical interfacial parameter in the Cahn–Hilliard model.
    \item The first term, $\|\hat{u}(\boldsymbol{x},t)\|^2_{L^2(D\times [0,T])}$, represents the standard $L^2$ error over the entire spatiotemporal domain, quantifying the average discrepancy in solution values.
    \item The second term, { $2\alpha\kappa \|\Delta \hat{u}\|^2_{L^2(D\times[0,T])}$}, measures the $L^2$ error of the Laplacian of the solution over the spatiotemporal domain, highlighting deviations in the predicted curvature and interfacial structure.
\end{itemize}
\end{definition}

{\textbf{Remark 4.1.} The total error for the Cahn-Hilliard equation, $\mathcal{E}^{CH}(\theta)^2$, is a tailored metric that differs significantly from the one for the Allen-Cahn equation. Its key distinction is the inclusion of an integrated error term for the \textbf{Laplacian} ($\Delta\hat{u}$), rather than the gradient ($\nabla\hat{u}$). 

This difference is a direct consequence of the mathematical structure of the Cahn-Hilliard equation itself. As a fourth-order PDE, its analysis naturally involves terms with second derivatives. Therefore, the Laplacian term is essential to properly account for and control the error in these higher-order derivatives. The weighting by the physical parameters $2\alpha\kappa$ is also directly inherited from the governing equation's formulation.}

\subsection{\textbf{Error Analysis}}
{ In this subsection, we derive bounds for the residual of the Cahn-Hilliard equations and establish that the total error is controlled by both generalization and training errors. The analysis is conducted for the converged network with stabilized weights, as justified by Theorem~\ref{thm:weight_convergence}. The pre-training weight update mechanism of Residuals-RAE-PINNs ensures that weights remain constant during each training iteration, making this theoretical framework applicable.}

\subsubsection{Bound on the Residuals}
Lemma \ref{lemma::Error_N} elucidates that a tanh neural network with two hidden layers is capable of reducing $\|\hat{u}\|_{H^2(D \times [0,T])}$ to an arbitrarily small value when the width of networks approaches to infinite, which provides a specific limit for this estimation.  Consequently, using this Lemma, we are able to demonstrate the existence of an upper bound for the residuals of Cahn--Hilliard equations. Subsequently, we present the findings related to the application of the Physics-Informed Neural Network (PINN) to the Cahn--Hilliard equation.

Here, we also want to emphasize that the Cahn--Hilliard equation includes high-order derivatives. To simplify it, a commonly used method is to decompose it into two coupled second-order partial differential equations, as shown below:
\begin{subequations}
\begin{align}
    &\label{eq::CH_decompostion} u_t-\nabla^2(-(\alpha \kappa) \mu+\kappa f(u))=0, \quad \mu=\nabla^2 u \quad t \in[0, T], \quad x \in D \subset \mathbb{R}^d \\ 
    & f(u)=u^3-u.
\end{align}
\end{subequations}

The decomposition presented here is similar to the second-order decomposition employed in the work conducted by Qian et al \cite{qian2023physics}. for general second-order PDEs in the temporal domain. Qian et al. divided the residual loss, which is determined by the interior points of second-order PDEs, into two components. By utilizing the error analysis theorem derived earlier in the above Cahn--Hilliard equation (\ref{eq::CH}), a new approximation for the error in Cahn--Hilliard equations can be derived. 
{ 
\begin{remark}[Error Propagation in the Decomposed Formulation]
\label{rem:error_propagation}
The introduction of the auxiliary variable $\mu = \Delta u$ raises the question of whether the approximation error in $\mu_\theta$ compounds the error in $u_\theta$. To address this, we decompose the auxiliary error $\hat{\mu} = \mu_\theta - \mu$ as follows:
\begin{equation}
\label{eq::error_decomposition_mu}
\hat{\mu} = (\mu_\theta - \Delta u_\theta) + (\Delta u_\theta - \Delta u) = R_{int1} + \Delta \hat{u}.
\end{equation}
This decomposition reveals that the auxiliary error consists of two additive components: the constraint violation $R_{int1}$ and the Laplacian of the primary error $\Delta \hat{u}$. Crucially, these contributions do not multiply; rather, they combine linearly. The generalization error \eqref{eq::CH2GenError} explicitly penalizes $\|R_{int1}\|_{L^2}^2$ through the term $(\alpha\kappa)\int_{0}^T\int_{D}|R_{int1}|^2 d\boldsymbol{x} dt$, while Definition 4.1 controls $\|\Delta \hat{u}\|_{L^2}^2$ directly. As training progresses and $R_{int1} \to 0$, the constraint $\mu_\theta = \Delta u_\theta$ is enforced, yielding $\hat{\mu} \approx \Delta \hat{u}$. Theorem~\ref{theorem::ch7} then guarantees that the total error $\mathcal{E}^{CH}(\theta)^2$, incorporating both $\|\hat{u}\|_{L^2}^2$ and $\|\Delta \hat{u}\|_{L^2}^2$, remains bounded by the generalization error. Hence, the decomposition does not introduce multiplicative error accumulation.
\end{remark}
}
To begin with, we give the definition of residuals as follows,
{ \begin{subequations}
\label{eq::CH_residual_decomposition}
\begin{align}
    &R_{int1}^{CH}=\hat{\mu}-\nabla^2 \hat{u},\\
    &R_{int2}^{CH}=\hat{u}_t+\alpha \kappa\nabla^2 \hat{\mu}-\kappa \nabla^2(f(u_\theta)-f(u)),\\
    &R_{tb}^{CH}=\hat{u}(x,0),\\
    &R_{sb1}^{CH}=\hat{u}(x,t)-\hat{u}(-x,t),\\
    &R_{sb2}^{CH}=\nabla\hat{u}(x,t)-\nabla\hat{u}(-x,t),
\end{align} 
\end{subequations}}
where $\hat{u} = u_\theta - u$, $\hat{\mu} = \mu_\theta - \mu$. We next redefine the generalization error as $\mathcal{E}_G^{CH}(\theta)^2$, which is expressed as follows:
\begin{equation}
\label{eq::CH2GenError}
    \begin{aligned}
        \mathcal{E}_G^{CH}(\theta)^2
        &= { \int_D|R_{tb}^{CH}|^2d\boldsymbol{x}
        +\int_0^{T}\int_D|R_{int2}^{CH}|^2 d\boldsymbol{x} dt
        +(\alpha\kappa)\int_0^{T}\int_D|R_{int1}^{CH}|^2 d\boldsymbol{x} dt}\\
        &{ +2(\alpha\kappa)\hat{M}\left( \left(\int_0^T\int_{\partial D}|R_{sb1}^{CH}|^2ds(\boldsymbol{x})dt\right)^{1/2}+\left(\int_0^T\int_{\partial D}|\Delta R_{sb1}^{CH}|^2ds(\boldsymbol{x})dt\right)^{1/2}\right)}\\
        &{ +2(\alpha\kappa)\hat{M}\left(
        \left(\int_0^T\int_{\partial D}|R_{sb2}^{CH}|^2 ds(\boldsymbol{x})dt\right)^{1/2}+
        \left(\int_0^T\int_{\partial D}|\Delta R_{sb2}^{CH}|^2 ds(\boldsymbol{x})dt\right)^{1/2}\right),}
    \end{aligned}
\end{equation}
where $\hat{M}=\max\{M_u,M_\mu\}$, and $M_u=|\partial{D}|^{1/2}\cdot(\|u\|_{C^1(\partial{D}\times \left[0,T\right])}+\|u_\theta \|_{C^1(\partial{D}\times \left[0,T\right])})$
$M_\mu=|\partial{D}|^{1/2}\cdot(\|\mu\|_{C^1(\partial{D}\times \left[0,T\right])}+\|\mu_\theta \|_{C^1(\partial{D}\times \left[0,T\right])})$. 

Since it is not feasible to directly minimize the generalization error $\mathcal{E}_G^{(CH)}(\theta)^2$ in practice, we approximate the integrals in definition (\ref{eq::CH2GenError}) using {  Monte-Carlo (LHS) sampling}. {  The training error below uses the same residual notation ($R_{tb}^{CH}$, $R_{int1}^{CH}$, $R_{int2}^{CH}$, $R_{sb1}^{CH}$, $R_{sb2}^{CH}$, and $\Delta R_{sb1}^{CH}$, $\Delta R_{sb2}^{CH}$) as in the generalization error (\ref{eq::CH2GenError}); the only difference is that integrals are replaced by weighted sums over sample points.} Therefore, we define the training error as follows,
\begin{equation}
\label{eq::ch2_training_error}
    \begin{aligned}
        \mathcal{E}_T^{CH}(\theta,\mathcal{S})^2
        & = \sum_{i=1}^{M_{tb}}\lambda^{(i)}_{tb}|R_{tb}^{CH}(\boldsymbol{x}_{tb}^{(i)})|^2 + \sum_{i=1}^{M_{int}}\lambda^{(i)}_{int}|R_{int2}^{CH}(\boldsymbol{x}_{int}^{(i)},t_{int}^{(i)})|^2+(\alpha\kappa)\sum_{i=1}^{M_{int}}\lambda^{(i)}_{int}|R_{int1}^{CH}(\boldsymbol{x}_{int}^{(i)},t_{int}^{(i)})|^2\\
        & + 2(\alpha\kappa)\hat{M}\left(\left(\sum_{i=1}^{M_{sb}}\lambda^{(i)}_{sb}|R_{sb1}^{CH}(\boldsymbol{x}_{sb}^{(i)},t_{sb}^{(i)})|^2\right)^{1/2} + \left(\sum_{i=1}^{M_{sb}}\lambda^{(i)}_{sb}|\Delta R_{sb1}^{CH}(\boldsymbol{x}_{sb}^{(i)},t_{sb}^{(i)})|^2\right)^{1/2}\right)\\
        & + 2(\alpha\kappa)\hat{M}\left( \left(\sum_{i=1}^{M_{sb}}\lambda^{(i)}_{sb}|R_{sb2}^{CH}(\boldsymbol{x}_{sb}^{(i)},t_{sb}^{(i)})|^2\right)^{1/2} + \left(\sum_{i=1}^{M_{sb}}\lambda^{(i)}_{sb}|\Delta R_{sb2}^{CH}(\boldsymbol{x}_{sb}^{(i)},t_{sb}^{(i)})|^2\right)^{1/2}\right).
    \end{aligned}
\end{equation}

Here, $\mathcal{S}^{CH}=\mathcal{S}^{CH}_{int}\cup\mathcal{S}^{CH}_{tb}\cup\mathcal{S}^{CH}_{sb}$ represents the collection of training points. Specifically, $\mathcal{S}^{CH}_{int}=\{(\boldsymbol{x}_{int}^{(i)},t_{int}^{(i)})\}_{i=1}^{M_{int}}$ , $\mathcal{S}^{CH}_{tb}=\{(\boldsymbol{x}_{tb}^{(i)})\}_{i=1}^{M_{tb}}$ , $\mathcal{S}^{CH}_{sb}=\{(\boldsymbol{x}_{sb}^{(i)},t_{sb}^{(i)})\}_{i=1}^{M_{sb}}$ correspond to the control equation, initial conditions, and boundary conditions, respectively.

The forthcoming theorem aims to perform an error analysis on the decoupled Cahn--Hilliard (CH) equation. 
This approach entails estimating the residuals and total error associated with the CH equation, with the objective of obtaining conclusions that are consistent with those derived in the previous theorems. 
{\begin{theorem}\label{theorem::ch6}
Let $n\geq 2$, $d, r, k \in \mathbb{N}$ with $k \geq 5$. Let $\psi \in H^{r+4k}(D)$ with $r > \frac{d}{2} + 4k$, and assume $f \in C^k(\mathbb{R})$ with bounded derivatives up to order $k$. If $f^{(j)}$ has polynomial growth of degree at most $p-j$ for $j = 0, 1, \ldots, k$ and some $p \geq 1$. For every integer $N>5$, there are tanh neural networks $u_\theta$ and $\mu_\theta$, each with two hidden layers, where the widths are at most $3\left\lceil\frac{k+n-2}{2}\right\rceil\left|P_{k-1, d+2}\right|+\lceil N T\rceil+d(N-1)$ and $3\left\lceil\frac{d+n+1}{2}\right\rceil\left|P_{d+2, d+2}\right|\lceil N T\rceil N^d$, such that 
    \begin{subequations}
    \label{eq::CH_residualBound1}
    \begin{align}
    &\left\|R_{int1}^{CH}\right\|_{L^2(\Omega)}\lesssim \ln^2 \left(N\right) N^{-k+2},\\
    & \left\|R_{int2}^{CH}\right\|_{L^2(\Omega)}\lesssim \ln^2 (N) N^{-k+4}, \\ 
    & \left\|R_{s b 2}^{CH}\right\|_{L^2(\partial D \times[0, T])}\lesssim \ln^2 (N) N^{-k+2}, \\ 
    & \left\|R_{s b 1}^{CH}\right\|_{L^2(\partial D \times[0, T])},\left\|R_{tb }^{CH}\right\|_{L^2(D)}  \lesssim \ln (N) N^{-k+1} .
    \end{align}
    \end{subequations}
\end{theorem}
}

\textbf{Proof of Theorem \ref{theorem::ch6}.} 
Due to the nonlinear term of the Cahn--Hilliard equation $f(u)$ and Lemma \ref{chle2}, there exists a positive constant $M$ such that $|\nabla^2f(u_1)-\nabla^2f(u_2)|\leq M|u_1-u_2|, \forall u_1,u_2\in C^2\left(\Omega\right)$. Therefore, we have 
\begin{equation}\label{ineq::ch}
\begin{aligned}
    R_{\text {int2 }}^{CH}(\boldsymbol{x}, t)
    & =\partial_t \hat{u} + (\alpha \kappa) \nabla^2 \hat{\mu}+\kappa \nabla^2 (f(u_\theta)-f(u)))\\
    & \leq
    \partial_t \hat{u}+(\alpha \kappa) \nabla^2 \hat{\mu}+\kappa M|\hat{u}|
    \end{aligned}
\end{equation}

Based on Lemma \ref{chle2} and Lemma \ref{lemma::Error_N}, there exists a neural network $u_\theta$ and $\mu_\theta$ with two same hidden layers and widths $3\left\lceil\frac{k+n-2}{2}\right\rceil\left|P_{k-1, d+2}\right|+\lceil N T\rceil+d(N-1)$ and $3\lceil \frac{d+n+1}{2} \rceil \left| P_{d+2,d+2} \right|\left[ NT \right] N^d$, such that for every $0 \leq l \leq 2$ and $0 \leq s \leq 2$,
\begin{equation}
\label{eq::uErrorAboutN2}
\begin{aligned}
    &\left\|u_\theta-u\right\|_{H^l(\Omega)} \leq C_{l, k, d+1, u} \lambda_{l, u}(N) N^{-k+l},\quad 0 \leq l \leq 2\\
    &\left\|\mu_\theta-\mu\right\|_{H^{l}(\Omega)} \leq C_{s, k-2, d+1, u} \lambda_{s, \mu}(N) N^{-k+s+2},\quad 0 \leq s \leq 2
\end{aligned}
\end{equation}
where $\lambda_{l, u}=2^l 3^{d+1}(1+\sigma) \ln ^l\left(\beta_{l, \sigma, d+1, u} N^{d+k+3}\right)$, $\sigma=\frac{1}{100}$, $\lambda_{s, u}=2^{s} 3^{d+1}(1+\sigma) \ln ^{s}\left(\beta_{l, \sigma, d+1, u} N^{d+k+1}\right)$ and the definition for the other constants can be found in Lemma \ref{lemma::Error_N}.

In light of Lemma \ref{lemma::trace},we can bound the terms of the PINN residual as follows:
\begin{subequations}
\begin{align}
	&\left\| \hat{u}_t \right\| _{L^2(\Omega )}\le \parallel \hat{u}\parallel _{H^1(\Omega )},\\
	&\parallel \nabla ^2\hat{u}\parallel _{L^2(\Omega )}\le \parallel \hat{u}\parallel _{H^2(\Omega )},\parallel \nabla ^2\hat{\mu}\parallel _{L^2(\Omega )}\le \parallel \hat{\mu}\parallel _{H^2(\Omega )},\\
	&\parallel \hat{u}\parallel _{L^2(D)}\le \parallel \hat{u}\parallel _{L^2(\partial \Omega )}\le C_{h_{\Omega},d+1,\rho _{\Omega}}\parallel \hat{u}\parallel _{H^1(\Omega )},\\
	&\parallel \nabla \hat{u}\parallel _{L^2(D)}\le \parallel \nabla \hat{u}\parallel _{L^2(\partial \Omega )}\le C_{h_{\Omega},d+1,\rho _{\Omega}}\parallel \hat{u}\parallel _{H^2(\Omega )},\\
	&\parallel \nabla \hat{u}\parallel _{L^2(\partial D\times [0,T])}\le \parallel \nabla \hat{u}\parallel _{L^2(\partial \Omega )}\le C_{h_{\Omega},d+1,\rho _{\Omega}}\parallel \hat{u}\parallel _{H^2(\Omega )}.
\end{align}
\end{subequations}

By combining these relationships with inequalities (\ref{ineq::ch}) and (\ref{eq::uErrorAboutN2}), we can obtain
\begin{equation*}
\begin{aligned}
        \parallel R_{int1}^{CH}\parallel 
        &\leq \parallel \mu \parallel _{L^2(\Omega )}+\parallel \nabla ^2 u \parallel_{L^2(\Omega )}\\
        &=\parallel \mu \parallel_{H^0(\Omega )}+ \parallel\hat{u} \parallel_{H^2(\Omega )}\\
        &\lesssim C_{0.k-2,d+1,\mu}\lambda _{0,\mu}(N)\cdot N^{-k+2}+C_{2.k,d+1,u}\lambda _{2,u}(N)\cdot N^{-k+2}\\
        &\lesssim \ln ^2N\cdot N^{-k+2}
\end{aligned}
\end{equation*}
\begin{equation*}
   \begin{aligned}
	\parallel R_{int2}^{CH}\parallel _{L^2(\Omega )}&\le \parallel \hat{u}_t+(\alpha \kappa )\nabla ^2(\hat{\mu})+M\kappa |\hat{u}|\parallel _{L^2(\Omega )}\\
	&\le \parallel \hat{u}\parallel _{L^2(\Omega )}+|\alpha \kappa |\parallel \nabla ^2\mu \parallel _{L^2\left( \Omega \right)}+M|\kappa |\parallel \hat{u}\parallel _{L^2(\Omega )}\\
	&\le \parallel \hat{u}\parallel _{H^1(\Omega )}+|\alpha \kappa |\parallel \hat{\mu}\parallel _{H^2(\Omega )}+M|\kappa |\parallel \hat{u}\parallel _{H^0(\Omega )}\\
	&\le C_{1,k,d+1,u}\lambda _{1,u}(N)N^{-k+1}+|\alpha \kappa |C_{2,k-2,d+1,\mu}\lambda _{2,\mu}(N)N^{-k+4}\\
    & \quad+M|\kappa |C_{0,k,d+1,u}\lambda _{0,u}(N)N^{-k}\\
	&\lesssim \ln ^2(N)N^{-k+4},\\
\end{aligned}
\end{equation*}
\begin{equation*}
\left\|R_{t b }^{CH}\right\|_{L^2(D)} \leq C_{h_{\Omega}, d+1,
\rho_{\Omega}}\|\hat{u}\|_{H^1(\Omega)} \lesssim \ln( N) N^{-k+1},
\end{equation*}
\begin{equation*}
\left\|R_{s b 1}^{CH}\right\|_{L^2(\partial D \times[0, T])} 
\leq 2\cdot C_{h_{\Omega}, d+1, \rho_{\Omega}}\|\hat{u}\|_{H^1(\Omega)} \lesssim \ln( N) N^{-k+1},
\end{equation*}
\begin{equation*}
\left\|R_{s b 2}^{CH}\right\|_{L^2(\partial D \times[0, T])} 
\leq 2\cdot C_{h_{\Omega}, d+1, \rho_{\Omega}}\|\hat{u}\|_{H^2(\Omega)} \lesssim \ln ^2 (N) N^{-k+2} .
\end{equation*}


{In light of the established bounds equation \eqref{eq::CH_residualBound1}, it becomes evident that selecting a sufficiently large number
of $N$ ($N$ is a theoretical parameter from function approximation theory that controls the network's architectural complexity), we can effectively minimize both the residuals of the Physics-Informed Neural Network (PINN) as
stated in equation \eqref{eq::CH_residual_decomposition}, and reduce the generalization error to extremely small levels. This observation
is crucial in illustrating the capacity of PINNs to adapt and converge effectively in response to the
complexity of the network configuration.

\subsubsection{Bounds on the Total Approximation Error}
We next demonstrate that the total error $\mathcal{E}^{CH}(\theta)^2$ remains sufficient small when the generalization error $\mathcal{E}_G(\theta)^2$ is small using the PINN approximation.}
{
\begin{theorem}\label{theorem::ch7}
Let $d \in \mathbb{N}, u \in C^1(\Omega)$ and $\mu\in C^1(\Omega)$ be the classical solution to the Cahn--Hilliard equations. Let $u_\theta$ and $\mu_\theta$ denote the PINN approximation solutions with parameter $\theta$. Then the following relation holds,
{ 
\begin{equation}
\mathcal{E}^{CH}(\theta)^2\leq C_G\,\frac{\exp((1+2\kappa{M})T)-1}{1+2\kappa{M}}.
\end{equation}
}
where $M$ is a constant and $C_G$ is given by
\begin{equation*}
{  C_G =\int_0^{T}\int_D|R_{int2}|^2 d\boldsymbol{x} dt+\left(\alpha\kappa\right)\int_{0}^T\int_{D}|R_{int1}|^2d\boldsymbol{x} dt+2\left(\alpha\kappa\right)\hat{M}\hat{C}}
\end{equation*}
with the components $\hat{M}$ and $\hat{C}$ defined as
\begin{equation}
\begin{split}
\hat{M} & = \max_{g \in \{u,\mu \} }\{|\partial{D}|^{1/2}\cdot\left(\|g\|_{C^1(\partial{D}\times[0,T])}+\|g_\theta \|_{C^1(\partial{D}\times [0,T])}\right)\},\\
\hat{C} & = \sum_{i\in \{ sb_1, sb_2\}}\left(\int_0^T\int_{\partial{D}}|\Delta{R_{i}}|^2ds(\boldsymbol{x})dt\right)^{1/2}+\sum_{j\in \{ sb_1, sb_2\}}\left(\int_0^T\int_{\partial{D}}|{R_{j}}|^2ds(\boldsymbol{x})dt\right)^{1/2}.
\end{split}
\end{equation}
\end{theorem}}
\textbf{Proof of Theorem \ref{theorem::ch7}.} Considering 
\begin{equation}
\label{eq::CH_int}
\begin{aligned}
\int_D R_{int2}^{CH}\cdot\hat{u} d\bm{x} 
& = \int_D \hat{u}_t\hat{u}d\bm{x}-\kappa\int_D\nabla^2[f(u_\theta)-f(u)]\hat{u}d\bm{x} + (\alpha\kappa)\int_D\nabla^2(\hat{\mu})\hat{u} d\bm{x}\\
& = { \frac{1}{2}\frac{d}{dt}}\|\hat{u}\|_{L^2(D)}^2-\kappa\int_D\nabla^2[f(u_\theta)-f(u)]\hat{u} d\bm{x} + (\alpha\kappa)\int_D\nabla^2(\hat{\mu})\hat{u} d\bm{x}\\
& = { \frac{1}{2}\frac{d}{dt}}\|\hat{u}\|_{L^2(D)}^2-\kappa{I}_1+(\alpha\kappa)I_2,
\end{aligned}
\end{equation}
with
\begin{subequations}
\label{eq::CH_I}
\begin{align}
    I_1
    &=\int_D\nabla^2[f(u_\theta)-f(u)]\hat{u} d\bm{x} \leq  M\int_D|\hat{u}|^2 d\bm{x}\\
    \nonumber
    I_2
    & = \int_D\nabla^2\hat{\mu}\hat{u} dx\\
        \nonumber
    & = -\int_{D}\nabla\hat{\mu}\nabla\hat{u} d\bm{x}+\int_{\partial{D}}\nabla\hat{\mu}\hat{u}\cdot\bm{n}ds(\bm{x})\\
    & =\int_D\hat{\mu}\cdot\Delta\hat{u} d\bm{x} + \int_{\partial{D}}\nabla\hat{\mu}\cdot\hat{u}\cdot\bm{n}ds(\bm{x})-\int_{\partial{D}}\hat{\mu}\nabla{\hat{u}}\cdot\bm{n}ds(\bm{x}),\label{eq::ch2}
\end{align}
\end{subequations}
where $M$ is a Lipschitz constant of $\nabla^2{f}$, and equation (\ref{eq::ch2}) are  established because of the divergence integral theorem.

By combining Eq. (\ref{eq::CH_int}) and (\ref{eq::CH_I}), we have
\begin{equation}
\label{eq::CH_int_dt}
    \begin{aligned}
    &\quad\quad { \frac{1}{2}\frac{d}{dt}}\|\hat{u}\|_{L^2(D)}^2+(\alpha\kappa)\int_D|\hat{\mu}|^2 d\bm{x}\\
    & = \int_D R_{int2}^{CH}\hat{u} d\bm{x}+\kappa{I}_1
    +(\alpha\kappa)\int_{D}R^{CH}_{int1}
    \cdot\hat{\mu} d\bm{x} +(\alpha\kappa)\int_{\partial{D}}\nabla
    \hat{\mu}\cdot\hat{u}\cdot\bm{n}ds(\bm{x})
    -(\alpha\kappa)\int_{\partial{D}}
    \Delta\hat{u}\nabla{\hat{u}}\cdot
    \bm{n}ds(\bm{x})
    \end{aligned}
\end{equation}

Based on Cauchy-Schwarz's inequality, we have 
\begin{subequations}
\label{eq::CH_sb}
    \begin{align}\int_{\partial{D}}\Delta\hat{u}\nabla{\hat{u}}\cdot\bm{n}ds(\bm{x})
    &\leq M_1\left(\int_{\partial{D}}|\Delta{R_{sb1}^{CH}}|^2ds(\bm{x})\right)^{1/2}+M_2\left(\int_{\partial{D}}|{R_{sb2}^{CH}}|^2ds(\bm{x})\right)^{1/2},\\
    \int_{\partial{D}}\nabla \hat{\mu}\cdot\hat{u}\cdot\bm{n}ds(\bm{x})
    &\leq M_3\left(\int_{\partial{D}}|{R_{sb1}^{CH}}|^2 ds(\bm{x})\right)^{1/2}+M_4\left(\int_{\partial{D}}|{\Delta{R}_{sb2}^{CH}|^2 ds(\bm{x})}\right)^{1/2}
    \end{align}
\end{subequations}
where $M_1=|\partial{D}|^{1/2}\cdot(\|{u}\|_{C^1\left(\partial{D}\times \left[0,T\right]\right)}+
\|{u_\theta}\|_{C^1\left(\partial{D}\times \left[0,T\right]\right)}),M_2=|\partial{D}|^{1/2}\cdot\left(\|{\mu}\|_{C\left(\partial{D}\times \left[0,T\right]\right)}+
\|{ \mu_\theta}\|_{C(\partial{D}\times \left[0,T\right])}\right)$,
$M_3=|\partial{D}|^{1/2}\cdot\left(\|{\mu}\|_{C^1\left(\partial{D}\times \left[0,T\right]\right)}+\|{\mu_\theta}\|_{C^1\left(\partial{D}\times \left[0,T\right]\right)}\right),M_4=|\partial{D}|^{1/2}\cdot\left(\|{u}\|_{C\left(\partial{D}\times \left[0,T\right]\right)}+\|{u_\theta}\|_{C\left(\partial{D}\times \left[0,T\right]\right)}\right).$

Let $\hat{M}=\max\{M_1,M_2,M_3,M_4\}$. By combining Eq. (\ref{eq::CH_int_dt}) and (\ref{eq::CH_sb}),we have
\begin{equation}
\label{eq::CHdint}
\begin{aligned}
    &\quad\quad { \frac{1}{2}\frac{d}{dt}}\|\hat{u}\|_{L^2(D)}^2+(\alpha\kappa)\int_D|\Delta\hat{u}|^2 d\bm{x}\\
    & \leq \frac{1}{2}\int_D|R_{int}^{CH}|^2 d\bm{x}+\frac{1+2\kappa M}{2}\int_D|\hat{u}|^2 d\bm{x}+\frac{1}{2}(\alpha\kappa)\int_D \mu^2 d\bm{x}
    +(\alpha\kappa)\hat{M}\left(\int_{\partial{D}}|\Delta{R_{sb1}^{CH}}|^2ds(\bm{x})\right)^{1/2}\\
    &\quad\quad+(\alpha\kappa)\hat{M}\left(\int_{\partial{D}}|{R_{sb2}^{CH}}|^2ds(\bm{x})\right)^{1/2}+(\alpha\kappa)\hat{M}\left(\int_{\partial{D}}|{R_{sb1}^{CH}}|^2 ds(\bm{x})\right)^{1/2}+(\alpha\kappa)\hat{M}\left(\int_{\partial{D}}|{\Delta{R}_{sb2}^{CH}|^2 ds(\bm{x})}\right)^{1/2}
    \end{aligned}    
\end{equation}

By integrating inequality (\ref{eq::CHdint}) over $[0,T']$ for any $T'\leq T$ and applying the Cauchy-Schwarz inequality, we can derive the following result:
{ 
$$
    \begin{aligned}
        &\quad\quad\quad\int_D|\hat{u}(\bm{x},T')|^2d\bm{x}+(2\alpha\kappa)\int_0^{T'}\int_D|\hat{\mu}|^2 d\bm{x}dt\\
        &\leq\int_0^{T}\int_D|R_{int2}^{CH}|^2 d\bm{x} dt+(2\alpha\kappa)\hat{M}T^{1/2}{\hat{C}}+(\alpha \kappa)\int_{0}^T\int_D |R^{CH}_{int1}|^2d\bm{x}dt\\
        &\quad +(1+2\kappa{M})\int_0^{T'}\int_D|\hat{u}|^2d\bm{x}dt+(1+2\kappa M)(2\alpha \kappa)\int_{0}^{T'}\int_D |\hat{\mu}|^2 d\bm{x}dt.
    \end{aligned}
$$
}
(The superfluous term $\int_D|R_{tb}^{CH}|^2d\bm{x}$ has been removed.)
where 
$$
    \begin{aligned}
        \hat{C}= \left(\int_0^T\int_{\partial{D}}|\Delta{R_{sb1}^{CH}}|^2ds(\bm{x})dt\right)^{1/2}+\left(\int_0^T\int_{\partial{D}}|{R_{sb2}^{CH}}|^2ds(\bm{x})dt\right)^{1/2}
        & +\left(\int_0^T\int_{\partial{D}}|{R_{sb1}^{CH}}|^2 ds(\bm{x})dt\right)^{1/2}\\
        & +
\left(\int_0^T\int_{\partial{D}}|{\Delta{R}_{sb2}^{CH}|^2 ds(\bm{x})}dt\right)^{1/2}.
    \end{aligned}
$$

To simplify the expression, let us define $C_G$ to represent the combined error terms:
{ 
$$
C_G = (\alpha \kappa)\int_{0}^T\int_D |R^{CH}_{int1}|^2d\bm{x}dt+\int_0^{T}\int_D|R_{int2}^{CH}|^2 d\bm{x} dt
+2(\alpha\kappa)\hat{M}T^{\frac{1}{2}}\hat{C}.
$$
}
(The term $\int_D|R_{tb}^{CH}|^2d\bm{x}$ is omitted from $C_G$ as superfluous.) Setting
\begin{equation*}
    y(\xi)=\int_D|\hat{u}(\bm{x},\xi)|^2d\bm{x}+(2\alpha\kappa)\int_0^{\xi}\int_D|\hat{\mu}|^2 d\bm{x}d\tau,
\end{equation*}
we have
$$y(T')\leq C_G+(1+2\kappa{M})\int_0^{T'}y(t) dt.$$

{  Applying the integral form of the Grönwall inequality (Lemma~\ref{lemma::Gronwall}), we obtain $y(t)\leq C_G\cdot \exp((1+2\kappa{M})t)$ for each $t\in[0,T]$. Integrating over $t\in[0,T]$ yields
$$ \int_0^T y(t)\,dt \leq C_G\int_0^T \exp((1+2\kappa{M})t)\,dt = C_G\,\frac{\exp((1+2\kappa{M})T)-1}{1+2\kappa{M}}. $$
The total error $\mathcal{E}^{CH}(\theta)^2$ is bounded above by a constant multiple of $\int_0^T y(t)\,dt$, hence
$$ \mathcal{E}^{CH}(\theta)^2 \leq C_G\,\frac{\exp((1+2\kappa{M})T)-1}{1+2\kappa{M}}. $$
This yields the Theorem.}

 
{Theorem \ref{theorem::ch7} establishes a bridge between the total error $\mathcal{E}^{CH}(\theta)^2$ and the generalization error $\mathcal{E}^{CH}_G(\theta)^2$, which is defined by the computable residual terms.

The core implication of this theorem is that it ensures that within our framework, a learned solution cannot be arbitrarily distant from the true solution, provided that the PINN training successfully minimizes the residuals defined by the governing equation, boundary, and initial conditions.

This lays the theoretical foundation for our subsequent result, which connects the total error to the training error.}
{\begin{theorem}
\label{theorem::CH2_TER_Numerical}
Let $d \in \mathbb{N}$ and $T>0$. Let $u \in C^4(D\times[0,T])$ and $\mu\in C^2(D\times[0,T])$. $u_\theta$ and $\mu_\theta$ denote the PINN approximation with parameter $\theta $. Then the total error satisfies
\begin{equation}
\begin{aligned}
\mathcal{E}^{CH}(\theta)^2 &{  \leq C_T\,\frac{\exp((1+2\kappa{M})T)-1}{1+2\kappa{M}}}\\
&=\mathcal{O}\left(\mathcal{E}^{CH}_T\left(\theta,\mathcal{S}^{CH}\right)^2 + {  M_{t b}^{-\frac{1}{2}} + M_{i n t}^{-\frac{1}{2}} + M_{s b}^{-\frac{1}{2}}}\right), 
\end{aligned}
\end{equation}
where $M$ is a constant and $C_T$ is given by
$$
\begin{aligned}
C_T= & {  C_{\left[(R_{t b }^{CH})^2\right]} M_{t b}^{-\frac{1}{2}}}+\mathcal{Q}_{M_{t b}}^D\left[\left(R_{t b }^{CH}\right)^2\right]+  {  C_{\left[\left(R_{i n t2 }^{CH}\right)^2\right]} M_{i n t}^{-\frac{1}{2}}}+\mathcal{Q}_{M_{i n t}}^{\Omega}\left[\left(R_{i n t 2}^{CH}\right)^2\right] \\
& +\left(\alpha \kappa\right)\left( {  C_{\left[\left(R_{i n t1 }^{CH}\right)^2\right]} M_{i n t}^{-\frac{1}{2}}}+\mathcal{Q}_{M_{i n t}}^{\Omega}\left[\left(R_{i n t 1}^{CH}\right)^2\right] \right)\\
&+2\left(\alpha\kappa\right)\hat{M}\left({  C_{\left[\left(R_{s b 1}^{CH}\right)^2\right]} M_{s b}^{-\frac{1}{2}}}+\mathcal{Q}_{M_{s b}}^{\partial{D}\times \left[0,T\right]}\left[\left(R_{s b 1}^{CH}\right)^2\right]\right)^{1/2}\\
&+2\left(\alpha\kappa\right)\hat{M}\left({  C_{\left[\left(R_{s b 2}^{CH}\right)^2\right]} M_{s b}^{-\frac{1}{2}}}+\mathcal{Q}_{M_{s b}}^{\partial{D}\times \left[0,T\right]}\left[\left(R_{s b 2}^{CH}\right)^2\right]\right)^{1/2}\\
& + 2\left(\alpha\kappa\right)\hat{M}\left({  C_{\left[\left(\Delta R_{s b 1}^{CH}\right)^2\right]} M_{s b}^{-\frac{1}{2}}}+\mathcal{Q}_{M_{s b}}^{\partial{D}\times \left[0,T\right]}\left[\left(\Delta R_{s b 1}^{CH}\right)^2\right]\right)^{1/2}\\
&+2\left(\alpha\kappa\right)\hat{M}\left({  C_{\left[\left(\Delta R_{s b 2}^{CH}\right)^2\right]} M_{s b}^{-\frac{1}{2}}}+\mathcal{Q}_{M_{s b}}^{\partial{D}\times \left[0,T\right]}\left[\left(\Delta R_{s b 2}^{CH}\right)^2\right]\right)^{1/2},
\end{aligned}
$$
with the component $\hat{M}$ defined as
\begin{equation*}
    \hat{M} = \max_{g \in \{u,\mu \} }\{|\partial{D}|^{1/2}\cdot\left(\|g\|_{C^1(\partial{D}\times[0,T])}+\|g_\theta \|_{C^1(\partial{D}\times [0,T])}\right)\}.
\end{equation*}
\end{theorem}}
\textbf{Proof of Theorem \ref{theorem::CH2_TER_Numerical}} {  By combining Theorem \ref{theorem::ch7} with the Monte-Carlo (LHS) sampling error formula (\ref{equ::quadrules}), we can establish a bound for the generalization error constant $C_G$ in terms of the training error.} Recall from Theorem \ref{theorem::ch7} that $C_G$ is defined by the true integrals of the various residual terms. {  We now apply the Monte-Carlo (LHS) sampling error formula to each integral term that constitutes $C_G$, yielding the following bounds (with $O(M^{-1/2})$ sampling error):}
\begin{equation}\label{eq:quadrature_bounds_block_ch}
\begin{gathered}
\begin{aligned}
\int_\Omega\left|R_{t b }^{CH}\right|^2 \mathrm{d} \boldsymbol{x} & \leq {  C_{\left[\left(R_{t b }^{CH}\right)^2\right]} M_{t b}^{-\frac{1}{2}}}+\mathcal{Q}_{M_{t b}}^D\left[\left(R_{t b }^{CH}\right)^2\right]
\end{aligned}
\\
\begin{aligned}
\int_0^T\int_{D}\left|R_{int1 }^{CH}\right|^2 \mathrm{d} \boldsymbol{x} \mathrm{d} t & \leq {  C_{\left[\left(R_{int1}^{CH}\right)^2\right]} M_{int}^{-\frac{1}{2}}}+\mathcal{Q}_{M_{int}}^{\Omega}\left[\left(R_{int1}^{CH}\right)^2\right]
\end{aligned}
\\
\begin{aligned}
\int_0^T\int_{D}\left|R_{int2 }^{CH}\right|^2 \mathrm{d} \boldsymbol{x} \mathrm{d} t & \leq {  C_{\left[\left(R_{int2}^{CH}\right)^2\right]} M_{int}^{-\frac{1}{2}}}+\mathcal{Q}_{M_{int}}^{\Omega}\left[\left(R_{int2}^{CH}\right)^2\right]
\end{aligned}
\\
\begin{aligned}
    \int_0^T\int_{\partial D}\left|R_{sb1 }^{CH}\right|^2 ds(\boldsymbol{x}) \mathrm{d} t & \leq {  C_{\left[\left(R_{sb1}^{CH}\right)^2\right]} M_{sb}^{-\frac{1}{2}}}+\mathcal{Q}_{M_{sb}}^{\partial D\times[0,T]}\left[\left(R_{sb1}^{CH}\right)^2\right]
\end{aligned}
\\
\begin{aligned}
    \int_0^T\int_{\partial D}\left|\Delta R_{sb1 }^{CH}\right|^2 ds(\boldsymbol{x}) \mathrm{d} t & \leq {  C_{\left[\left(\Delta R_{sb1}^{CH}\right)^2\right]} M_{sb}^{-\frac{1}{2}}}+\mathcal{Q}_{M_{sb}}^{\partial D\times[0,T]}\left[\left(\Delta R_{sb1}^{CH}\right)^2\right]
\end{aligned}
\\
\begin{aligned}
    \int_0^T\int_{\partial D}\left| R_{sb2 }^{CH}\right|^2 ds(\boldsymbol{x}) \mathrm{d} t
    & \leq {  C_{\left[\left(R_{sb2}^{CH}\right)^2\right]} M_{sb}^{-\frac{1}{2}}}+\mathcal{Q}_{M_{sb}}^{\partial D\times[0,T]}\left[\left( R_{sb2 }^{CH}\right)^2\right]
\end{aligned}
\\
\begin{aligned}
    \int_0^T\int_{\partial D}\left|\Delta R_{sb2 }^{CH}\right|^2 ds(\boldsymbol{x}) \mathrm{d} t
    & \leq {  C_{\left[\left(\Delta R_{sb2}^{CH}\right)^2\right]} M_{sb}^{-\frac{1}{2}}}+\mathcal{Q}_{M_{sb}}^{\partial D\times[0,T]}\left[\left(\Delta R_{sb2 }^{CH}\right)^2\right]
\end{aligned}
\end{gathered}
\end{equation}

By substituting these individual bounds back into the definition of $C_G$, we find that $C_G$ is bounded by a new quantity, which we define as $C_T$. This relationship, $C_G \le C_T$, allows us to state the final error bound in terms of trainable quantities. By Theorem \ref{theorem::ch7}, it therefore holds that
{ 
$$
\mathcal{E}(\theta)^2=||\hat{u}(\bm{x},t)||^2_{L^2(D \times [0,T])}+2\alpha\kappa\int_0^T||\hat{\mu}(\bm{x},t)||^2_{L^2(D)} dt \leq C_T\,\frac{\exp((1+2\kappa M)T)-1}{1+2\kappa M}.
$$
}
where $C_T$ is given by
$$
\begin{aligned}
C_T :=~ & {  C_{\left[{\left(R_{t b }^{CH}\right)}^2\right]} M_{t b}^{-\frac{1}{2}}}+\mathcal{Q}_{M_{t b}}^D\left[{\left(R_{t b }^{CH}\right)}^2\right] + (\alpha\kappa) \left( {  C_{\left[\left(R_{int1 }^{CH}\right)^2\right]} M_{i n t}^{-\frac{1}{2}}}+\mathcal{Q}_{M_{i n t}}^{\Omega}\left[\left(R_{i n t1 }^{CH}\right)^2\right] \right) \\
& +{  C_{\left[\left(R_{int2}^{CH}\right)^2\right]} M_{i n t}^{-\frac{1}{2}}}+\mathcal{Q}_{M_{i n t}}^{\Omega}\left[\left(R_{i n t2 }^{CH}\right)^2\right] \\
& +2(\alpha\kappa)\hat{M}\left({  C_{\left[\left(R_{s b 1}^{CH}\right)^2\right]} M_{s b}^{-\frac{1}{2}}}+\mathcal{Q}_{M_{s b}}^{\partial D\times[0,T]}\left[\left(R_{s b 1}^{CH}\right)^2\right]\right)^{1/2}\\
&+2(\alpha\kappa)\hat{M}\left({  C_{\left[\left(R_{s b 2}^{CH}\right)^2\right]} M_{s b}^{-\frac{1}{2}}}+\mathcal{Q}_{M_{s b}}^{\partial D\times[0,T]}\left[\left(R_{s b 2}^{CH}\right)^2\right]\right)^{1/2} \\
& + 2(\alpha\kappa)\hat{M}\left({  C_{\left[\left(\Delta R_{s b 1}^{CH}\right)^2\right]} M_{s b}^{-\frac{1}{2}}}+\mathcal{Q}_{M_{s b}}^{\partial D\times[0,T]}\left[\left(\Delta R_{s b 1}^{CH}\right)^2\right]\right)^{1/2}\\
&+2(\alpha\kappa)\hat{M}\left({  C_{\left[\left(\Delta R_{s b 2}^{CH}\right)^2\right]} M_{s b}^{-\frac{1}{2}}}+\mathcal{Q}_{M_{s b}}^{\partial D\times[0,T]}\left[\left(\Delta R_{s b 2}^{CH}\right)^2\right]\right)^{1/2},
\end{aligned}
$$
and $\hat{M}=\max\{M_u,M_\mu\}$, $M_u=|\partial{D}|^{1/2}\cdot(\|u\|_{C^1(\partial{D}\times \left[0,T\right])}+\|u_\theta \|_{C^1(\partial{D}\times \left[0,T\right])})$

Note that the complexities of $C_{\left[( R_{p}^{CH})^2\right]} , p=tb,int1,int2,sb1,sb2$ and $C_{\left[(\Delta R_{q}^{CH})^2\right]} , q=sb1,sb2$ can be derived in a manner similar to that in Theorem \ref{theorem::ac5}.

\subsection{\textbf{Numerical Examples}}
\subsubsection{1D Cahn--Hilliard equation}
We will first examine the performance of the 1D Cahn--Hilliard equation, which describes the process of phase separation and interface development. The spatial-temporal domain is taken as $(x,t) \in D \times [0,T] = [-1,1]\times[0,0.5]$. Then the equation is given by:
\begin{subequations}
\begin{align}
\label{eq::ch1dnum}
& u_t-\left(\kappa\left(u^3-u\right)-\alpha\kappa u_{x x}\right)_{x x}=0,   \qquad (x,t) \in[-1,1]\times[0,0.5] \\
& u(x,0)=-\cos \left(2 \pi x
\right), \qquad (x,t)\in[-1,1]\times\{0\}\\
& u(-x,t)=u(x,t), \qquad (x,t)\in \partial D\times[0,0.5]\\
& \nabla u(-x,t)=\nabla u(x,t), \qquad (x,t)\in \partial D \times[0,0.5] 
\end{align}
\end{subequations}
where the parameter $\kappa = 1$ serves as the mobility parameter, influencing the rate of diffusion of components, while the parameter $\alpha = 0.02$ is associated with the surface tension at the interface. To implement the Residuals-RAE-PINN scheme for the 1D Cahn--Hilliard equation, we first have the following modified loss function 
\begin{align}
\nonumber
\mathcal{L}(\theta) &=  \gamma_{int} \cdot \mathcal{L}_{int}(\theta) + \gamma_{tb} \cdot \mathcal{L}_{tb}(\theta)+ \gamma_{sb} \cdot \mathcal{L}_{sb}(\theta) \\
\nonumber
&= \gamma_{int} \cdot \frac{1}{M_{int}} \sum_{i=1}^{M_{int}} \hat{\lambda}_{int}^{(i)} \left| u_t(x_{int}^{(i)}, t_{int}^{(i)})-\left(\kappa\left(u(x_{int}^{(i)}, t_{int}^{(i)})^3-u(x_{int}^{(i)}, t_{int}^{(i)})\right)-\alpha\kappa u(x_{int}^{(i)}, t_{int}^{(i)})_{x x}\right)_{x x}\right| ^2\\
\nonumber
&\quad+  \gamma_{tb} \cdot \frac{1}{M_{tb}} \sum_{i=1}^{M_{tb}} \left|u(x_{tb}^{(i)}, 0) - u_{0}(x_{tb}^{(i)}) \right| ^2\\
\nonumber
&\quad+  \gamma_{sb}\left( \cdot \frac{1}{M_{sb}} \sum_{i=1}^{M_{sb}}  (\left|u(-x_{sb}^{(i)}, t_{sb}^{(i)}) - u(x_{sb}^{(i)}, t_{sb}^{(i)})\right| ^2 \right)^{1/2}\\
\nonumber
&\quad+\gamma_{sb} \cdot\left( \frac{1}{M_{sb}} \sum_{i=1}^{M_{sb}} \left|
\nabla u(-x_{sb}^{(i)}, t_{sb}^{(i)}) - \nabla u(x_{sb}^{(i)}, t_{sb}^{(i)})\right|^2\right)^{1/2}.
\end{align}

For solving the 1D Cahn--Hilliard equation, we assign the weights for different loss terms at $\gamma_{int} = 1, \gamma_{sb} = 1, \gamma_{tb} = 100$. Furthermore, we incorporate a Residuals-RAE weighting scheme for the pointwise weights of residual points, denoted by $\hat{\lambda}_{int}^{(i)}$.  The impact of the nearest points on these weights is controlled by the hyper-parameter $k_{int}$, which is set to a value of 50. This parameter plays a pivotal role in determining the extent to which neighboring points influence the calculation of pointwise weights. In this case, we use the training points set $\{(x_{int}^{(i)}, t_{int}^{(i)})\}$, $\{(-1, t_{sb}^{(i)})\}$, $\{(x_{tb}^{(i)})\}$ obtained from the Latin hypercube sampling approach. Specifically, we have $M_{int} = 10,000$ training points for the residual, $M_{sb} = 256$ points for the boundary conditions, and $M_{tb} = 512$ points for the initial conditions.

\begin{figure}
\centering
    \subcaptionbox{Results for 1D Cahn--Hilliard equation using Residuals-RAE PINNs (with 200 neurons per layer).\label{fig::1dch_error} }{\includegraphics[width=1.0\linewidth, height=0.45\linewidth]{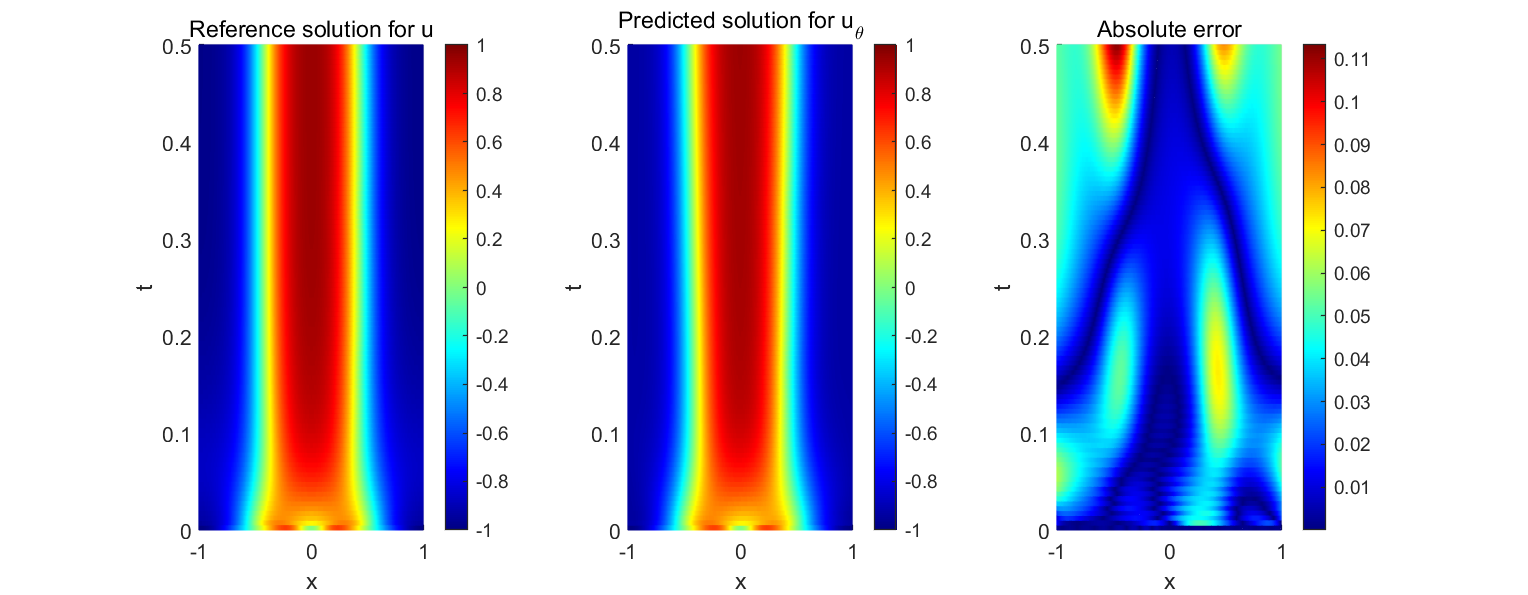}}

    \subcaptionbox{Reference and Residuals-RAE-PINN predicted solutions (neurons per layer $\#1$ and  $\#2$) of the 1D Cahn--Hilliard equation (initial condition $\#$ 2) at different time snapshots (a) $t=0$, (b) $t=0.25$, (c) $t=0.5$.\label{fig::1dCH_difft}}{\includegraphics[width=1.0\linewidth]{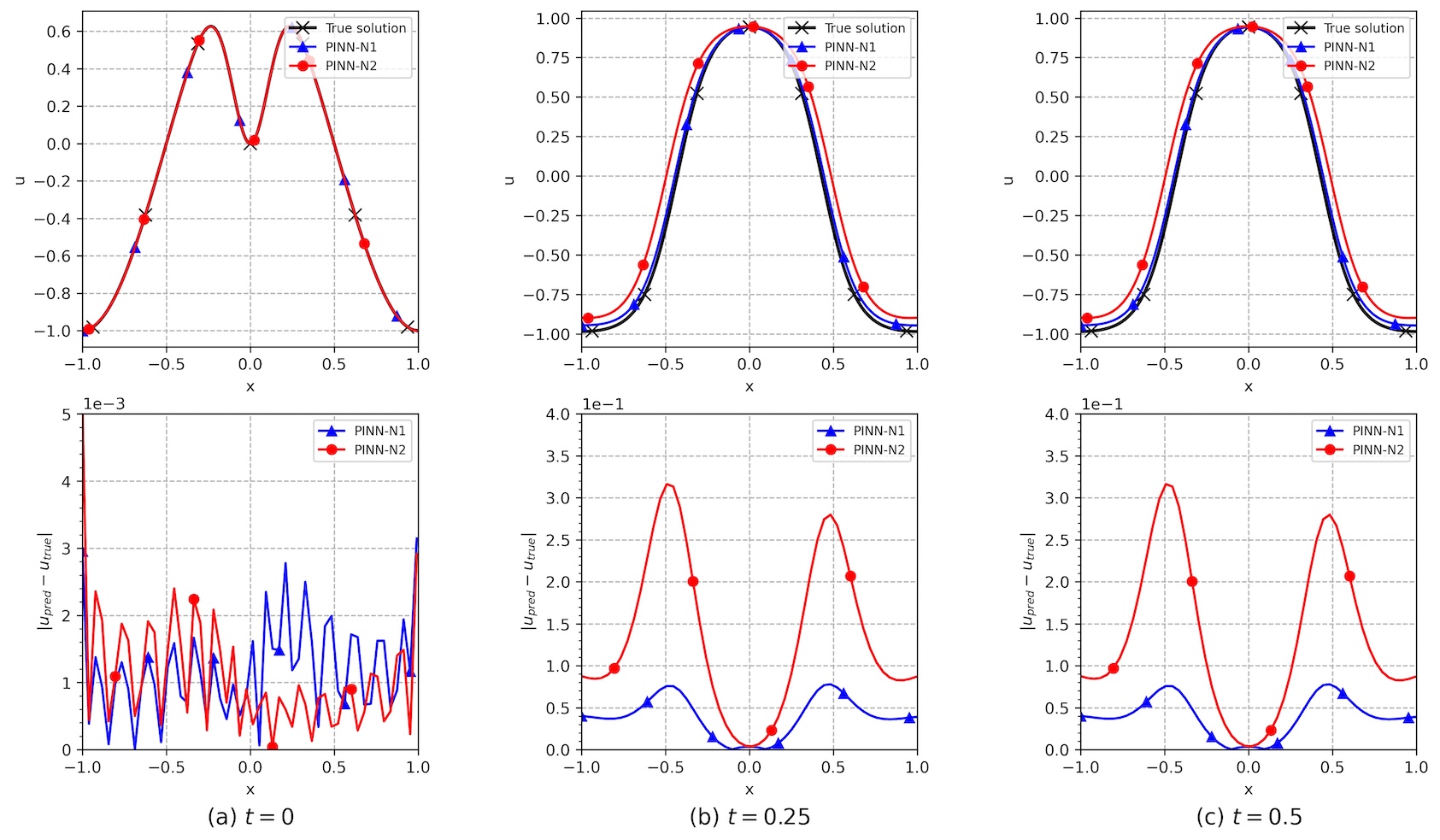}}
    \caption{\textbf{\textbf{Results for solving 1D Cahn--Hilliard using Residuals-RAE-PINN.}} \label{fig::diff_methods_ch}}
\end{figure}

\begin{table}[h]\caption{{\bf {\color{black}{{Description of training data for 1D Cahn--Hilliard equation.}}} \label{tab::data_ch1d}}}
$$
\begin{array}{lll}
\hline M_{tb} & \text { Initial collocation points } & 512 \\
M_{sb} & \text { Boundary collocation points } & 256 \\
M_{int} & \text { Residual collocation points } & 10,000 \\
\hline
\end{array}
$$
\end{table}

The training process utilizes a neural network architecture with two layers, each consisting of 256 nodes. The network is initialized using the Xavier scheme. The optimization of the loss function involves two stages: to begin with, the network undergoes 5,000 iterations using the Adam optimizer, followed by a refinement stage of 1,000 iterations with the L-BFGS algorithm. The training is deemed complete when the reduction in loss between successive epochs is below a certain threshold, specifically less than $10^{-7}$. The effectiveness of this Residuals-RAE method is quantified by the relative $L^2$ error, with which yielding errors in the order of $5\times 10^{-2}$. The pointwise error for the results obtained from the PINN at different time snapshots is plotted in Fig. \ref{fig::1dCH_difft}. It can be observed that the results of PINN-NN1 with 512 neurons and PINN-NN2 256 neurons per layer produce the displayed results. 

\begin{figure}
\centering
\includegraphics[width=0.7\linewidth, height = 0.40\linewidth]{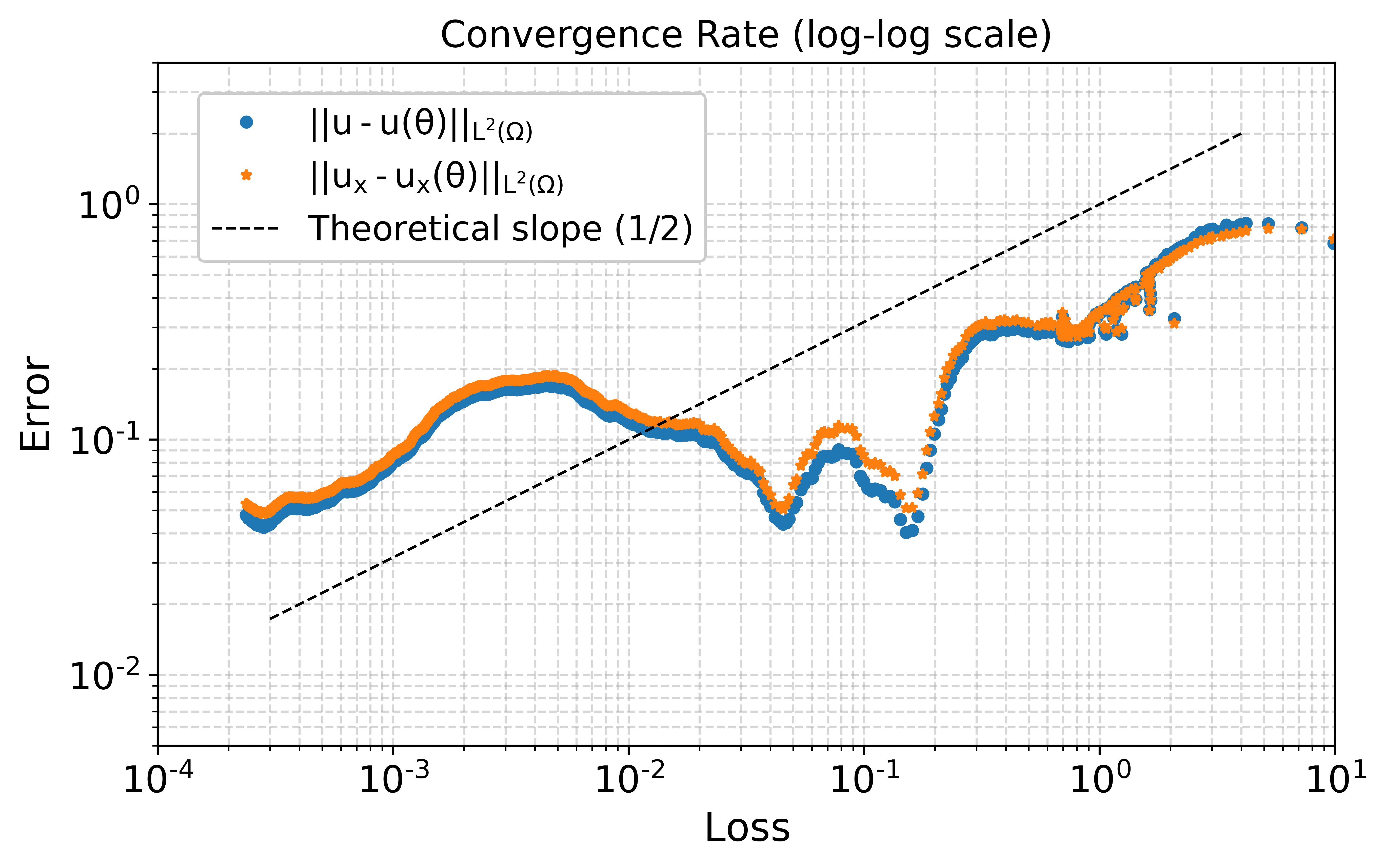}
\caption{\textbf{1D Cahn--Hilliard equation: The dynamics of $l^2$ errors of $u$ and $u_x$ as a function of the training loss value using Residuals-RAE.}\label{fig::1dCH_loss}}
\end{figure}

Fig. \ref{fig::1dCH_loss} illustrates the solution error versus training loss relationship during the solution process of the Cahn--Hilliard equation using a Physics-Informed Neural Network (PINN). The numerical results reveal that overall trend of the total error is approximately proportional to the square root of the training loss. This trend is consistent with the expected theoretical behavior as stated in Theorem \ref{theorem::CH2_TER_Numerical}. The fluctuations observed in the middle range of loss values appear to scale with a power slightly greater than 1/2. These fluctuations are believed to be a result of the involvement of higher-order derivatives in the Cahn--Hilliard equation, which can introduce complexities in the training dynamics. 
{While the training dynamics for the Cahn-Hilliard equation exhibit significant fluctuations, likely due to the complexity of the fourth-order derivatives, the overall trajectory demonstrates that the solution error decreases as the training loss is minimized. This successful convergence, despite the challenging optimization landscape, validates the effectiveness of our approach for this intricate PDE.}

{ Adhering to the unified experimental protocol established in Section \ref{sec:1D Allen--Cahn equation}, the comparative results for the Cahn-Hilliard equation (Table \ref{tab::comp_ch1d}) were obtained using the standard architecture (e.g., 256 neurons), consistent with the unified settings detailed in Table \ref{tab:hyperparameters_detailed} (in Appendix \ref{sec:Detailed Experimental Settings}). This ensures that the performance comparison against baselines remains strictly fair and unaffected by variations in model capacity.}


\begin{table}[htbp]
\centering
\caption{{ Relative $l_2$-error and $l_\infty$-error comparison between different methods for 1D Cahn--Hilliard Equations using residuals-RAE-PINNs} \label{tab::comp_ch1d}}
\begin{tabular}{ccccc}
\toprule
\textnormal{PINN method}
& \multicolumn{1}{@{\hskip 0pt}l}{\hspace{0.8em}Relative $l_2$-error } 
& \multicolumn{2}{@{\hskip 0pt}l}{\hspace{0.8em}Relative $l_\infty$-error} \\
\midrule
baseline \cite{raissi2019physics} & 3.725e-01 & \hspace{1.2em} 1.183e-01 \\
SA-PINNs \cite{mcclenny2023self}  & 2.896e-01 & \hspace{1.2em} 9.252e-02 \\
PINNs-WE \cite{liu2024discontinuity}   & 4.724e-01 & \hspace{1.2em} 1.398e-01 \\
{\bf Residuals-RAE-PINNs}   & 4.829e-02 & \hspace{1.2em} 2.002e-02  \\
\bottomrule
\end{tabular}
\end{table}


\subsubsection{2D Cahn--Hilliard equation}
We proceed to evaluate the effectiveness of the Residuals-RAE-PINN algorithm in solving the 2D Cahn--Hilliard equation, as described in \cite{mattey2022novel}. In this study, the authors employ a re-training approach using the same neural network as in previous time segments, similar to the concept of sequence-to-sequence learning \cite{krishnapriyan2021characterizing}.  The spatial-temporal domain is taken as $(x,y,t) \in D \times [0,T] = [0,1]^2\times[0,0.003]$, and the initial-boundary value problem
\begin{subequations}
   \begin{align}
& \partial_t u=\kappa \Delta\left(-\alpha^2 \Delta u+u^3-u\right), \qquad  (x,y,t) \in D \times [0,T]\\
& u(x, y,0)= 0.4 \cos{(3\pi x)}\cos{(3 \pi y)}
, \qquad (x,y,t)\in D\times[0,T]\\
&u^{(d-1)}(x, y, t)=u^{(d-1)}(-x,-y, t),\qquad  (x,y,t)\in \partial D\times[0,T]\quad(d=1,2)
\end{align} 
\end{subequations}
where $\kappa = 1, \alpha = 0.02$, { $u^{(0)}$
was intended to represent the zeroth derivative, which is the function value $u$ and $u^{(1)}$ was intended to represent the first derivative, specifically the gradient $\nabla u$.}
To solve the problem by Residuals-RAE-PINN, we first consider the modified loss function as below
\begin{equation}
\begin{aligned}
\mathcal{L}(\theta) &=  \gamma_{int} \cdot \mathcal{L}_{int}(\theta)  + \gamma_{tb} \cdot \mathcal{L}_{tb}(\theta)
+ \gamma_{sb} \cdot \mathcal{L}_{sb}(\theta) \\
&= \gamma_{int} \cdot \frac{1}{M_{int}} \sum_{i=1}^{M_{int}} \hat{\lambda}_{int}^{(i)} \left| u_t(x_{int}^{(i)},y_{int}^{(i)}, t_{int}^{(i)})+\alpha\kappa \Delta^2 u(x_{int}^{(i)}, y_{int}^{(i)},t_{int}^{(i)})-\kappa\Delta (u^3(x_{int}^{(i)},y_{int}^{(i)}, t_{int}^{(i)})-u(x_{int}^{(i)},y_{int}^{(i)}, t_{int}^{(i)}))\right| ^2\\
&\quad+  \gamma_{tb} \cdot \frac{1}{M_{tb}} \sum_{i=1}^{M_{tb}} \left|u(x_{tb}^{(i)},y_{tb}^{(i)}, 0) -0.4\cos{(3\pi x_{tb}^{(i)})\cos{(3\pi y_{tb}^{(i)}})} \right| ^2\\
&\quad+  \gamma_{sb} \cdot \frac{1}{M_{sb}} \sum_{i=1}^{M_{sb}} \left( \sum_{d=1}^2\left|u^{(d-1)} (x_{sb}^{(i)}, y_{sb}^{(i)}, t_{sb}^{(i)})-u^{(d-1)}(-x_{sb}^{(i)}, -y_{sb}^{(i)}, t_{sb}^{(i)})\right|^2\right)^{1/2},
\end{aligned}    
\end{equation}
where the weights for the loss terms  $\gamma_{int}, \gamma_{sb}, \gamma_{tb}$, are set to $1, 100, 100$ respecitively. The input layer of the PINN structure {\color{blue}contains 3 neurons}, which is the same configuration used for solving the 2D Allen--Cahn equation. The activation function chosen for this neural network is the hyperbolic tangent function ($tanh$).

To generate the collocation points, we randomly sample $M_{int} = 10,000$, $M_{sb} = 256$, $M_{tb} = 512$ points from the computational domain using the Latin hypercube sampling technique. In order to obtain pointwise weights, we apply the Residuals-RAE weighting scheme with a hyper-parameter $k_{int}$ set to $50$, controlling the influence of nearby points on the information obtained.

\begin{table}[h]\caption{{\bf {\color{black}{{Description of training data for 2D Cahn--Hilliard equation.}}} \label{tab::data_ch2d}}}
$$
\begin{array}{lll}
\hline M_{tb} & \text { Initial collocation points } & 512 \\
M_{sb} & \text { Boundary collocation points } & 256 \\
M_{int} & \text { Residual collocation points } & 10,000 \\
\hline
\end{array}
$$
\end{table}

\begin{figure}
\centering
\includegraphics[width=1.0\linewidth]{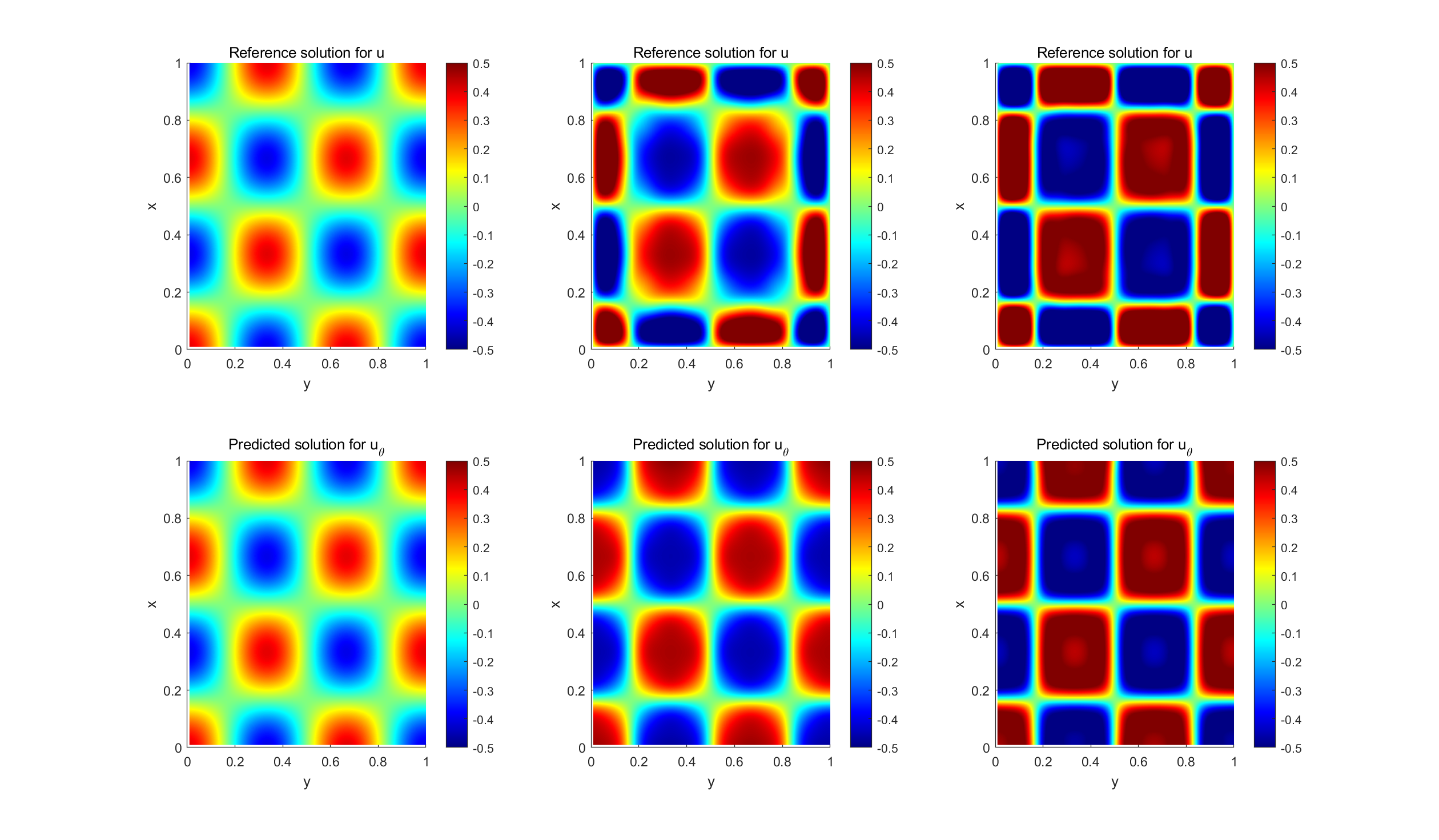}
\caption{\textbf{Results for 2D Cahn--Hilliard equation using Residuals-RAE PINNs.}\label{fig::2dch_difft}}
\end{figure}

During the training process, the Xavier-normal method is employed to initialize trainable parameters in the network. To minimize the modified loss function, $50,000$ ADAM and $10,000$ LBFGS iterations are applied. Additionally, the learning rate exponential decay method is employed, with an exponential decay rate of $0.95$ every 5000 epochs in the numerical experiments. It's important to note that due to the strong non-linearity and higher-order nature of the Cahn--Hilliard equation, accurately capturing the solution with a shallow neural network (e.g. depth = 2) is challenging. As a result, several approaches have been proposed in the field of Physics-Informed Neural Networks (PINNs), such as self-adaptive samples, domain decompositions, and self-adaptive weights. Furthermore, in order to ensure the representability of PINNs, all of these approaches utilize deep or deeper neural networks with multiple layers ($L \geq 5$). From the observation of Fig. \ref{fig::2dch_difft}, it can be seen that the Residuals-RAE-PINN is capable of capturing the solutions for the Cahn--Hilliard equation without requiring domain decomposition in the time dimension, whereas the method proposed in the bc-PINN \cite{mattey2022novel} necessitates it.

\begin{figure}
\centering
\includegraphics[width=0.7\linewidth, height = 0.40\linewidth]{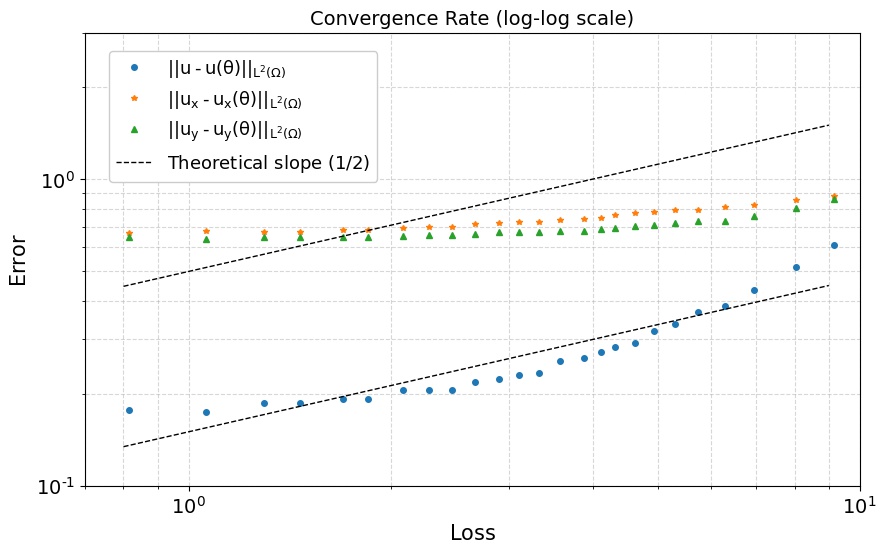}
\caption{\textbf{2D Cahn--Hilliard equation: The dynamics of $l^2$ errors of $u$, $u_x$ and $u_y$ as a function of the training loss value using Residuals-RAE.}\label{fig::2dCH_loss}}
\end{figure}

To further investigate the relationship between solution errors and the training loss dynamics, we present a log-log plot of $u$, $u_x$, and $u_y$ in Fig. \ref{fig::2dCH_loss}. {This type of plot is commonly used in previous studies to analyze the training dynamics of PINNs. In our figure, the reference line with a slope of $1/2$ serves as a theoretical benchmark, indicating the expected square-root scaling of the solution error with respect to the training loss. By examining the convergence trend relative to this reference, we can observe that the training dynamics broadly follow the anticipated behavior predicted by the theory.} 

{ Comparisons among different methods for this case are summarized in Table \ref{tab::comp_ch2d}. It can be observed that, despite the use of self-adaptive weighting schemes, the solutions obtained by various methods remain insufficiently accurate. This is likely due to the inherent complexity of the 2D Cahn–Hilliard equation, particularly the presence of high-order derivatives, which poses significant challenges for shallow neural networks with only two layers to capture effectively. To ensure a fair comparison, all shared hyperparameters—including the neural network architecture, number of training iterations, and number of collocation points—are kept consistent across all experiments. Furthermore, adhering to the unified evaluation protocol established in Section \ref{sec:1D Allen--Cahn equation}, we applied the consistent method-specific hyperparameters listed in Table \ref{tab:hyperparameters_detailed} (in Appendix \ref{sec:Detailed Experimental Settings}) for all competing methods, guaranteeing that the comparative results in Table \ref{tab::comp_ch2d} reflect the optimized performance of each baseline.}


\begin{table}[htbp]
\centering
\caption{{ Relative $l_2$-error and $l_\infty$-error comparison between different methods for 2D Cahn--Hilliard Equations using residuals-RAE-PINNs}\label{tab::comp_ch2d}}
\begin{tabular}{ccccc}
\toprule
\textnormal{PINN method}
& \multicolumn{1}{@{\hskip 0pt}l}{\hspace{0.8em}Relative $l_2$-error } 
& \multicolumn{2}{@{\hskip 0pt}l}{\hspace{0.8em}Relative $l_\infty$-error} \\
\midrule
baseline \cite{raissi2019physics} & 9.213e-01 & \hspace{1.2em} 7.629e-01 \\
SA-PINNs \cite{mcclenny2023self}  & 5.370e-01 & \hspace{1.2em}  2.060e-02 \\
PINNs-WE \cite{liu2024discontinuity}   & 3.049e-01 & \hspace{1.2em} 1.395e-02 \\
{\bf Residuals-RAE-PINNs}   & 2.951e-01 & \hspace{1.2em} 1.390e-02 \\
\bottomrule
\end{tabular}
\end{table}

\section{Conclusion}\label{section::5}
This paper investigates error estimation when using physics-informed neural networks (PINN) to solve the Allen--Cahn (AC) and Cahn--Hilliard (CH) equations, which are commonly used in phase field modeling research. The primary discoveries are as follows:
\begin{itemize}
    \item Error bounds are derived for the total error and approximation error between the PINN solution and exact solution for both the AC and CH equations, assuming certain smoothness conditions, {  with integrals discretized by Monte-Carlo (LHS) sampling ($O(M^{-1/2})$ sampling error).} These results are presented in Theorem \ref{theorem::ac4}, Theorem \ref{theorem::ac5}, Theorem \ref{theorem::ch7} and Theorem \ref{theorem::CH2_TER_Numerical}.
    \item For the AC equation, it is shown that the PINN solution converges at a rate of $O\left(\ln^2\left(N\right)N^{-k+2}\right)$ as { $N$ increases, which demands a corresponding increase in the network's width and architectural complexity. As established in Theorem \ref{theorem::ac4} and Theorem \ref{theorem::ac5}, the parameter $k$ represents the continuity of the exact solution.} 
    {\item For the CH equation, the PINN solution converges as the network complexity parameter $N$ increases, as demonstrated in Theorem \ref{theorem::ch7} and Theorem \ref{theorem::CH2_TER_Numerical}, at a rate of $O\left(\ln^2\left(N\right)N^{-k+4}\right)$.}
    \item In our numerical experiment, we employ a pointwise weighting scheme called Residuals-RAE to address the challenge of solving highly nonlinear PDEs in PINNs.

    \item Numerical results validate the error bounds and convergence rates derived in the theorems.
\end{itemize}

Overall, the paper offers a comprehensive analysis of the errors involved in utilizing PINNs to approximate solutions for the nonlinear AC and CH equations in phase field modeling. The main contributions lie in the derivation of error bounds and convergence rates for the application of the PINN methodology to these specific problems.

\section*{Acknowledgement}
This work received support from SandGold AI Research under the Auto-PINN Research program. Ieng Tak Leong is supported by the Science and Technology Development
Fund, Macau SAR (File no. 0132/2020/A3). Reviewers and editors' comments on improving this study and manuscript will be highly appreciated. High-Performance Computing resources were provided by SandGold AI Research.

\newpage
\begin{appendices}

In this section, we list out several fundamental yet important results in the area of PDE and functional analysis. 
First and foremost, we aim to provide a comprehensive overview of the notations employed in this paper.

\section{Notation}
\label{Notation}

\subsection{Basic notation}
For $n \in \mathbb{N}$, we denote $\{1,2, \ldots n\}$ by $[n]$ for simplicity in the paper.
For two Banach spaces $X, Y, \mathscr{L}(X, Y)$ refers to the set of continuous linear operator mapping from $X$ to $Y$.

For a mapping $f: X \rightarrow Y$, and $u, v \in X, d f(u, v)$ denotes Gateaux differential of $f$ at $u$ in the direction of $v$, and $f^{\prime}(u)$ denotes Fréchet derivative of $f$ at $u$.

For $r>0$, and $x_0 \in X$, where $X$ is a Banach space equipped with norm $\|\cdot\|_X, B_{int}\left(x_0\right)$ refers to $\left\{x \in X:\left\|x-x_0\right\|_X<r\right\}$.

For a function $f: X \rightarrow \mathbb{R}$, where $X$ is a measurable space. We denote by supp $f$ the support set of $f$, i.e. the closure of $\{x \in X: f(x) \neq 0\}$.

\subsection{Multi-index notation}
For $n \in \mathbb{N}$, we call an $n$-tuple of non-negative integers $\alpha \in \mathbb{N}^n$ a multi-index. We use the notation $|\alpha|=\sum_{i=1}^n \alpha_i, \alpha !=\Pi_{i=1}^n \alpha_{i} !$. For $x=\left(x_1, x_2, \ldots x_n\right) \in \mathbb{R}^n$, we denote by $x^\alpha=\Pi_{i=1}^n x_i^{\alpha_i}$ the corresponding multinomial. Given two multi-indices $\alpha, \beta \in \mathbb{N}^d$, we say $\alpha \leq \beta$ if and only if $\alpha_i \leq \beta_i, \forall i \in[n]$.

For an open set $\Omega \subset \mathbb{R}^n, T \in \mathbb{R}^{+}$and a function $f(x): \Omega \rightarrow \mathbb{R}$ or $f(x, t): \Omega \times[0, T] \rightarrow \mathbb{R}$, we denote by
\begin{align}
    D^\alpha f=\frac{\partial^{|\alpha|} f}{\partial x_1^{\alpha_1} \ldots \partial x_n^{\alpha_n}}
\end{align}
the classical or weak derivative of $f$.

For $k \in \mathbb{R}^n$, we denote by $D^k f$ (or $\nabla^k f$ ) the vector whose components are $D^\alpha f$ for all $|\alpha|=k$, and we abbreviate $D^1 f$ as $D f$.

\subsection{Space notation}
 Let $m\in \mathbb{N}$, $1\leq p \leq \infty$ and let $\Omega \subset \mathbb{R}^n $ be open. We denote by $L^p(\Omega) $ the usual Lebesgue space. 

The Sobolev space $W ^{m,p}(\Omega)$ is defined as
\begin{align}
    \left\{f(x) \in L^p(\Omega): D^\alpha f \in L^p(\Omega), \forall \alpha \in \mathbb{N}^n \text { with }|\alpha| \leq m\right\}
\end{align}
The norm we define as
\begin{align}
    \|f\|_{W^{m, p}(\Omega)}:=\left(\sum_{|\alpha| \leq m}\left\|D^\alpha f\right\|_{L^p(\Omega)}^p\right)^{\frac{1}{p}}
\end{align}
for $1\leq p < \infty$, and 
\begin{align}
    \|f\|_{W^{m, \infty}(\Omega)}:=\max _{|\alpha| \leq m}\left\|D^\alpha f\right\|_{L^{\infty}(\Omega)}
\end{align}
for $p=\infty$.

If $p=2$, we  write
\begin{align}
    H^k(U)=W^{k, 2}(U) \quad(k=0,1, \ldots)
\end{align}

The space $C([0, T] ; X)$
comprises all continuous functions $\mathbf{u}:[0, T] \rightarrow X$ with
\begin{align}
 \|\mathbf{u}\|_{C([0, T] ; X)}:=\max _{0 \leq t \leq T}\|\mathbf{u}(t)\|<\infty .
\end{align}

\section{Auxiliary results}
\label{Auxiliary results}
In this section, we will present several fundamental findings and outline the key theorems that have been derived from previous research.
\begin{lemma}[Sobolev embedding Theorem]
\label{lemma::embedding}
Let $\Omega$ be an open set in $\mathbb{R}^n , p \in[1, \infty]$, and $m \leq k$ be a non-negative integer.

(i) If $\frac{1}{p}-\frac{k}{n}>0$, and set $q=\frac{n p}{n-p k}$, then $W^{m, p} \subset W^{m-k, q}$ and the embedding is continuous, i.e. there exists a constant $c>0$ such that $\|u\|_{m-k, q} \leq c\|u\|_{m, p}, \forall u \in W^{m, p}$.

(ii) If $\frac{1}{p}-\frac{k}{n} \leq 0$, then for any $q \in[1, \infty), W^{m, p} \subset W^{m-k, q}$ and the embedding is continuous.
\end{lemma}

\begin{lemma}[Cauchy–Schwarz inequality]
\label{lemma::CS}
 We define the functions $f$ and $g$. Suppose $f\in L^2{(\Omega)}$ and $g\in L^2{\left ( \Omega \right ) }$, then 
\begin{align}
    \left(\int f(\boldsymbol{x}) g(\boldsymbol{x}) d x\right)^2 \leq \int f^2(\boldsymbol{x}) d x \int g^2(\boldsymbol{x}) d x.
\end{align}
\end{lemma}
\begin{lemma}[Grönwall inequality]
\label{lemma::Gronwall}

(i) Let $\xi(x)$ be a nonnegative summable function on $\left [ 0,T \right ] $ which satisfies for a.e. $t$ the integral inequality
\begin{align}
    \xi(t) \leq C_1 \int_0^t \xi(s) \mathrm{d} s+C_2
\end{align}

for constants $C_1 $, $C_2 \ge 0$. Then 
\begin{align}
    \xi(t) \leq C_2\left(1+C_1 t e^{C_1 t}\right)
\end{align}
for a.e. $0\leq t\leq T$.

(ii) In particular, if
\begin{align*}
    \xi(t)\leq C_1 \int_0^t \xi(s)ds
\end{align*}
for a.e. $0\leq t\leq T$, then 
\begin{align*}
    \xi(t) =0 \quad a.e.
\end{align*}
\end{lemma}

\begin{lemma}[Multiplicative trace inequality,\cite{Ryck2022Estimates}]
\label{lemma::trace}

Let $d \ge 2$, $\Omega \subset \mathbb{R}^d$ be a Lipschitz domain, and let $ \gamma_0 : H^1(\Omega) \to L^2(\partial \Omega)$ : $u \longmapsto  u|_{\partial \Omega} $ be the trace operator. Denote by $h_\Omega $ the diameter of $\Omega$ and by $\rho_\Omega$ the radius of the largest d-dimensional ball that can be inscribed into $\Omega$. Then it holds that
\begin{align}
    \left\|\gamma_0 u\right\|_{L^2(\partial \Omega)} \leq C_{h_{\Omega}, d, \rho_{\Omega}}\|u\|_{H^1(\Omega)},
\end{align}
where $C_{h_{\Omega}, d, \rho_{\Omega}}=\sqrt{2 \max \left\{2 h_{\Omega}, d\right\} / \rho_{\Omega}} $
\end{lemma}

\begin{lemma}[\cite{Ryck2022Estimates}]
\label{lemma::u_NetworkError}
Let $d,n,L,W\in \mathbb{N}$ and let $u_\theta:\mathbb{R}^{d+1}\to \mathbb{R}^{d+1}$ be a neural  network with $\theta\in \Theta 
$ for $L\geq 2,R,W\geq 1
$, cf. Definition \ref{def::net}. Assume that $\|\sigma\|_{C^n}\geq 1$. Then it holds for $1\leq j\leq d+1 $ that 
\begin{align}
    \|\left(u_\theta\right)_j\|_{C^n}\leq 16^L\left(d+1\right)^{2n}\left(e^2n^4W^3R^n\|\sigma\|_{C^n}\right)^{nL}.
\end{align}
\end{lemma}
\begin{lemma}[{  Monte-Carlo (LHS) sampling}]
\label{lemma::quadrules}

Given $\varLambda \subset \mathbb{R^d}$ and $f \in L^1(\varLambda)$, we want to approximate $\int_\varLambda  f(x) d\boldsymbol{x}$
. {  Monte-Carlo (Latin hypercube) sampling provides such an approximation by drawing $M$ sample points $\boldsymbol{x}_n \in \varLambda$ ($1\leq n\leq M$) via Latin hypercube sampling and using the sample mean with appropriate weights $\lambda^{(n)}$ ($1\leq n\leq M$), i.e.,}
\begin{align}
    \frac{1}{M}\sum_{n=1}^M \lambda^{(n)} f(\boldsymbol{x}_n) \approx \int_\varLambda f(\boldsymbol{x}) d\boldsymbol{x}.
\end{align}
{  The approximation accuracy is given by}
\begin{align}
\label{equ::quadrules}
\mathcal{Q}^\varLambda_M\left[f\right]:=\frac{1}{M}\sum_{n=1}^M \lambda^{(n)}f(\boldsymbol{x}_n), \quad {  |\int_\varLambda f(\boldsymbol{x})d\boldsymbol{x}-\mathcal{Q}_M^\varLambda\left[f\right]|\leq C_f M^{-\frac{1}{2}},}
\end{align}
{  where $C_f$ depends on the variance of $f$ and the domain (e.g. $C_f \lesssim \|f\|_{L^2}$ under standard assumptions).}
\end{lemma}

\section{Proofs of Main Lemma and Theorem}
\label{Proofs}
\textbf{Proof of Lemma \ref{acle2}.}
For the base case, since $u_t - \epsilon^2 \nabla^2 u = -f(u)$ and $\psi \in H^r(D)$, standard parabolic theory gives $u \in C([0,T]; H^r(D))$. The time derivative regularity $u_t \in C([0,T]; H^{r-2}(D))$ follows from the equation structure.

For the inductive step with $k \geq 1$: Assume $u \in H^{r+2(k-1)}(D \times [0,T])$. Since $r > \frac{d}{2} + 2k \geq \frac{d}{2} + 2$, by Sobolev embedding, $u$ is bounded, so $f(u) \in H^{r+2(k-1)}(D \times [0,T])$ by the chain rule and smoothness assumptions on $f$. 

Differentiating the Allen-Cahn equation $2k-2$ times in space (using the fact that spatial and time derivatives commute for this equation), we get that $\nabla^{2k-2} u_t - \epsilon^2 \nabla^{2k} u = -\nabla^{2k-2} f(u)$. Since the right-hand side is in $H^{r}(D \times [0,T])$ and the initial data $\nabla^{2k-2} \psi \in H^{r+2}(D)$, parabolic regularity theory gives $\nabla^{2k-2} u \in H^{r+2}(D \times [0,T])$, which implies $u \in H^{r+2k}(D \times [0,T])$.

The condition $r > \frac{d}{2} + 2k$ ensures that $H^{r+2k-2}(D) \hookrightarrow C^{2k-2}(D)$, which is needed for the nonlinear term estimates.

\textbf{Proof of Lemma \ref{chle2}.} 
We proceed by induction on $k$. The Cahn-Hilliard equation can be written as:
$$u_t = \nabla^2(\kappa f(u) - \alpha \kappa \nabla^2 u) = \kappa \nabla^2 f(u) - \alpha \kappa \nabla^4 u$$

\textbf{Base Case ($k = 0$):} This follows from Lemma \ref{chle21} (basic regularity for fourth-order parabolic equations).

\textbf{Base Case ($k = 1$):} We have $\psi \in H^{r+4}(D)$ with $r > \frac{d}{2} + 4$, and $f \in C^1(\mathbb{R})$ with $|f'(s)| \leq C(1 + |s|^{p-1})$.

From the basic regularity result, there exists $T_0 > 0$ such that $u \in C([0,T_0]; H^r(D)) \cap C^1([0,T_0]; H^{r-4}(D))$.

Since $r > \frac{d}{2} + 4 > \frac{d}{2}$, we have $H^r(D) \hookrightarrow L^\infty(D)$, so there exists $M > 0$ such that $\|u(t)\|_{L^\infty} \leq M$ for $t \in [0,T_0]$ (possibly reducing $T_0$).

For $f(u) \in C([0,T_0]; H^r(D))$, we use the chain rule. Since $u \in C([0,T_0]; H^r(D))$ and $f \in C^1(\mathbb{R})$ with controlled growth:
\begin{align}
\|f(u(t))\|_{H^r} &\leq C\|f'(u(t))\|_{L^\infty} \|u(t)\|_{H^r} + C\|f(u(t))\|_{L^\infty} \\
&\leq C(1 + M^{p-1})\|u(t)\|_{H^r} + C(1 + M^p) \\
&\leq C'
\end{align}

Therefore, $\nabla^2 f(u) \in C([0,T_0]; H^{r-2}(D))$. From the equation $u_t = \kappa \nabla^2 f(u) - \alpha \kappa \nabla^4 u$: $\nabla^4 u \in C([0,T_0]; H^{r-4}(D))$ since $u \in C([0,T_0]; H^r(D))$ and $\nabla^2 f(u) \in C([0,T_0]; H^{r-2}(D))$. Thus, $u_t \in C([0,T_0]; H^{r-4}(D))$.

Taking the time derivative of the equation, we have:
$$u_{tt} = \kappa \nabla^2(f'(u) u_t) - \alpha \kappa \nabla^4 u_t.$$

Since $u_t \in C([0,T_0]; H^{r-4}(D))$ and $r-4 \geq \frac{d}{2}$ (since $r > \frac{d}{2} + 4$), we have $u_t$ bounded in $L^\infty$. Using the product rule and chain rule, we have:
$$\nabla^2(f'(u) u_t) = f''(u)|\nabla u|^2 + f'(u)\nabla^2 u_t + 2\nabla f'(u) \cdot \nabla u_t.$$

All terms on the right are in $C([0,T_0]; H^{r-4}(D))$ by our regularity assumptions and Sobolev embeddings. Therefore, $u_{tt} \in C([0,T_0]; H^{r-8}(D))$.

For any multi-index $\alpha$ with $|\alpha| = 4$, differentiate the original equation, we get:
$$(\partial^\alpha u)_t = \kappa \partial^\alpha \nabla^2 f(u) - \alpha \kappa \partial^\alpha \nabla^4 u.$$

Using Leibniz rule for $\partial^\alpha \nabla^2 f(u) = \partial^\alpha \nabla^2 f(u)$, we get:
$$\partial^\alpha f(u) = f'(u) \partial^\alpha u + \sum_{0 < \beta \leq \alpha} \binom{\alpha}{\beta} (\partial^{\alpha-\beta} f'(u)) \partial^\beta u.$$

Since $|\alpha-\beta| \leq 4$ and $|\beta| \leq 4$, all derivatives of $u$ up to order 4 are in appropriate Sobolev spaces from our assumptions. The derivatives $\partial^{\alpha-\beta} f'(u)$ involve compositions of $f'$ with $u$ and its derivatives, which can be controlled using $f' \in C^0$ with polynomial growth, chain rule estimates in Sobolev spaces and the fact that $r > \frac{d}{2} + 4$ ensures all necessary embeddings

Since $\psi \in H^{r+4}(D)$, we have $\partial^\alpha \psi \in H^r(D)$. By fourth-order parabolic regularity theory applied to:
$$(\partial^\alpha u)_t + \alpha \kappa \nabla^4(\partial^\alpha u) = \kappa \partial^\alpha \nabla^2 f(u)$$
with initial data $\partial^\alpha u(0) = \partial^\alpha \psi \in H^r(D)$ and right-hand side in $C([0,T_0]; H^{r-2}(D))$, we obtain:
$$\partial^\alpha u \in C([0,T_0]; H^r(D)) \cap C^1([0,T_0]; H^{r-4}(D)).$$

Taking the union over all $|\alpha| = 4$, we get $u \in H^{r+4}(D \times [0,T_0])$ and $u_t \in H^r(D \times [0,T_0])$. Assume the result holds for $k-1$, i.e., there exists $T_{k-1} > 0$ such that:
$$u \in H^{r+4(k-1)}(D \times [0,T_{k-1}]) \quad \text{and} \quad u_t \in H^{r+4(k-1)-4}(D \times [0,T_{k-1}])$$

We want to show that for some $T_k \leq T_{k-1}$:
$$u \in H^{r+4k}(D \times [0,T_k]) \quad \text{and} \quad u_t \in H^{r+4k-4}(D \times [0,T_k])$$

Since $r > \frac{d}{2} + 4k \geq \frac{d}{2} + 4(k-1) + 4$, we have $r + 4(k-1) > \frac{d}{2}$, which ensuring:
$$H^{r+4(k-1)}(D) \hookrightarrow L^\infty(D).$$

Therefore, $\|u(t)\|_{L^\infty} \leq C\|u(t)\|_{H^{r+4(k-1)}} \leq C$ for $t \in [0,T_{k-1}]$. For any multi-index $\beta$ with $|\beta| = 4(k-1)$, we need to show that $\partial^\beta f(u)$ has the required regularity.

Using Faà di Bruno's formula:
$$\partial^\beta f(u) = \sum_{j=1}^{|\beta|} f^{(j)}(u) \cdot P_j(\partial^{\gamma_1} u, \ldots, \partial^{\gamma_j} u),$$
where $P_j$ are universal polynomials and $|\gamma_i| \leq |\beta| = 4(k-1)$.

Since $f \in C^k(\mathbb{R})$ and $|\beta| = 4(k-1) \leq 4k-4 < k$ (for $k \geq 1$), all required derivatives $f^{(j)}$ are bounded on the range of $u$. By the inductive hypothesis, all $\partial^{\gamma_i} u$ with $|\gamma_i| \leq 4(k-1)$ belong to appropriate Sobolev spaces with regularity index at least $r$. Since $r > \frac{d}{2}$, all polynomial combinations $P_j$ remain bounded in the required Sobolev spaces.

Therefore, $\partial^\beta f(u) \in H^r(D \times [0,T_{k-1}])$ for all $|\beta| = 4(k-1)$.

Differentiating the Cahn-Hilliard equation $4(k-1)$ times in space, we have:
$$\partial^\beta u_t = \kappa \partial^\beta \nabla^2 f(u) - \alpha \kappa \partial^\beta \nabla^4 u.$$

Since $\partial^\beta \psi \in H^{r+4k-4(k-1)} = H^{r+4}(D)$ and the right-hand side belongs to $H^r(D \times [0,T_{k-1}])$ after applying $\nabla^2$ to the first term, fourth-order parabolic regularity gives:
$$\partial^\beta u \in C([0,T_k]; H^{r+4}(D)) \cap C^1([0,T_k]; H^r(D)).$$

Taking the union over all $|\beta| = 4(k-1)$, we have:
$$u \in H^{r+4+4(k-1)}(D \times [0,T_k]) = H^{r+4k}(D \times [0,T_k]),$$
and
$$u_t \in H^{r+4(k-1)}(D \times [0,T_k]) = H^{r+4k-4}(D \times [0,T_k]).$$

This completes the inductive step and the proof.

\section{Detailed Experimental Settings}
\label{sec:Detailed Experimental Settings}
{ To ensure transparency and facilitate replicability, the detailed hyperparameters for the proposed Residuals-RAE method and the baseline methods (SA-PINNs, PINNs-WE, and Standard PINN) are summarized in Table \ref{tab:hyperparameters_detailed}.

\begin{table}[H]
\centering
\caption{Detailed hyperparameter settings for all methods used in the comparative experiments. Note that hyperparameters for baseline methods (SA-PINNs and PINNs-WE) were selected based on the optimal settings reported in their original literature to ensure fair comparison.}
\label{tab:hyperparameters_detailed}
\resizebox{\textwidth}{!}{%
\begin{tabular}{lcccc}
\toprule
\textbf{Parameter / Method} & \textbf{Residuals-RAE (Ours)} & \textbf{SA-PINNs} & \textbf{PINNs-WE} & \textbf{Standard PINN} \\
\midrule
\textbf{Network Architecture} & 2 layers, 128/256 neurons & 2 layers, 128/256 neurons & 2 layers, 128/256 neurons & 2 layers, 128/256 neurons \\
\textbf{Activation Function} & Tanh & Tanh & Tanh & Tanh \\
\textbf{Optimizer} & Adam + L-BFGS & Adam + L-BFGS & Adam + L-BFGS & Adam + L-BFGS \\
\textbf{Learning Rate (LR)} & $1.0 \times 10^{-3}$ & $1.0 \times 10^{-3}$ & $1.0 \times 10^{-3}$ & $1.0 \times 10^{-3}$ \\
\textbf{Weight Initialization} & Xavier & Xavier & Xavier & Xavier \\
\midrule
\textbf{Special Hyperparameters} & 
\makecell[c]{$k_{int} = 50$ \\ (neighbors for K-NN) \\ $\beta = 0.9$ (momentum)} & 
\makecell[c]{\textbf{Mask LR} $= \mathbf{10^{-2}}$ \\ (mask update rate) \\ Init Mask = 1.0} & 
\makecell[c]{Scaling $k = 1.0$ \\ \textbf{Update Freq} $= \mathbf{100}$ \\ (weighting parameter)} & 
N/A \\
\bottomrule
\end{tabular}%
}
\end{table}
Note: Network width is set to 128 for Allen-Cahn and 256 for Cahn-Hilliard experiments to accommodate problem complexity, applied consistently across all methods.}

\section{Detailed Training Algorithm (Residual-RAE)}
\label{sec:appendix_algorithm}
{ To explicitly illustrate the decoupled nature of our training strategy and ensure the reproducibility of the proposed Residuals-RAE-PINNs, we provide the complete pseudocode in Algorithm \ref{alg:residuals_rae}. Unlike conventional self-adaptive methods where weights are optimized via gradient descent, our approach computes weights \textit{a priori} using the current PDE residuals and spatial distribution information. This separation ensures that the network parameters are updated under a fixed weight landscape within each iteration, enhancing numerical stability.

\begin{algorithm}[H]
    \caption{Training Procedure for Residuals-RAE-PINNs}
    \label{alg:residuals_rae}
    
    \SetKwInOut{Input}{Input}
    \SetKwInOut{Output}{Output}

    \Input{
        Training datasets $\mathcal{D}_{int}$; \\
        Hyperparameters: $K_{max}, \eta, \beta, k_{int}$; \\
        Model: Shallow Neural Network $u_{\theta}$ ($L=2$).
    }
    \Output{Optimized network parameters $\theta^*$.}
    
    \BlankLine
    \textbf{Initialization:} Initialize { $\theta^{0}$} and weights { $\boldsymbol{\lambda}_{int}^{0} \leftarrow \mathbf{1}$}\;
    \BlankLine

    {  \tcc{Notation: iteration index $k$ and $\theta^k$, $\boldsymbol{\lambda}_{int}^k$ match Section~2 and the weight-convergence proof.}}
    \For{$k = 1$ \KwTo $K_{max}$}{
        \tcc{Phase 1: Pre-training Weight Update (Decoupled)}
        Compute interior PDE residuals using { $\theta^{k-1}$}: \\
        $R_{int}^{(i)} \leftarrow \left| \mathcal{N}[u_{\theta^{k-1}}](\boldsymbol{x}_{int}^{(i)}, t_{int}^{(i)}) \right|$ for all $i \in \mathcal{D}_{int}$\;
        
        Compute normalized simple weights $w_{int}$ via Eq. \eqref{eq::simpleRAE}: \\
        $w_{int}^{(i)} \leftarrow \frac{R_{int}^{(i)}}{\sum_{j} R_{int}^{(j)}} \cdot M_{int}$\;
        
        Apply K-Nearest Algorithm for spatial smoothing: \\
        Find set $\boldsymbol{S}_{k_{int}}^{(i)}$ (nearest neighbors) and compute $\lambda_{Knear,r}^{(i)}$ for all $i$, then form the vector $\boldsymbol{\lambda}_{Knear,r}^k$ via Eq. \eqref{eq::lknear}\;
        
        Update adaptive weights with momentum: \\
        { $\boldsymbol{\lambda}_{int}^{k} \leftarrow \beta \cdot \boldsymbol{\lambda}_{Knear,r}^k + (1-\beta) \cdot \boldsymbol{\lambda}_{int}^{k-1}$}\;
        
        \textbf{Freeze} { $\boldsymbol{\lambda}_{int}^{k}$} (treat as constants)\;
        
        \BlankLine
        \tcc{Phase 2: Network Optimization}
        Construct modified loss $\mathcal{J}$ via Eq. \eqref{eq::modifiedloss} using { $\boldsymbol{\lambda}_{int}^{k}$}\;
        
        Compute Gradients: $\mathbf{g} \leftarrow \nabla_{\theta} \mathcal{J}$\;
        Update Parameters: { $\theta^{k} \leftarrow \theta^{k-1} - \eta \cdot \mathbf{g}$}\;
    }
    \Return{{ $u_{\theta^{K_{max}}}$}}
\end{algorithm}
\paragraph{Implementation Details} 
In our practical implementation, the K-Nearest Neighbors (K-NN) search step (Line 5 in Algorithm \ref{alg:residuals_rae}) is efficiently computed using the \texttt{scikit-learn} library (or \texttt{torch.cdist} depending on your code) to minimize computational overhead. The weights $\boldsymbol{\lambda}_{int}$ are normalized at each step to maintain a consistent scale for the loss function. For the specific values of hyperparameters used in our experiments (such as $k_{int}$ and $\beta$), please refer to the detailed summary in \textbf{Table \ref{tab:hyperparameters_detailed} (Appendix \ref{sec:Detailed Experimental Settings})}.}
\end{appendices}

\renewcommand\refname{\vskip -1cm}

\renewcommand\refname{\vskip -1cm}
\newpage
\section*{References}
\bibliography{references}

\end{document}